\documentclass[12pt]{amsart}

\usepackage{amsmath, amssymb, fullpage}
\usepackage[colorlinks=true,citecolor=black,linkcolor=black,urlcolor=blue]{hyperref}

\usepackage{silence}
\usepackage{thmtools}
\usepackage{thm-restate}

\usepackage{comment}

\usepackage{tikz}
\tikzstyle{vertex}=[circle, draw=black, inner sep=0pt, minimum size=5pt]

\newcommand{\vertex}{\node[vertex]}
\tikzstyle{vtx}=[circle, draw, inner sep=0pt, minimum size=8pt]

\newtheorem{theorem}{Theorem}[section]
\newtheorem{proposition}[theorem]{Proposition}
\newtheorem{corollary}[theorem]{Corollary}
\newtheorem{lemma}[theorem]{Lemma}
\theoremstyle{definition}

\newtheorem{definition}[theorem]{Definition}
\newtheorem{example}[theorem]{Example}
\newtheorem{remark}[theorem]{Remark}

\newtheorem{question}[theorem]{Question}

\definecolor{applegreen}{rgb}{0.55,0.71,0.0}

\definecolor{darkcandyapplered}{rgb}{0.64, 0.0, 0.0}

\renewcommand{\hat}{\widehat}
\def\A{\mathcal{A}}
\def\E{\mathcal{E}}
\def\F{\mathcal{F}}

\def\P{\mathcal{P}}
\def\C{\ensuremath{S}}
\def\H{\ensuremath{C}}
\def\des{\ensuremath{\mathsf{des}}}
\def\skipstat{\ensuremath{\mathsf{skip}}}
\def\Skip{\ensuremath{\mathsf{SKIP}}}
\def\tail{\ensuremath{\mathsf{tail}}}
\def\head{\ensuremath{\mathsf{head}}}
\def\act{\ensuremath{\mathsf{ACT}}}
\def\inact{\ensuremath{\mathsf{INACT}}}

\def\indeg{\ensuremath{\mathrm{indeg}}}
\def\zero{\ensuremath{\boldsymbol{0}}}
\def\one{\ensuremath{\boldsymbol{1}}}
\def\Cat{\ensuremath{\mathsf{Cat}}}

\def\RR{\mathbb{R}}
\def\ZZ{\mathbb{Z}}

\def\aa{\mathbf{a}}
\def\bb{\mathbf{b}}
\def\ee{\mathbf{e}}

\def\kk{\mathbf{k}}

\def\ss{\mathbf{s}}
\def\tt{\mathbf{t}}

\def\vv{\mathbf{v}}

\def\zz{\mathbf{z}}

\DeclareMathOperator{\vol}{vol}
\DeclareMathOperator{\stat}{\ensuremath{\mathsf{stat}}}

\title{Column convex matrices,\\ $G$-cyclic orders, and flow polytopes}

\author[Gonz\'alez D'Le\'on]{Rafael S. Gonz\'alez D'Le\'on }
\address[R.\ S.\ Gonz\'alez D'Le\'on]{Escuela de Ciencias Exactas e Ingenier\'ia\\Universidad Sergio Arboleda\\Bogot\'a\\Colombia} 
\email{rafael.gonzalezl@usa.edu.co}
\urladdr{\url{http://dleon.combinatoria.co}}

\author[Hanusa]{Christopher R.\ H.\ Hanusa}
\address[C.\ R.\ H.\ Hanusa]{Department of Mathematics \\ Queens College (CUNY) \\ 65-30 Kissena Blvd. \\ Flushing, NY 11367\\ United States}
\email{\href{mailto:chanusa@qc.cuny.edu}{\texttt{chanusa@qc.cuny.edu}}}
\urladdr{\url{http://qc.edu/~chanusa/}}

\author[Morales]{Alejandro H. Morales}
\address[A.\ H.\ Morales]{Department of Mathematics and Statistics, University of Massachusetts, Amherst, MA, 01003, United States} 
\email{ahmorales@math.umass.edu}
\urladdr{\url{http://people.math.umass.edu/~ahmorales/}}

\author[Yip]{Martha Yip}
\address[M.\ Yip]{Department of Mathematics, University of Kentucky, 719 Patterson Office Tower, Lexington, KY 40506-0027, United States}
\email{martha.yip@uky.edu}
\urladdr{\url{http://www.ms.uky.edu/~myip/}}


\makeatletter
\@namedef{subjclassname@2020}{%
  \textup{2020} Mathematics Subject Classification}
\makeatother


\begin{document}
\maketitle

\begin{abstract}
We study polytopes defined by inequalities of the form $\sum_{i\in I} z_{i}\leq 1$ for $I\subseteq [d]$ and nonnegative $z_i$ where the inequalities can be reordered into a matrix inequality involving a column-convex $\{0,1\}$-matrix. These generalize polytopes studied by Stanley, and the consecutive coordinate polytopes of Ayyer, Josuat-Verg\`es, and Ramassamy. We prove an integral equivalence between these polytopes and flow polytopes of directed acyclic graphs $G$ with a Hamiltonian path, which we call spinal graphs. 
We show that the volume of these flow polytopes is the number of extensions of a set of partial cyclic orders defined by the graph $G$. As a special case we recover results on volumes of consecutive coordinate polytopes.

We study the combinatorics of $k$-Euler numbers, which are generalizations of the classical Euler numbers, and which arise as volumes of flow polytopes of a special family of spinal graphs. We show that their refinements, Ramassamy's $k$-Entringer numbers, can be realized as values of a Kostant partition function, satisfy a family of generalized boustrophedon recurrences, and are log concave along root directions.

Finally, via our main integral equivalence and the known formula for the $h^*$-polynomial of consecutive coordinate polytopes, we give a combinatorial formula for the $h^*$-polynomial of flow polytopes of non-nested spinal graphs.  
For spinal graphs in general, we present a conjecture on upper and lower bounds for their $h^*$-polynomial.
\end{abstract}


\section{Introduction} \label{section:introduction}

\subsection{A motivating question of Stanley}
Stanley proposes the question \cite[Exercise 4.56\;(d)]{Stanley2012} of finding a formula for the volume of the polytope $\mathcal{C}_{d,k}$ in $\RR^d$ defined by the inequalities $z_i \ge 0$ for all $i=1,\dots, d$, and 
\begin{equation}\label{eqn.Cdk}
    z_{i+1}+z_{i+2}+\cdots + z_{i+k}\le 1, 
\end{equation}
for all $i=0,\ldots, d-k$. 

A polytope of this form is an \emph{integral polytope}: one whose vertices all lie in $\ZZ^d$. It is well-known that the Euclidean volume of a $d$-dimensional integral polytope $\mathcal{P}$ is of the form $\vol\mathcal{P}/d!$ where $\vol \mathcal{P} \in \ZZ_{\geq0}$. The quantity $\vol\mathcal{P}$ is commonly known as the \emph{normalized volume} of $\mathcal{P}$, and this is the notion of volume we use in this article.

Stanley proves that $\vol \mathcal{C}_{d,2}$ is the $d$-th Euler number $E_d$ by showing that $\mathcal{C}_{d,2}$ is the chain polytope of the zigzag poset on $d$ elements \cite[Exercise 4.56 (b) and (c)]{Stanley2012}, whose volume is given by the number of linear extensions of the poset \cite{Stanley86}. 
For general~$k$, Stanley gives a set of difference equations~\cite[Exercise 4.56\;(d)]{Stanley2012} that can be used to recursively compute $\vol\mathcal{C}_{d,k}$ but leaves open the problem of finding a direct combinatorial formula for this volume.

Ayyer, Josuat-Verg\`es, and Ramassamy~\cite{AyyerJosuatvergesRamassamy2018} give a beautiful answer to Stanley's question by showing that the volume of $\mathcal{C}_{d,k}$ is the number of total cyclic orders that are extensions of a collection of partial cyclic orders, which has the same spirit as the result of Stanley in the case $k=2$. 
Furthermore, they apply their method to compute the volumes of polytopes belonging to the larger family of \emph{consecutive coordinate polytopes}, in which the defining inequalities are of the form
\begin{equation}\label{eqn.Cdk2}
    z_{i}+z_{i+1}+\cdots + z_{i'}\le 1
\end{equation}
for integers $i<i'$. This motivates studying polytopes whose defining inequalities are not comprised of consecutive coordinates.  

We define the polytope $\mathcal{B}_{\C}$ for a collection $\C$ of subsets of $[d]:=\{1,\ldots,d\}$ to be
\begin{align*}
  \mathcal{B}_{\C}=\bigg\{(z_1,\ldots,z_d) \in \mathbb{R}_{\geq0}^d ~\bigg|~ 
\sum_{i\in I}z_{i} \leq 1  \hbox{ for } I\in \C \bigg\}. 
\end{align*}
Consecutive coordinate polytopes are examples of such polytopes in the special case when $\C$ is a collection of \emph{intervals} $[i,i']:=\{i,i+1,\hdots,i'\}$.

It is known that the polytopes $\mathcal{C}_{d,k}$ are not chain polytopes of posets; see \cite[Remark 2.6]{AyyerJosuatvergesRamassamy2018}. Instead, we are able to show that a subfamily of polytopes of the form $\mathcal{B}_\C$ (which includes $\mathcal{C}_{d,k}$ and all consecutive coordinate polytopes) are integrally equivalent to flow polytopes of a special type of graph. 


\subsection{Convex \texorpdfstring{$\{0,1\}$}{}-matrices and flow polytopes}
A collection of $n$ inequalities of the form \eqref{eqn.Cdk} and \eqref{eqn.Cdk2} can be written as a matrix inequality of the form $M\zz\leq \bb$ for an $n\times d$ $\{0,1\}$-matrix $M$ and a nonnegative vector $\bb\in \mathbb{R}_{\geq0}^n$. We define polytopes $\mathcal{B}_{M,\bb}$ as those of the form 
$$\mathcal{B}_{M,\bb} = \Big\{\zz \in \mathbb{R}_{\geq0}^d ~\Big|~  M\mathbf{z}\leq \mathbf{b} \Big\},$$
where $M$ has no columns that are identically zero, and develop their theory in Section~\ref{sec:integralequivalence}.

Every consecutive coordinate polytope is a polytope of the form $\mathcal{B}_{M,\one}$ where $\one=(1,\hdots,1)$ and $M$ is row convex---that is, the non-zero entries in every row are contiguous. (Similarly, column-convex matrices are those in which the non-zero entries in every column are contiguous.) In fact, in Lemma~\ref{lemma:row_convex_is_doubly_convex} we show through matrix operations that the class of consecutive coordinate polytopes is the same as the class of polytopes $\mathcal{B}_{M,\one}$ for matrices $M$ that are simultaneously row and column convex (sometimes called \emph{interval matrices} \cite[Chapter~19]{schrijver1998theory}).

We prove that in the more general setting where $M$ is a column convex matrix, the polytope $\mathcal{B}_{M,\bb}$ is integrally equivalent to a flow polytope of a directed acyclic graph $G$ with a Hamiltonian path, which we call a {\em spinal graph}.

\begin{restatable*}[]{theorem}{integralequivalencethm}
\label{thm:integral_equivalence}
Let $\bb\in \ZZ^n$, $M$ an $n\times d$ column-convex matrix, and $G$ its associated graph. The polytope $\mathcal{B}_{M,\bb}$ is integrally equivalent to the flow polytope $\mathcal{F}_{G}(\mathbf{a})$ where \[\mathbf{a}=(b_1,b_2-b_1,\dots,b_{n}-b_{n-1},-b_n).\]
\end{restatable*}

When $\bb=\one$, then $\aa=\ee_1-\ee_{n+1}$, and we write $\mathcal{F}_G$ for the flow polytope $\mathcal{F}_{G}(\mathbf{a})$. Section~\ref{sec:integralequivalence} concludes by discussing operations on graphs that preserve integral equivalence of the associated flow polytopes.

Theorem \ref{thm:integral_equivalence} has the benefit of transporting all the machinery of triangulations of flow polytopes \cite{DKK,MMS,KMS19} and  enumeration of lattice points and volumes of flow polytopes \cite{BV,MM19,BGHHKMY} to the family of polytopes associated to column-convex matrices, that in turn contains the family of consecutive coordinate polytopes.


\subsection{The combinatorics of flow polytope volumes}
\smallskip
With our new integral equivalence, it is natural to develop combinatorial tools to calculate the volume of the polytopes, which we do in Sections~\ref{sec:volume}--\ref{sec:distance}.

\smallskip
Techniques to compute the volume and lattice points of flow polytopes have been extensively studied in recent literature.  Baldoni and Vergne \cite{BV} give a set of Lidskii formulas to calculate the volume and lattices points of flow polytopes. Postnikov and Stanley (unpublished) describe a triangulation that can be used to provide a different proof of the Lidskii formulas when $\aa=\ee_1-\ee_{n+1}$. Meszaros and Morales \cite{MM19} extend Postnikov and Stanley's triangulation to any netflow vector $\aa$.  In \cite{BGHHKMY} Benedetti et al.\ introduce \emph{gravity diagrams} as a family of combinatorial objects whose enumeration can be used as a tool to calculate the volume. Section~\ref{sec:volume} presents the Lidskii volume formulas that are relevant to our present work. 

\smallskip
In Section~\ref{sec:G-cyclic} we discuss total cyclic orders and we introduce a pair of new combinatorial objects called upper and lower $G$-cyclic orders for any spinal graph $G$. We use the Lidskii formulas to prove that the enumeration of either gives the volume of the flow polytope $\mathcal{F}_G$.

\begin{restatable*}[]{theorem}{volumethm}
\label{thm.volume}
For a spinal graph $G$, the volume of the flow polytope $\F_{G}$ is the number of upper (or lower) $G$-cyclic orders. In other words,
$$\vol \F_{G} = A^{\uparrow}_{G}=A^{\downarrow}_G.$$
\end{restatable*}

A significant proportion of flow polytopes that have been studied in the literature are cases of flow polytopes of spinal graphs, so Theorem \ref{thm.volume} has a broad scope of application. Some examples of these graphs are  Pitman-Stanley graphs \cite{BV}, graphs used in M\'esz\'aros' product formulas for volumes of flow polytopes \cite{Mes15}, Corteel-Kim-M\'esz\'aros graphs~\cite{CKM}, zigzag graphs~\cite{BGHHKMY}, and caracol graphs and their multigraph generalizations~\cite{BGHHKMY, MM19, BGMY}.

In Section~\ref{sec:consec} we show that when restricting to consecutive coordinate polytopes, of which there are Catalan many by Proposition~\ref{prop:number nonnesting}, the corresponding spinal graphs satisfy a non-nested condition. For those graphs Proposition~\ref{prop.ofmaintheorem} shows that the upper and lower $G$-cyclic orders coincide. Furthermore, Proposition~\ref{prop:Gcompatibility_is_extension} shows that these $G$-cyclic orders are exactly the same as the total cyclic extensions of partial cyclic orders of  \cite{AyyerJosuatvergesRamassamy2018}.  In this sense, Theorem~\ref{thm.volume} can be seen as a generalization of Ayyer, Josuat-Verg\`es, and Ramassamy's result that uses techniques of flow polytopes instead of using polytope triangulations and a transfer map. 


\subsection{Distance graphs, \texorpdfstring{$k$}{}-Euler numbers, and \texorpdfstring{$k$}{}-Entringer numbers}

Section~\ref{sec:distance} applies our work to the family of \emph{distance graphs} $G(k,d+k)$, which have vertex set $[d+k]$ and the edges of the form $(i,i+1)$  and $(i,i+k)$.  These graphs generalize the \emph{zigzag graphs} $G(2,d+2)$ and correspond to the consecutive coordinate polytopes $\mathcal{C}_{d,k}$ studied in \cite{AyyerJosuatvergesRamassamy2018}.

In Proposition~\ref{prop: vertices G(k,d+k)} we prove that the number of vertices of $\mathcal{F}_{G(k,d+k)}$ are given by a generalization of Fibonacci numbers. The volume $\vol \mathcal{F}_{G(k,d+k)}$ can be seen as a $k$-generalization of the \emph{Euler numbers} since $\vol \mathcal{F}_{G(2,d+2)} = E_{d}$. Just as \emph{Entringer numbers} refine Euler numbers, Ramassamy \cite{Ramassamy2017} and Ayyer, Josuat-Verg\`es, and Ramassamy \cite{AyyerJosuatvergesRamassamy2018}, define \emph{$k$-Entringer numbers} $E_{(s_1,\ldots,s_k)}$ that refine $\vol \mathcal{F}_{G(k,d+k)}$. 
We realize $k$-Entringer numbers as the number of certain integer flows on $G(k,d+k)$ and as evaluations of Kostant partition functions.
By exploiting the recursive nature of distance graphs, we show in Theorem~\ref{thm:kBoustrophedon} that the $k$-Entringer numbers can be computed recursively on $k$ levels, extending the boustrophedon recursion of~\cite[Theorem 7.4]{AyyerJosuatvergesRamassamy2018}. 

Another benefit of the flow polytope perspective is the following log-concavity result for $k$-Entringer numbers. It is proved by presenting another way to realize the $k$-Entringer numbers as evaluations of Kostant partition functions which enables us to apply a result of Huh et al.~\cite[Proposition 11]{HMMStD} related to the Alexandrov--Fenchel inequalities of mixed volumes of polytopes \cite{Alexandrov, Fe2, Fe1}.

\begin{restatable*}[]{theorem}{logconcavitythm}
\label{thm.log-concavity}
Let $k\geq 2$ and $N=d-k+1\geq 0$. Given $\ss=(s_1,\ldots,s_k)\in T_N^k$, then the numbers $E_{\ss}$ are log-concave along root directions. That is, 
\[
E_{\ss}^2 \geq E_{\ss -\ee_i+\ee_j}\cdot E_{\ss +\ee_i-\ee_j}.
\]
\end{restatable*}


\subsection{The \texorpdfstring{$h^*$}{}-polynomial of \texorpdfstring{$\mathcal{F}_G$}{}}

In Section~\ref{sec:h_star} we study the $h^*$-polynomial of $\mathcal{F}_G$, whose coefficients sum to the volume of the polytope. In \cite{AyyerJosuatvergesRamassamy2018}, the authors found that the $h^*$-polynomial of the consecutive coordinate polytope is equal to the generating polynomial of the descent statistic on the family of total cyclic extensions of the partial cyclic order determined by its collection of intervals. 

When translated to flow polytopes, this result becomes Theorem~\ref{theorem:hstar_non_nested} which states that when $G$ is a non-nested graph, the $h^*$-polynomial of the flow polytope $\mathcal{F}_{G}$ is
\[
h^*_{\mathcal{F}_{G}}(z) = \sum_{ \gamma \in \mathcal{A}_G} z^{\des(\pi(\gamma))}.
\]
This reinterpreted formula is not true for a general spinal graph $G$ because the upper and lower $G$-cyclic orders are not the same. However we are able to conjecture dominance bounds on $h_{\F_G}^*(z)$ that have been verified for all simple spinal graphs $G$ with up to $7$ vertices. 
 
\begin{restatable*}[]{conjecture}{hstarconjecture}
\label{conj: hstar}
Given a spinal graph $G$ we have that 
\begin{equation*} 
P_{\mathcal{A}^\downarrow_G,\des}(z)\lhd h_{\F_G}^*(z)\lhd 
P_{\mathcal{A}^{\uparrow}_G,\des}(z).
\end{equation*}
\end{restatable*}
\noindent
In the above expression, $P_{\mathcal{A},\des}(z) =\sum_{\gamma\in \mathcal{A}}z^{\des(\pi(\gamma))}$.

Finally, Section~\ref{sec:further_work} assembles directions of future work stemming from this paper, including questions about the generating functions for $k$-Euler and $k$-Entringer numbers.


\section{Column-convex matrices and flow polytopes}\label{sec:integralequivalence}

\subsection{A key integral equivalence}
In this section we prove an integral equivalence between flow polytopes on spinal graphs and polytopes from column-convex matrices.

\begin{definition}
A directed graph $G$ on vertex set $V(G)=[n+1]$ is said to be a \emph{spinal graph} if its edge multiset $E(G)$ contains at least one edge of the form $(i,i+1)$ for each $i=1,\dots, n$ and where every other edge in $E(G)$ is directed from a smaller vertex to a larger vertex.  
The former are called \emph{slack edges} which we label $x_i$ and are said to make up the {\em spine} of the graph (and if there are multiple edges of the form $(i,i+1)$, only one of them will be considered slack), while the latter are called {\em non-slack edges} and are labeled $y_j$. For an edge $e$ we write $\tail(e)$ for the initial vertex (or tail) of $e$ and $\head(e)$ for the terminal vertex (or head) of $e$.
\end{definition}

All graphs in this article are spinal graphs and all matrices in this article are $\{0,1\}$-matrices; this will be assumed going forward. 

A matrix is said to be \emph{column convex} if it satisfies the property:
\begin{equation*}
  \text{If  $i<i'$ and entries $M_{i,j}=M_{i',j}=1$ then entry $M_{k,j}=1$ for all $i\le k \le i'$.}
\end{equation*}
A matrix is said to be \emph{row convex} if it satisfies the property:
\begin{equation*}
  \text{If  $j<j'$ and entries $M_{i,j}=M_{i,j'}=1$ then entry $M_{i,k}=1$ for all $j\le k \le j'$.}
\end{equation*}
Matrices that are both row and column convex are said to be \emph{doubly convex}. These matrices are also called \emph{interval matrices} \cite[Example 7, Chapter 19]{schrijver1998theory}.

\begin{definition}\label{def:associated}
Given a column-convex $n\times d$ 
matrix $M$ with no columns identically zero, define its \emph{associated graph} $G$ to have vertex set $[n+1]$ and $n+d$ edges of two types: An edge $x_i:(i,i+1)$ for $i=1,\dots, n$ and an edge $y_j:(i_j,i_j')$ for every column $j$ of $M$ in which the non-zero entries span from row $i_j$ to $i_j'-1$. 

Given a spinal graph $G$ on $n+1$ vertices with $d$ non-slack edges, define its \emph{associated matrix} $M$ with dimensions $n\times d$ as follows.  Draw the vertices of the graph $G$ in ascending order from left to right along a horizontal line.  Place a $1$ in position $(i,j)$ of $M$ if edge $y_j$ crosses the vertical line drawn between vertices $i$ and $i+1$ and a $0$ otherwise. 
\end{definition}

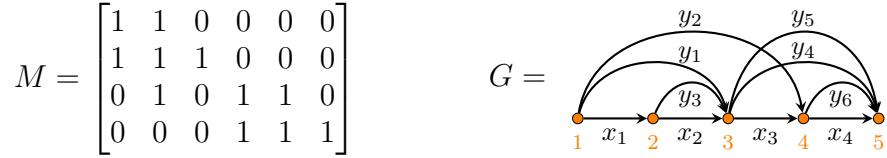
\begin{figure}[t]
\begin{center}
\begin{tikzpicture}
\begin{scope}[scale=0.7]
\node at (0,0){
$M=\begin{bmatrix}
1&1&0&0&0&0\\   
1&1&1&0&0&0\\
0&1&0&1&1&0\\
0&0&0&1&1&1
\end{bmatrix}$
};
\end{scope}
\begin{scope}[xshift=120, yshift=-15, yscale=1.2]
\node at (0.2,0.48){$G=$};
	\vertex[fill=orange, minimum size=4pt, label=below:{\tiny\textcolor{orange}{$1$}}](v1) at (1,0) {};
	\vertex[fill=orange, minimum size=4pt, label=below:{\tiny\textcolor{orange}{$2$}}](v2) at (2,0) {};
	\vertex[fill=orange, minimum size=4pt, label=below:{\tiny\textcolor{orange}{$3$}}](v3) at (3,0) {};
	\vertex[fill=orange, minimum size=4pt, label=below:{\tiny\textcolor{orange}{$4$}}](v4) at (4,0) {};
	\vertex[fill=orange, minimum size=4pt, label=below:{\tiny\textcolor{orange}{$5$}}](v5) at (5,0) {};

	\draw[-stealth, thick] (v1)--(v2);
	\draw[-stealth, thick] (v2)--(v3);
	\draw[-stealth, thick] (v3)--(v4);
	\draw[-stealth, thick] (v4)--(v5);
	\draw[-stealth, thick] (v1) .. controls (1.25,0.8) and (2.75,0.8) .. (v3);
	\draw[-stealth, thick] (v1) .. controls (1.25,1.25) and (3.75,1.25) .. (v4);	
	\draw[-stealth, thick] (v2) .. controls (2.25,0.5) and (2.75,0.5) .. (v3);
	\draw[-stealth, thick] (v3) .. controls (3.25,0.8) and (4.75,0.8) .. (v5);
	\draw[-stealth, thick] (v3) .. controls (3.25,1.25) and (4.75,1.25) .. (v5);
	\draw[-stealth, thick] (v4) .. controls (4.25,0.5) and (4.75,0.5) .. (v5);
	
	\node at (1.5, -0.2){\footnotesize$x_1$};
	\node at (2.5, -0.2){\footnotesize$x_2$};
	\node at (3.5, -0.2){\footnotesize$x_3$};
	\node at (4.5, -0.2){\footnotesize$x_4$};
	\node at (2.5, 0.7){\footnotesize$y_1$};
	\node at (2.5, 1.1){\footnotesize$y_2$};
	\node at (2.5, 0.22){\footnotesize$y_3$};
	\node at (4, 0.75){\footnotesize$y_4$};
	\node at (4, 1.1){\footnotesize$y_5$};
	\node at (4.5, 0.22){\footnotesize$y_6$};
\end{scope}
\end{tikzpicture}
\end{center}
\caption{For the column-convex matrix $M$, its associated graph $G$ has $10$ edges, including the four slack edges of the spine and the six non-slack edges from the columns of $M$.}
\label{fig:GraphFromMatrix}
\end{figure}

These two operations are inverses of each other, which we state as the following proposition.  See Figure~\ref{fig:GraphFromMatrix} for an example.

\begin{proposition}\label{proposition:bijection_col_convex_matrices_graphs}
There is a bijection between column-convex $n \times d$ matrices and spinal graphs on $n+1$ vertices with $n+d$ edges. 
\end{proposition}

\begin{definition}
For an $n\times d$ matrix $M$ and an integer vector $\bb\in \ZZ^{n}$, the polytope $\mathcal{B}_{M,\bb}$ is defined as $$\mathcal{B}_{M,\bb} = \big\{\mathbf{x}\in \mathbb{R}_{\geq0}^d  ~\big|~ 
M\mathbf{x}\leq \mathbf{b} \big\}.$$
\end{definition}

\begin{definition}\label{definition:flow_polytope}
Consider a graph $G$ on $[n+1]$ and an integer {\em net flow vector} $\aa\in \ZZ^{n+1}$ whose entries sum to zero. An {\em $\aa$-flow} on $G$ is a tuple $(f_e)_{e\in E(G)}$ of nonnegative real numbers for which flow is conserved at every internal vertex. Mathematically, 
$$\sum_{\substack{e\in E(G)\\ \tail(e)=i}}f_{e} + a_i = \sum_{\substack{e\in E(G)\\ \head(e)=i}}f_e$$
for $i=1,\ldots, n+1$. 
The {\em flow polytope} $\F_G(\aa)$ is defined as the set of $\aa$-flows on $G$.
When the flow vector $\aa$ equals $\ee_1-\ee_{n+1}$, we will write $\F_G= \F_G(\ee_1-\ee_{n+1})$ for convenience.
\end{definition}

The vertices of $\mathcal{F}_G$ are characterized as follows.

\begin{proposition}[{\cite[Corollary 3.1]{GS}}] \label{prop: char verties F_G}
Let $G$ be a graph on $[n+1]$. The vertices of $\mathcal{F}_G$ correspond to
unit flows along paths from vertex $1$ to vertex $n+1$.
\end{proposition}

Two lattice polytopes $\mathcal{P}\subset\RR^m$ and $\mathcal{Q}\subset\RR^n$ are {\em integrally equivalent} if there exists an affine transformation $\varphi:\RR^m\rightarrow\RR^n$ whose restriction to $\mathcal{P}$ is a bijection $\varphi:\mathcal{P}\rightarrow \mathcal{Q}$ that preserves the lattice.
A key observation is that the polytopes associated to convex-column matrices are integrally equivalent to flow polytopes.

\integralequivalencethm
\begin{proof}
Given an $\mathbf{a}$-flow $(f_e)_{e\in E(G)}$, 
conservation of flow at vertex $i$ is given by 
\begin{equation}\label{eq:i}
    a_i + f_{x_{i-1}} + \sum_{\substack{j \in [d]\\\head(y_j)=i}}f_{y_j} = f_{x_i} + \sum_{\substack{j \in [d]\\\tail(y_j)=i}} f_{y_j}.
\end{equation}
For the set of $n+1$ equations of the form of Equation~\eqref{eq:i}, the flow conservation equation at $n+1$ is redundant.
An equivalent set of equations describing  $\mathcal{F}_{G}(\mathbf{a})$ can be obtained by adding  Equation~\eqref{eq:i} for vertices $1$ through~$i$ for $i\in[n]$:
\begin{equation*}
    b_i=a_1+\cdots+ a_i = f_{x_i} + \sum_{\substack{j \in [d]\\\tail(y_j)=i}} f_{y_j}.
\end{equation*}
Define a map $\varphi: \mathcal{B}_{M,\bb}\rightarrow \RR^{n+d}$
that sends the point $\zz=(z_1,\ldots, z_{d}) \in \mathcal{B}_{M,\bb}$ to a flow on $G$ by 
$$
f_{x_i}= b_i-\sum_{\substack{j \in [d]\\\tail(y_j)=i}} z_j\quad\hbox{ and }\quad f_{y_j} = z_j .$$
Writing $\varphi(\zz)$ as the vector $(f_{x_1},\ldots, f_{x_n},f_{y_1},\ldots, f_{y_{d}})$, the map $\varphi$ can be written as the affine linear map
$$\varphi(\zz)=\begin{bmatrix}
-M\\
I
\end{bmatrix}\zz+ \begin{bmatrix}
\mathbf{b}\\
{\zero}
\end{bmatrix},$$
where $I$ denotes the $d \times d$ identity matrix and  $\zero=(0,\hdots,0)$. 
The map $\varphi$ is an injection whose image is $\mathcal{F}_{G}(\mathbf{a})$ and preserves the affine lattices generated by each polytope. Therefore, $\F_{G}(\mathbf{a})$ and $\mathcal{B}_{M,\bb}$ are integrally equivalent.
\end{proof}

\begin{corollary}
Let $M$ be an $n \times d$ column-convex matrix with associated graph $G$ and let $\bb\in \ZZ^n$. 
We have that
$$\vol \mathcal{B}_{M,\bb}=\vol \mathcal{F}_{G}(\mathbf{a})$$
where $\mathbf{a}=(b_1,b_2-b_1,\dots,b_{n}-b_{n-1},-b_n)$.
\end{corollary}


\subsection{Integrally equivalent flow polytopes through matrix operations} \label{sec: operations and redundancies}

The integral equivalence in Theorem~\ref{thm:integral_equivalence} allows us to establish the integral equivalence of a family of flow polytopes of different graphs related by transformations on their associated matrices. 
While a reordering of the columns of the matrix only corresponds to a relabeling of the non-slack edges, a reordering of the rows of the matrix is a more fruitful transformation. 

\begin{proposition}
Let $G$ and $G'$ be two graphs whose associated matrices $M$ and $M'$ differ by a reordering of their rows. If $\bb$ and $\bb'$ in $\mathbb{Z}^n$ are related by the same reordering and $\aa$ and $\aa'$ are defined by $\aa=(b_1,b_2-b_1,\hdots,b_n-b_{n-1},-b_n)$ and $\aa'=(b_1',b_2'-b_1',\hdots,b_n'-b_{n-1}',-b_n')$ then $\mathcal{F}_{G}(\mathbf{a})$ and $\mathcal{F}_{G'}(\mathbf{a}')$ are integrally equivalent.
\end{proposition}
\begin{proof}
A reordering of the rows of $M$ (and the corresponding entries of $\bb$) does not change the polytope $\mathcal{B}_{M,\bb}$ because the defining inequalities remain the same. The proposition follows by Theorem~\ref{thm:integral_equivalence} and transitivity.
\end{proof}
When $\bb=\kk= (k,k,\hdots,k)$ for some positive integer $k$, then $\aa=\aa'=(k,0,\hdots,0,-k)$ and we have the following. 

\begin{corollary}\label{cor:rowcolumn}
Let $G$ and $G'$ be two graphs whose associated matrices $M$ and $M'$ differ by a reordering of their rows and/or columns. Then $\F_G(k,0,\hdots,0,-k)$ is integrally equivalent to $\F_{G'}(k,0,\hdots,0,-k)$. In particular, $\F_G$ is integrally equivalent to $\F_{G'}$.
\end{corollary}

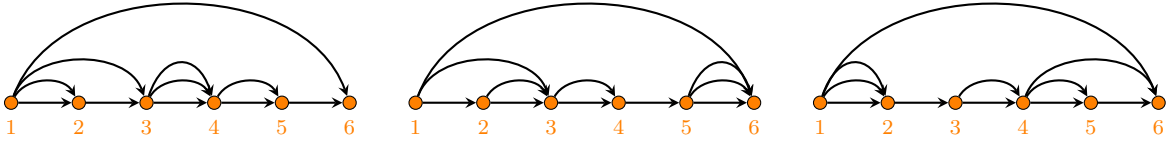
\begin{figure}[t]
    \centering
\begin{tikzpicture}[scale=0.9]

\begin{scope}[xshift=0, yshift=0]
    \vertex[fill=orange, label=below:{\tiny\textcolor{orange}{$1$}}](v1) at (3,0) {};
	\vertex[fill=orange, label=below:{\tiny\textcolor{orange}{$2$}}](v2) at (4,0) {};
	\vertex[fill=orange, label=below:{\tiny\textcolor{orange}{$3$}}](v3) at (5,0) {};
	\vertex[fill=orange, label=below:{\tiny\textcolor{orange}{$4$}}](v4) at (6,0) {};
	\vertex[fill=orange, label=below:{\tiny\textcolor{orange}{$5$}}](v5) at (7,0) {};
	\vertex[fill=orange, label=below:{\tiny\textcolor{orange}{$6$}}](v6) at (8,0) {};
    \draw[thick, -stealth] (v1) to (v2);
    \draw[thick, -stealth] (v2) to (v3);
    \draw[thick, -stealth] (v3) to (v4);
    \draw[thick, -stealth] (v4) to (v5);
    \draw[thick, -stealth] (v5) to (v6);
    \draw[thick, -stealth] (v1) to [out=70,in=110] (v6);
    \draw[thick, -stealth] (v1) to [out=70,in=110] (v3);
    \draw[thick, -stealth] (v1) to [out=60,in=120] (v2);
    \draw[thick, -stealth] (v4) to [out=60,in=120] (v5);
    \draw[thick, -stealth] (v3) to [out=60,in=120] (v4);
    \draw[-stealth, thick] (v3) .. controls (5.25,0.75) and (5.75,0.75) .. (v4);
    \end{scope}
    \begin{scope}[xshift=170, yshift=0]
    \vertex[fill=orange, label=below:{\tiny\textcolor{orange}{$1$}}](v1) at (3,0) {};
	\vertex[fill=orange, label=below:{\tiny\textcolor{orange}{$2$}}](v2) at (4,0) {};
	\vertex[fill=orange, label=below:{\tiny\textcolor{orange}{$3$}}](v3) at (5,0) {};
	\vertex[fill=orange, label=below:{\tiny\textcolor{orange}{$4$}}](v4) at (6,0) {};
	\vertex[fill=orange, label=below:{\tiny\textcolor{orange}{$5$}}](v5) at (7,0) {};
	\vertex[fill=orange, label=below:{\tiny\textcolor{orange}{$6$}}](v6) at (8,0) {};
    \draw[thick, -stealth] (v1) to (v2);
    \draw[thick, -stealth] (v2) to (v3);
    \draw[thick, -stealth] (v3) to (v4);
    \draw[thick, -stealth] (v4) to (v5);
    \draw[thick, -stealth] (v5) to (v6);
    \draw[thick, -stealth] (v1) to [out=70,in=110] (v6);
    \draw[thick, -stealth] (v1) to [out=70,in=110] (v3);
    \draw[thick, -stealth] (v2) to [out=60,in=120] (v3);
    \draw[thick, -stealth] (v3) to [out=60,in=120] (v4);
    \draw[thick, -stealth] (v5) to [out=60,in=120] (v6);
    \draw[-stealth, thick] (v5) .. controls (7.25,0.75) and (7.75,0.75) .. (v6);
    \end{scope}
\begin{scope}[xshift=340, yshift=0]
    \vertex[fill=orange, label=below:{\tiny\textcolor{orange}{$1$}}](v1) at (3,0) {};
	\vertex[fill=orange, label=below:{\tiny\textcolor{orange}{$2$}}](v2) at (4,0) {};
	\vertex[fill=orange, label=below:{\tiny\textcolor{orange}{$3$}}](v3) at (5,0) {};
	\vertex[fill=orange, label=below:{\tiny\textcolor{orange}{$4$}}](v4) at (6,0) {};
	\vertex[fill=orange, label=below:{\tiny\textcolor{orange}{$5$}}](v5) at (7,0) {};
	\vertex[fill=orange, label=below:{\tiny\textcolor{orange}{$6$}}](v6) at (8,0) {};
    \draw[thick, -stealth] (v1) to (v2);
    \draw[thick, -stealth] (v2) to (v3);
    \draw[thick, -stealth] (v3) to (v4);
    \draw[thick, -stealth] (v4) to (v5);
    \draw[thick, -stealth] (v5) to (v6);
    \draw[thick, -stealth] (v1) to [out=70,in=110] (v6);
    \draw[thick, -stealth] (v4) to [out=70,in=110] (v6);
    \draw[thick, -stealth] (v4) to [out=60,in=120] (v5);
    \draw[thick, -stealth] (v3) to [out=60,in=120] (v4);
    \draw[thick, -stealth] (v1) to [out=60,in=120] (v2);
    \draw[-stealth, thick] (v1) .. controls (3.25,0.75) and (3.75,0.75) .. (v2);
    \end{scope}
    \end{tikzpicture}
    \caption{Three graphs $G$ whose flow polytopes $\mathcal{F}_G$ are integrally equivalent. The last two graphs are reverses of each other.}
    \label{fig:NOTYET}
\end{figure}

Notably, Corollary~\ref{cor:rowcolumn} applies when $G'$ is the reverse of $G$. Figure~\ref{fig:NOTYET} shows three graphs whose corresponding flow polytopes are proved to be integrally equivalent because of Corollary~\ref{cor:rowcolumn}. 

\medskip
Another operation on matrices that does not change the polytope $\mathcal{B}_{M,\bb}$ when $\bb$ has constant entries is the introduction or removal of redundant rows.

Let $i$ and $i'$ be two rows of an $n\times d$ matrix $M$ and let $I$ and $I'$ be the subsets of $[d]$ that are the indices of non-zero entries of $M$ in rows $i$ and $i'$, respectively. If $I\subseteq I'$, we say that row $i$ is \emph{redundant}. 

\begin{proposition}\label{prop:redundantM}
Let $M$ be an $n\times d$ matrix and let $i$ be a redundant row. If $\widehat{M}$ is the matrix formed by removing row $i$ from $M$ then $\mathcal{B}_{M,\kk}=\mathcal{B}_{\widehat{M},\kk}$
\end{proposition}

\begin{proof}
Viewing $\mathcal{B}_{M,\kk}$ as an intersection of halfspaces, we see that the inequality $\sum_{i\in I} x_i \leq k$ is implied by an inequality $\sum_{i\in I'} x_i \leq k$ when $I\subseteq I'$ because $x_i\geq 0$ for all $i$. As a consequence, removing the former inequality does not change the polytope.
\end{proof}

Removing a row of a column-convex matrix preserves column convexity so we can ask how this operation impacts the associated graph $G$ and flow polytope $\mathcal{F}_{G}(k,0,\hdots,0,-k)$. 

\begin{proposition}
\label{prop:matrixrowremoval}
Let $M$ be a column-convex matrix with associated graph $G$. Let $\hat{M}$ be the matrix resulting by removing row $i$ from $M$ (and subsequently removing any identically zero columns) and let $\widehat{G}$ be its associated graph. The graph $\widehat{G}$ is formed from $G$ by deleting any multiple copies of $(i,i+1)$ and then contracting $(i,i+1)$. Furthermore, row $i$ of $M$ is redundant if and only if $G$ does not have simultaneously non-slack edges that originate at vertex $i$ and non-slack edges that terminate at vertex $i+1$. [This includes multiple edges of the form $(i,i+1)$.] 
\end{proposition}

\begin{proof}
That $\widehat{G}$ is formed from $G$ as described follows directly from Definition~\ref{def:associated}. 

If row $i$ is redundant, its non-zero indices are a subset of the non-zero indices of another row $i^*$. This implies that all non-slack edges that traverse the vertical line between $i$ and $i+1$ must also traverse the vertical line between $i^*$ and $i^*+1$. This implies that $G$ does not have both non-slack edges that originate at vertex $i$ and non-slack edges that terminate at vertex $i+1$.

If $G$ does not have both non-slack edges that originate at vertex $i$ and non-slack edges that terminate at vertex $i+1$, then either (a) all edges that originate at or before vertex $i$ do not terminate before vertex $i+2$ or (b) all edges that terminate at or after vertex $i+1$ do not originate before vertex $i-1$. In the former case, the entries of row $i$ in $M$ are a subset of the entries of row $i+1$ in $M$; in the latter, they are a subset of the entries of row $i-1$. In both cases row $i$ is redundant.
\end{proof}

A direct consequence of Propositions~\ref{prop:redundantM} and \ref{prop:matrixrowremoval} is that certain edge contractions (or their inverse operations, vertex expansions) yield integrally equivalent flow polytopes.  
A special case of the edge contraction result appears in~\cite[Lemma~2.2]{MS20}.

\begin{corollary}\label{cor:contractedflow}
Let $i\in[n]$. Let $G$ be a graph on $n+1$ vertices and let $\widehat{G}$ be the graph formed by contracting edge $(i,i+1)$. The flow polytopes $\mathcal{F}_{G}(k,0,\hdots,0,-k)$ and $\mathcal{F}_{\widehat{G}}(k,0,\hdots,0,-k)$ are integrally equivalent if and only if $G$ does not have simultaneously non-slack edges that originate at vertex $i$ and non-slack edges that terminate at vertex $i+1$.
\end{corollary}

\begin{remark}
The condition on $G$ in Propositions~\ref{prop:redundantM} and Corollary~\ref{cor:contractedflow} becomes intuitive in the flow polytope setting. Consider the following subgraph of a graph $G$ that has an edge $e_1$ entering vertex $i+1$ and an edge $e_2$ leaving vertex $i$. 
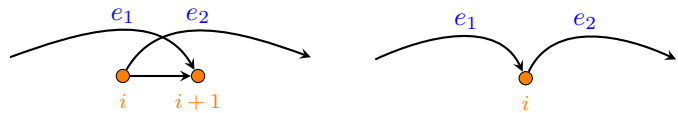
\medskip
\begin{center}
\begin{tikzpicture}[scale=1]
\begin{scope}[xshift=230, yshift=0]
    \vertex[fill=orange, label=below:{\tiny\textcolor{orange}{$i$}}](v1) at (5,0) {};
	\vertex[fill=orange, label=below:{\tiny\textcolor{orange}{$i+1$}}](v2) at (6,0) {};
    \draw[thick, -stealth] (v1) to (v2);
    \draw[thick, -stealth] (3.5,0.25) to [out=20,in=120] (v2);
    \draw[thick, -stealth] (v1) to [out=60,in=160] (7.5,0.25);
    \node at (5,0.8){\footnotesize\textcolor{blue}{$e_1$}};
    \node at (6,0.8){\footnotesize\textcolor{blue}{$e_2$}};
\end{scope}
\end{tikzpicture}\qquad
\begin{tikzpicture}[scale=1]
\begin{scope}[xshift=230, yshift=0]
    \vertex[fill=orange, label=below:{\tiny\textcolor{orange}{$i$}}](v1) at (5.5,0) {};
    \draw[thick, -stealth] (3.5,0.25) to [out=25,in=115] (v1);
    \draw[thick, -stealth] (v1) to [out=65,in=155] (7.5,0.25);
    \node at (4.71,0.76){\footnotesize\textcolor{blue}{$e_1$}};
    \node at (6.29,0.76){\footnotesize\textcolor{blue}{$e_2$}};
\end{scope}
\end{tikzpicture}
\end{center}

The number of integer flows in the two graphs is different because in the contracted graph (on the right) there can be a unit flow passing through both edges $e_1$ and $e_2$; however, this does not exist in the original graph (on the left). Similarly, if $(i,i+1)$ is a multiple edge in $G$, then the removal of copies of $(i,i+1)$ and its subsequent contraction also changes the number of integer flows.
\end{remark}

\begin{figure}[t]
\begin{tikzpicture}[scale=1]
\begin{scope}[xshift=0, yshift=0]

	\vertex[fill=orange, label=below:{\tiny\textcolor{orange}{$1$}}](v1) at (1,0) {};
	\vertex[fill=orange, label=below:{\tiny\textcolor{orange}{$2$}}](v2) at (2,0) {};
	\vertex[fill=orange, label=below:{\tiny\textcolor{orange}{$3$}}](v3) at (3,0) {};
	\vertex[fill=orange, label=below:{\tiny\textcolor{orange}{$4$}}](v4) at (4,0) {};
	\vertex[fill=orange, label=below:{\tiny\textcolor{orange}{$5$}}](v5) at (5,0) {};
	\vertex[fill=orange, label=below:{\tiny\textcolor{orange}{$6$}}](v6) at (6,0) {};
	\vertex[fill=orange, label=below:{\tiny\textcolor{orange}{$7$}}](v7) at (7,0) {};
	\vertex[fill=orange, label=below:{\tiny\textcolor{orange}{$8$}}](v8) at (8,0) {};
	\vertex[fill=orange, label=below:{\tiny\textcolor{orange}{$9$}}](v9) at (9,0) {};
	\vertex[fill=orange, label=below:{\tiny\textcolor{orange}{$10$}}](v10) at (10,0) {};
	
    \draw[thick, -stealth] (v1) to (v2);
    \draw[thick, -stealth] (v2) to (v3);
    \draw[thick, -stealth] (v3) to (v4);
    \draw[thick, -stealth] (v4) to (v5);
    \draw[thick, -stealth] (v5) to (v6);
    \draw[thick, -stealth] (v6) to (v7);
    \draw[thick, -stealth] (v7) to (v8);
    \draw[thick, -stealth] (v8) to (v9);
    \draw[thick, -stealth] (v9) to (v10);
    \draw[thick, -stealth] (v1) to [out=60,in=120] (v4);
    \draw[thick, -stealth] (v2) to [out=60,in=120] (v5);
    \draw[thick, -stealth] (v3) to [out=60,in=120] (v6);
    \draw[thick, -stealth] (v4) to [out=60,in=120] (v7);
    \draw[thick, -stealth] (v5) to [out=60,in=120] (v8);
    \draw[thick, -stealth] (v6) to [out=60,in=120] (v9);
    \draw[thick, -stealth] (v7) to [out=60,in=120] (v10);
    
    \node at (2.5,1){\footnotesize{$y_1$}};
    \node at (3.5,1){\footnotesize{$y_2$}};
    \node at (4.5,1){\footnotesize{$y_3$}};
    \node at (5.5,1){\footnotesize{$y_4$}};
    \node at (6.5,1){\footnotesize{$y_5$}};
    \node at (7.5,1){\footnotesize{$y_6$}};
    \node at (8.5,1){\footnotesize{$y_7$}};
    
    \node at (1.5,-.2){\footnotesize{$x_1$}};
    \node at (2.5,-.2){\footnotesize{$x_2$}};
    \node at (3.5,-.2){\footnotesize{$x_3$}};    
    \node at (4.5,-.2){\footnotesize{$x_4$}};    
    \node at (5.5,-.2){\footnotesize{$x_5$}};    
    \node at (6.5,-.2){\footnotesize{$x_6$}};
    \node at (7.5,-.2){\footnotesize{$x_7$}};
    \node at (8.5,-.2){\footnotesize{$x_8$}};
    \node at (9.5,-.2){\footnotesize{$x_9$}};  
\end{scope}
\begin{scope}[xshift=230, yshift=0]
    \vertex[fill=orange, label=below:{\tiny\textcolor{orange}{$3$}}](v3) at (3,0) {};
	\vertex[fill=orange, label=below:{\tiny\textcolor{orange}{$4$}}](v4) at (4,0) {};
	\vertex[fill=orange, label=below:{\tiny\textcolor{orange}{$5$}}](v5) at (5,0) {};
	\vertex[fill=orange, label=below:{\tiny\textcolor{orange}{$6$}}](v6) at (6,0) {};
	\vertex[fill=orange, label=below:{\tiny\textcolor{orange}{$7$}}](v7) at (7,0) {};
	\vertex[fill=orange, label=below:{\tiny\textcolor{orange}{$8$}}](v8) at (8,0) {};
    \draw[thick, -stealth] (v3) to (v4);
    \draw[thick, -stealth] (v4) to (v5);
    \draw[thick, -stealth] (v5) to (v6);
    \draw[thick, -stealth] (v6) to (v7);
    \draw[thick, -stealth] (v7) to (v8);
    \draw[thick, -stealth] (v3) to [out=60,in=120] (v4);
    \draw[thick, -stealth] (v3) to [out=60,in=120] (v5);
    \draw[thick, -stealth] (v3) to [out=60,in=120] (v6);
    \draw[thick, -stealth] (v4) to [out=60,in=120] (v7);
    \draw[thick, -stealth] (v5) to [out=60,in=120] (v8);
    \draw[thick, -stealth] (v6) to [out=60,in=120] (v8);
    \draw[thick, -stealth] (v7) to [out=60,in=120] (v8);
    \node at (3.5,.2){\footnotesize{$y_1$}};
    \node at (4,.7){\footnotesize{$y_2$}};
    \node at (4.5,1){\footnotesize{$y_3$}};
    \node at (5.5,1){\footnotesize{$y_4$}};
    \node at (6.5,1){\footnotesize{$y_5$}};
    \node at (7,.7){\footnotesize{$y_6$}};
    \node at (7.5,.2){\footnotesize{$y_7$}};
    \node at (3.5,-.2){\footnotesize{$x_3$}};    
    \node at (4.5,-.2){\footnotesize{$x_4$}};    
    \node at (5.5,-.2){\footnotesize{$x_5$}};  
    \node at (6.5,-.2){\footnotesize{$x_6$}};
    \node at (7.5,-.2){\footnotesize{$x_7$}};
\end{scope}
\end{tikzpicture}

    \caption{These two graphs $G$ and $G'$ are related by slack edge contractions that do not change the crossing set, so they have integrally equivalent flow polytopes $\mathcal{F}_{G}(k,0,\hdots,0,-k)$ and $\mathcal{F}_{G'}(k,0,\hdots,0,-k)$.}
    \label{fig:G310}
\end{figure}
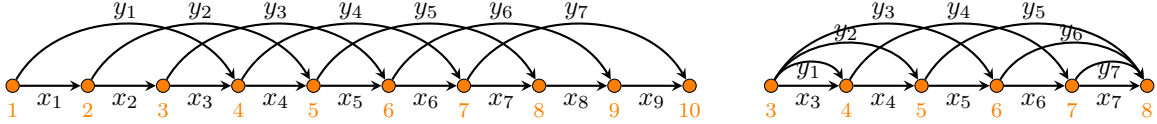

\begin{example}
Figure~\ref{fig:G310} shows a graph $G$ and a graph $G'$ that is the result of successively contracting slack edges $x_1$, $x_2$, $x_8$, and $x_9$. Because vertices $2$ and $3$ are not terminal vertices of any non-slack edges and vertices $8$ and $9$ are not originating vertices of any non-slack edges, the flow polytopes $\mathcal{F}_{G}(k,0,\hdots,0,-k)$ and $\mathcal{F}_{G'}(k,0,\hdots,0,-k)$ are integrally equivalent. The associated matrices $M$ and $M'$ are
\[M=\begin{bmatrix}
1 & 0 & 0 & 0 & 0 & 0 & 0 \\
1 & 1 & 0 & 0 & 0 & 0 & 0 \\
1 & 1 & 1 & 0 & 0 & 0 & 0 \\
0 & 1 & 1 & 1 & 0 & 0 & 0 \\
0 & 0 & 1 & 1 & 1 & 0 & 0 \\
0 & 0 & 0 & 1 & 1 & 1 & 0 \\
0 & 0 & 0 & 0 & 1 & 1 & 1 \\
0 & 0 & 0 & 0 & 0 & 1 & 1 \\
0 & 0 & 0 & 0 & 0 & 0 & 1 \\
\end{bmatrix}
\textup{ and }
M'=\begin{bmatrix}
1 & 1 & 1 & 0 & 0 & 0 & 0 \\
0 & 1 & 1 & 1 & 0 & 0 & 0 \\
0 & 0 & 1 & 1 & 1 & 0 & 0 \\
0 & 0 & 0 & 1 & 1 & 1 & 0 \\
0 & 0 & 0 & 0 & 1 & 1 & 1 \\
\end{bmatrix}.
\]
The graph $G$ is an example from the family of distance graphs, which are discussed in depth in Section~\ref{sec:distance}. 
\end{example}

Applying a sequence of edge contractions can be useful because it reduces the ambient dimension of the flow polytope. On the other hand, applying a sequence of vertex expansions can be useful because the resulting graph can be made to have no multiple edges and also ensure certain properties of the in-degree of vertices of the resulting graph. Furthermore, such operations can be used to ensure that the associated matrices have a desired form.

\medskip
We have seen that there can be a wide variety of graphs whose flow polytopes are integrally equivalent.  This leads to the following open question.

\begin{question}\label{question:equivalent}
Characterize all graphs $G$ that have integrally equivalent flow polytopes. \end{question}

\section{The volume of flow polytopes}
\label{sec:volume}

We now aim to determine the volume of the integrally equivalent polytopes of Theorem \ref{thm:integral_equivalence}. 
The starting point is a formula for the volume of flow polytopes as a Kostant partition function. This result and its generalization (Theorem~\ref{thm.genlidskii}) are known as the \emph{Lidskii volume formulas}, proved by Baldoni and Vergne~\cite[Theorem 38]{BV} via computations of residues, and also proved by both M\'esz\'aros and Morales~\cite[Theorem 1.1]{MM19} and by Kapoor-M\'esz\'aros-Setiabrata \cite{KMS19} using polytope subdivisions.

Let $G$ be a graph on $[n+1]$. For $i\in[n+1]$, let $\ee_i$ denote the $i$-th standard basis vector in $\RR^{n+1}$. For $i\in[n]$, let $\alpha_i = \ee_i - \ee_{i+1}$ denote the simple roots in the type A root system.
To each edge $e=(i,i')\in E(G)$, we associate the positive root
$$\alpha_e = \alpha_{(i,i')} = \alpha_i + \cdots + \alpha_{i'-1}. $$
Let $\Phi_G^+=\{\alpha_e \mid e\in E(G) \}$ denote the multiset of positive roots associated to $G$.
An $\aa$-flow on the graph $G$ is equivalent to expressing $\aa$ as a nonnegative linear combination of the positive roots associated to $G$.  
When the $\aa$-flow is integral, the flow is then equivalent to a vector partition of the vector $\aa$ with respect to the set of positive roots $\Phi_G^+$.
The number of integral $\aa$-flows on $G$ is the \emph{Kostant partition function} of $G$ evaluated at $\aa$, denoted by $K_G(\aa)$.
Note that $K_G(\aa)$ is also the number of integer points in $\F_G(\aa)$.

Theorem~\ref{thm.genlidskii} relates the volume of flow polytopes with Kostant partition functions. Given weak compositions $\ss=(s_1,\ldots,s_n)$ and $\tt = (t_1,\ldots,t_n)$ of $N$ we say $\ss$ \emph{dominates} $\tt$ and write $\ss \rhd \tt$ if $\sum_{i=1}^k s_i \geq \sum_{i=1}^k t_i$ for every $k\geq 1$.

\begin{theorem}[Baldoni and Vergne~{\cite[Theorem 38]{BV}}] \label{thm.genlidskii}
Let $G$ be a directed graph on the vertex set $[n+1]$ with $m$ edges, such that the out-degree of each vertex in $\{1,\ldots, n\}$ is at least one.
Let $t_i=\mathrm{outdeg}_G(i)-1$ for $i=1,\ldots, n$, and $\tt=(t_1,\ldots, t_n)$. 
If $\aa=(a_1,\ldots, a_n, -\sum_{i=1}^n a_i)$ where $a_1,\ldots,a_n$ are nonnegative integers, then the volume of the flow polytope $\F_G(\aa)$ is 
\begin{equation}
    \label{eq: lidskii outdegree}
\vol \F_G(\aa) = \sum_{\ss} \binom{m-n}{\ss} \cdot a_1^{s_1}\cdots a_n^{s_n} \cdot K_G( \ss-\tt),
\end{equation}
where the sum is over weak compositions $\ss=(s_1,\ldots,s_n)$ of $m-n$ that dominate $\tt$ and where \[K_G(\ss-\tt)=K_G(s_1-t_1,\ldots,s_n-t_n,0).\]
\end{theorem}

When $\aa=(1,0,0,\ldots,0,-1)$ the Lidskii volume formula has the following compact form (after reversing the graph, see~\cite[Corollary 1.4]{MM19}).

\begin{theorem}[{Postnikov-Stanley~\cite{Stanley_slides}, Baldoni and Vergne~\cite[Theorem 38]{BV}}] \label{thm.lidskii}
Let $G$ be a directed graph on the vertex set $[n+1]$ with $m$ edges, such that the in-degree of each vertex in $\{2,\ldots, n+1\}$ is at least one.  
Let $u_i=\mathrm{indeg}_G(i)-1$ for $i=2,\ldots, n+1$. 
The volume of the flow polytope $\F_G$ is
$$\vol\F_G = K_G\left(0,u_2,\ldots,u_n, - \sum_{i=2}^n u_i\right),$$
where $K_G$ is the Kostant partition function of $G$.
\end{theorem}

In other words, the volume of $\F_G$ (with unitary net flow $\aa=\ee_1-\ee_{n+1}$) is the number of integer points of $\F_G$ with net flow $\vv=(0,u_2,\ldots,u_n, -\sum_{i=2}^n u_i)$.

As a consequence of the Lidskii volume formula and Theorem~\ref{thm:integral_equivalence}, we have the following result.
\begin{corollary} 
\label{cor:volBM1}
Let $M$ be an $n\times d$ column-convex matrix and $\one=(1,\hdots,1)\in \mathbb{Z}^n$.  We have 
$$\vol \mathcal{B}_{M,\one} = \vol \mathcal{F}_G = K_G\left(0,u_2,\ldots,u_n, -\sum_{i=2}^n u_i\right).$$
\end{corollary}

This gives another motivation for Question~\ref{question:equivalent}:

\begin{corollary}
Let $G$ and $G'$ be graphs whose associated matrices differ by row and column reordering and/or by adding or removing redundant rows. Then
\[
K_{G}(\vv) = K_{G'}(\vv'),
\]
where $\vv$ and $\vv'$ are the in-degree vectors of $G$ and $G'$ respectively.
\end{corollary}

A special case of this identity is when the graph $G$ is reversed (see \cite[Corollary 1.4]{MM19} and \cite[Section 4]{MS}).

\section{\texorpdfstring{$G$}{}-cyclic orders}
\label{sec:G-cyclic}

Inspired by the work of Ayyer, Josuat-Verg\`es, and Ramassamy \cite{AyyerJosuatvergesRamassamy2018} on total cyclic extensions of partial cyclic orders, we introduce a new combinatorial object called a $G$-cyclic order whose enumeration gives the volume of the flow polytope $\F_G$ for any spinal graph $G$. We begin with the definition of total cyclic orders.

\begin{definition}
A \emph{partial cyclic order} on a set $X$ is a ternary relation $\gamma\subseteq X^3$ satisfying the following conditions:
\begin{enumerate}
    \item[(a)] $(x,y,z) \in \gamma$ implies $(y,z,x)\in \gamma$ (cyclicity),
    \item[(b)] $(x,y,z) \in \gamma$ implies $(z,y,x)\notin \gamma$ (asymmetry),
    \item[(c)] $(x,y,z) \in \gamma$ and $(x,z,u)\in \gamma$ implies $(x,y,u)\in \gamma$ (transitivity).
\end{enumerate}
A partial cyclic order is called a \emph{total cyclic order} if in addition it satisfies:
\begin{enumerate}
\item[(d)] for every $x,y,z \in X$, either $(x,y,z)\in \gamma$ or $(z,y,x)\in \gamma$ (comparability).
\end{enumerate}
\end{definition}

\begin{figure}[t]
\begin{center}
\begin{tikzpicture}[scale=1.0]

\begin{scope}[xshift=0, yshift=0]
\node at (0.2,0.28){$G\;\!=$};
	\vertex[fill=orange, minimum size=4pt, label=below:{\tiny\textcolor{orange}{$1$}}](v1) at (1,0) {};
	\vertex[fill=orange, minimum size=4pt, label=below:{\tiny\textcolor{orange}{$2$}}](v2) at (2,0) {};
	\vertex[fill=orange, minimum size=4pt, label=below:{\tiny\textcolor{orange}{$3$}}](v3) at (3,0) {};
	\vertex[fill=orange, minimum size=4pt, label=below:{\tiny\textcolor{orange}{$4$}}](v4) at (4,0) {};
	\vertex[fill=orange, minimum size=4pt, label=below:{\tiny\textcolor{orange}{$5$}}](v5) at (5,0) {};

	\draw[-stealth, thick] (v1)--(v2);
	\draw[-stealth, thick] (v2)--(v3);
	\draw[-stealth, thick] (v3)--(v4);
	\draw[-stealth, thick] (v4)--(v5);
	\draw[-stealth, thick] (v1) .. controls (1.25,0.8) and (2.75,0.8) .. (v3);
	\draw[-stealth, thick] (v1) .. controls (1.25,1.25) and (3.75,1.25) .. (v4);	
	\draw[-stealth, thick] (v2) .. controls (2.25,0.5) and (2.75,0.5) .. (v3);
	\draw[-stealth, thick] (v3) .. controls (3.25,0.8) and (4.75,0.8) .. (v5);
	\draw[-stealth, thick] (v3) .. controls (3.25,1.25) and (4.75,1.25) .. (v5);
	\draw[-stealth, thick] (v4) .. controls (4.25,0.5) and (4.75,0.5) .. (v5);
	
	\node at (1.5, -0.2){\footnotesize$x_1$};
	\node at (2.5, -0.2){\footnotesize$x_2$};
	\node at (3.5, -0.2){\footnotesize$x_3$};
	\node at (4.5, -0.2){\footnotesize$x_4$};
	\node at (2.5, 0.7){\footnotesize$y_1$};
	\node at (2.5, 1.1){\footnotesize$y_2$};
	\node at (2.5, 0.22){\footnotesize$y_3$};
	\node at (4, 0.75){\footnotesize$y_4$};
	\node at (4, 1.1){\footnotesize$y_5$};
	\node at (4.5, 0.22){\footnotesize$y_6$};
\end{scope}
\begin{scope}[xshift=250, yshift=10, scale=0.4]
\node at (-4.5,-.3){$\gamma=$};
\draw (0,0) circle (2cm);
	\vertex[fill, circle, inner sep=0, minimum size=4pt, label=above:\tiny$0$] at (0,2) {}; 
	\vertex[fill, circle, inner sep=0, minimum size=4pt,
		label=right:\textcolor{black}{\tiny$1$}] at (1.564,1.247) {};
	\vertex[fill, circle, inner sep=0, minimum size=4pt,
		label=right:\textcolor{black}{\tiny$2$}] at (1.95,-.445) {};
	\vertex[fill, circle, inner sep=0, minimum size=4pt,
		label=right:\textcolor{black}{\tiny$5$}] at (.868,-1.802) {};
	\vertex[fill, circle, inner sep=0, minimum size=4pt,
		label=left:\textcolor{black}{\tiny$3$}] at (-.868,-1.802) {};
	\vertex[fill, circle, inner sep=0, minimum size=4pt, 
		label=left:\textcolor{black}{\tiny$6$}] at (-1.95,-.445) {};
	\vertex[fill, circle, inner sep=0, minimum size=4pt,
		label=left:\textcolor{black}{\tiny$4$}] at (-1.564,1.247) {};
\end{scope}
\end{tikzpicture}
\end{center}
\caption{The graph $G$ and total cyclic order $\gamma$ used in Example~\ref{ex:active}.}
\label{fig:G_and_gamma}
\end{figure}
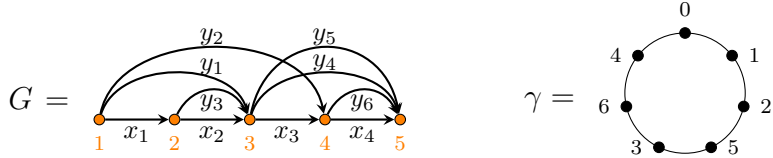

A total cyclic order can be represented visually by placing the elements of $X$ on a circle, as in the drawing of $\gamma$ in Figure~\ref{fig:G_and_gamma}.  We choose the convention of reading the elements in clockwise order, so we would read $\gamma$ as $(0,1,2,5,3,6,4)$. This sequence loops back around to the start, so that $0$ occurs directly after $4$. In a partial cyclic order one can also define the notion of chains. A \emph{chain} in $\gamma$ is a sequence $(j_1,j_2,\dots,j_k)$ such that $(j_1, j_i,j_{i+1}) \in \gamma$ for all $i = 2,\dots,k-1$. A total cyclic order is a partial cyclic order that has a unique maximal chain. Let $\H$ be a collection of chains. We write $\mathcal{A}_{\H}$ for the set of total cyclic extensions of the partial cyclic orders induced by $\H$ and $A_{\H}$ for its cardinality.  We refer the reader to \cite{AyyerJosuatvergesRamassamy2018} for an extended account on partial and total cyclic orders.

\medskip
Note that after a permutation of the edge indices, the non-slack edges of a graph $G$ can be labeled in a \emph{canonical order} by the lexicographic order on its vertex pairs, as in Figure~\ref{fig:G_and_gamma}. 
We will assume that going forward all graphs have their non-slack edges labeled in canonical order. Next, in preparation for defining total cyclic orders, compatible with a graph $G$ we need the following definitions.

\begin{definition}
Consider a graph $G$ with $d$ non-slack edges $y_1,\hdots,y_d$.  For each $j\in[d]$, we say that the index $k\in[j]$ is \emph{active at $j$}
if $\head(y_k)\le \tail(y_j)$
and say it is \emph{inactive at $j$} otherwise. Denote the set of all active indices at $j$ by $\act(j)$ and the set of inactive indices at $j$ by $\inact(j)$. By convention, let $0\in \act(j)$ for all $j\in[d]$. 
\end{definition}

Note that $\act(j)$ and $\inact(j)$ partition $\{0,1,\dots, j\}$ and 
$j\in\inact(j)$ for all $j\in[d]$.  

\begin{example}\label{ex:active}
The graph $G$ in Figure~\ref{fig:G_and_gamma} has six non-slack edges. This table shows the set of active indices $\act(j)$ and the set of inactive indices $\inact(j)$ for each $j=1,\ldots, 6$.
$$\begin{array}{c|cccccc}
j & 1 & 2 & 3 & 4 & 5 & 6\\ \hline
\act(j) & \{0\} & \{0\} & \{0\} & \{0,1,3\} & \{0,1,3\} & \{0,1,2,3\}  \\
\inact(j) & \{1\} & \{1,2\} & \{1,2,3\} & \{2,4\} & \{2,4,5\} & \{4,5,6\}\\
\end{array}
$$
\end{example}

\begin{definition} Let $G$ be a graph with $d$ non-slack edges and let $\gamma$ be a total cyclic order on $\{0,1,\dots,d\}$. For all $j\in[d]$, define the set
$$\Skip(j)=\{z\in \act(j)\mid (j-1,z,j)\in \gamma\}$$
and define the statistic $\skipstat(j)=|\Skip(j)|$.
\end{definition}
In other words, $\Skip(j)$ is the set of indices that are active at $j$ and lie between $j-1$ and $j$ in $\gamma$. Continuing Example~\ref{ex:active}, if we consider the total cyclic order $\gamma$ in Figure~\ref{fig:G_and_gamma}, we have calculated $\Skip(j)$ and $\skipstat(j)$ for each $j=1,\ldots,6$ and assembled them in the table below. As one should expect, $\Skip(j)\subseteq\act(j)$ for all $j$. Also notice that even though $2$ lies between $4$ and $5$ in $\gamma$, $2$ is not a member of $\Skip(5)$ because $2\notin\act(j)$.
$$\begin{array}{c|cccccc}
j & 1 & 2 & 3 & 4 & 5 & 6\\ \hline
\Skip(j) & \emptyset &  \emptyset & \emptyset & \emptyset & \{0,1\} & \{3\}  \\
\skipstat(j) & 0 & 0 & 0 & 0 & 2 & 1 
\end{array}
$$

\begin{definition}\label{def:skip<ACT}
We say that a total cyclic order $\gamma$ is a \emph{$G$-cyclic order} or that it is \emph{$G$-compatible} if for all $j\in[d]$ \begin{equation}\label{eq:skip}
    \sum_{k\in \inact(j)} \skipstat(k)<|\act(j)|.
\end{equation}
\end{definition}

The total cyclic order $\gamma$ in Figure~\ref{fig:G_and_gamma} is a $G$-cyclic order because Equation~\eqref{eq:skip} is satisfied for $j=1,\ldots,6$.

\begin{remark}
The intuition behind the $\skipstat$ statistic is that it models the way in which flow traveling through $G$ is temporarily ``captured'' when it enters a non-slack edge and ``released'' once the edge terminates. In the cyclic order, the indices skipped by $j$ become captured until $j$ becomes active. Definition~\ref{def:skip<ACT} is a translation of the capacity constraint on the amount of flow through all edges of an edge cut. These ideas play a central role in the proof of Thoerem~\ref{thm.volume}.
\end{remark}

For a total cyclic order $\gamma$ on $\{0,1,\hdots,d\}$ define $\gamma|_{j}$ to be the restriction of $\gamma$ to the numbers $\{0,1,\hdots,j\}$. There are two subfamilies of $G$-cyclic orders that are going to be of particular interest because their cardinality gives the volume of the flow polytope $\mathcal{F}_G$ according to Theorem \ref{thm.volume} below. These families are defined as follows based on where the entry $j$ occurs in $\gamma|_{j}$.

\begin{definition}\label{definition:upper_lower_G_cyclic_orders}
We say that $\gamma$ is an \emph{upper $G$-cyclic order} if entry $j$ in $\gamma|_{j}$ occurs immediately before an entry that is active at $j$. We say that $\gamma$ is a \emph{lower $G$-cyclic order} if entry $j$ in $\gamma|_{j}$ occurs immediately after $j-1$  or after an entry that is active at $j$.  

We denote by $\mathcal{A}^{\uparrow}_G$ the set of upper $G$-cyclic orders and $A^{\uparrow}_G$ its cardinality. We denote by $\mathcal{A}^{\downarrow}_G$ the set of lower $G$-cyclic orders and ${A}^{\downarrow}_G$ its cardinality.
\end{definition}

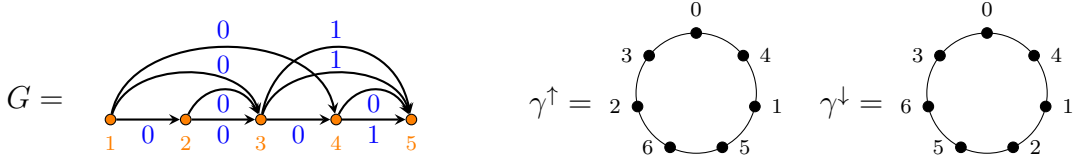
\begin{figure}[t]
\begin{center}
\begin{tikzpicture}[scale=1.0]
\begin{scope}[xshift=0, yshift=0]
\node at (0,0.3){$G=$};
	\vertex[fill=orange, minimum size=4pt, label=below:{\tiny\textcolor{orange}{$1$}}](v1) at (1,0) {};
	\vertex[fill=orange, minimum size=4pt, label=below:{\tiny\textcolor{orange}{$2$}}](v2) at (2,0) {};
	\vertex[fill=orange, minimum size=4pt, label=below:{\tiny\textcolor{orange}{$3$}}](v3) at (3,0) {};
	\vertex[fill=orange, minimum size=4pt, label=below:{\tiny\textcolor{orange}{$4$}}](v4) at (4,0) {};
	\vertex[fill=orange, minimum size=4pt, label=below:{\tiny\textcolor{orange}{$5$}}](v5) at (5,0) {};

	\draw[-stealth, thick] (v1)--(v2);
	\draw[-stealth, thick] (v2)--(v3);
	\draw[-stealth, thick] (v3)--(v4);
	\draw[-stealth, thick] (v4)--(v5);
	\draw[-stealth, thick] (v1) .. controls (1.25,0.8) and (2.75,0.8) .. (v3);
	\draw[-stealth, thick] (v1) .. controls (1.25,1.25) and (3.75,1.25) .. (v4);	
	\draw[-stealth, thick] (v2) .. controls (2.25,0.5) and (2.75,0.5) .. (v3);
	\draw[-stealth, thick] (v3) .. controls (3.25,0.8) and (4.75,0.8) .. (v5);
	\draw[-stealth, thick] (v3) .. controls (3.25,1.25) and (4.75,1.25) .. (v5);
	\draw[-stealth, thick] (v4) .. controls (4.25,0.5) and (4.75,0.5) .. (v5);
	
	\node at (1.5, -0.2){\textcolor{blue}{\footnotesize$0$}};
	\node at (2.5, -0.2){\textcolor{blue}{\footnotesize$0$}};
	\node at (3.5, -0.2){\textcolor{blue}{\footnotesize$0$}};
	\node at (4.5, -0.2){\textcolor{blue}{\footnotesize$1$}};
	\node at (2.5, 0.73){\textcolor{blue}{\footnotesize$0$}};
	\node at (2.5, 1.2){\textcolor{blue}{\footnotesize$0$}};
	\node at (2.5, 0.22){\textcolor{blue}{\footnotesize$0$}};
	\node at (4, 1.2){\textcolor{blue}{\footnotesize$1$}};
	\node at (4, 0.8){\textcolor{blue}{\footnotesize$1$}};
	\node at (4.5, 0.22){\textcolor{blue}{\footnotesize$0$}};
\end{scope}
\begin{scope}[xshift=250, yshift=10, scale=0.4]
\node at (-4.5,-.3){$\gamma^\uparrow=$};
\draw (0,0) circle (2cm);
	\vertex[fill, circle, inner sep=0, minimum size=4pt, label=above:\tiny$0$] at (0,2) {}; 
	\vertex[fill, circle, inner sep=0, minimum size=4pt,
		label=right:\textcolor{black}{\tiny$4$}] at (1.564,1.247) {};
	\vertex[fill, circle, inner sep=0, minimum size=4pt,
		label=right:\textcolor{black}{\tiny$1$}] at (1.95,-.445) {};
	\vertex[fill, circle, inner sep=0, minimum size=4pt,
		label=right:\textcolor{black}{\tiny$5$}] at (.868,-1.802) {};
	\vertex[fill, circle, inner sep=0, minimum size=4pt,
		label=left:\textcolor{black}{\tiny$6$}] at (-.868,-1.802) {};
	\vertex[fill, circle, inner sep=0, minimum size=4pt, 
		label=left:\textcolor{black}{\tiny$2$}] at (-1.95,-.445) {};
	\vertex[fill, circle, inner sep=0, minimum size=4pt,
		label=left:\textcolor{black}{\tiny$3$}] at (-1.564,1.247) {};
\end{scope}
\begin{scope}[xshift=360, yshift=10, scale=0.4]
\node at (-4.5,-.3){$\gamma^\downarrow=$};
\draw (0,0) circle (2cm);
	\vertex[fill, circle, inner sep=0, minimum size=4pt, label=above:\tiny$0$] at (0,2) {}; 
	\vertex[fill, circle, inner sep=0, minimum size=4pt,
		label=right:\textcolor{black}{\tiny$4$}] at (1.564,1.247) {};
	\vertex[fill, circle, inner sep=0, minimum size=4pt,
		label=right:\textcolor{black}{\tiny$1$}] at (1.95,-.445) {};
	\vertex[fill, circle, inner sep=0, minimum size=4pt,
		label=right:\textcolor{black}{\tiny$2$}] at (.868,-1.802) {};
	\vertex[fill, circle, inner sep=0, minimum size=4pt,
		label=left:\textcolor{black}{\tiny$5$}] at (-.868,-1.802) {};
	\vertex[fill, circle, inner sep=0, minimum size=4pt, 
		label=left:\textcolor{black}{\tiny$6$}] at (-1.95,-.445) {};
	\vertex[fill, circle, inner sep=0, minimum size=4pt,
		label=left:\textcolor{black}{\tiny$3$}] at (-1.564,1.247) {};
\end{scope}
\end{tikzpicture}
\end{center}
\caption{The bijections in the proof of Theorem~\ref{thm.volume} send the integer flow on the graph $G$ with net flow $\vv=(0,0,2,1,-3)$ to the upper $G$-cyclic order $\gamma^\uparrow$ and the lower $G$-cyclic order $\gamma^\downarrow$. 
}
\label{fig.flowcycbijection2}
\end{figure}

The total cyclic order in Figure~\ref{fig:G_and_gamma} is neither an upper $G$-cyclic order nor a lower $G$-cyclic order. It is not an upper $G$-cyclic order because $6$ occurs immediately before the inactive entry $4$ in $\gamma|_6$, and it is not a lower $G$-cyclic order because $5$ occurs immediately after the inactive entry $2$ in $\gamma|_5$.  On the other hand, in Figure~\ref{fig.flowcycbijection2}, $\gamma^{\uparrow}\in \mathcal{A}^{\uparrow}_G$ and $\gamma^{\downarrow}\in \mathcal{A}^{\downarrow}_G$.

We are now ready to state the main theorem of this section.

\volumethm

\begin{proof}
We prove a bijection between the set of integer flows in $\F_{G}(0,u_2,\ldots,u_n, -\sum_{i=2}^n u_i)$ and the set $\mathcal{A}^\uparrow_G$ of upper $G$-cyclic orders (and simultaneously with the set $\mathcal{A}^\downarrow_G$ of lower $G$-cyclic orders), where  $u_i = \mathrm{indeg}_G(i) -1$ for $i=2,\ldots, n$. The result then follows from Corollary~\ref{cor:volBM1}. 

We let $G$ have $n$ slack edges of the form $(i,i+1)$ labeled $x_i$ for $i=1,\ldots, n$, and $d$ non-slack edges of the form $(\tail(y_j),\head(y_j))$ for $j=1,\ldots, d$.
Given an integer flow $(f_e)_{e\in E(G)}$, we place the numbers $0,1,\ldots, d$ in a cyclic arrangement $\gamma$ in the following way.

\begin{itemize}
\item[(1)] Place $0$ to start the cyclic arrangement $\gamma|_0$.
\item[(2)] Successively obtain $\gamma|_j$ from $\gamma|_{j-1}$ by inserting $j$ as follows. Insert $j$ immediately before an element of $\act(j)$ so that $\skipstat(j)=f_{y_j}$. (In the case of lower $G$-cyclic orders, instead insert $j$ immediately after $j-1$ or an element of $\act(j)$ so that $\skipstat(j)=f_{y_j}$.) 
\item[(3)] Define $\gamma=\gamma|_{d}$.
\end{itemize}

For this to be well defined, we need to ensure for each $j$ that $f_{y_j}<|\act(j)|$, which equals one more than the number of non-slack edges that have terminated at or before $\head(y_j)$ (because $0$ is always active). Numerically, we have
\begin{align*}
    |\act(j)|&=1+\sum_{i=1}^{\head(y_j)}\big(\indeg(i)-1\big)=1+\sum_{i=1}^{\head(y_j)}u_i    >\sum_{i=1}^{\head(y_j)}u_i
    \ge f_{y_j},
\end{align*}
where the last inequality is satisfied by any valid flow $(f_e)_{e\in E(G)}$.

We now show that the resulting $\gamma$ is $G$-compatible. By construction, for all $j\in[d]$, $j$ occurs immediately before (immediately after) an entry that is active at $j$. Let $j^*$ be the largest index of an edge leaving vertex $\head(y_j)$. Every index that is inactive at $j$ is also inactive at $j^*$. Since $\skipstat(j)=f_{y_j}$, we have
$$\sum_{k\in \inact(j)}\!\!\!\skipstat(k)\le \sum_{k\in \inact(j^*)}\!\!\!\skipstat(k)=\sum_{k\in \inact(j^*)}f_{y_k}.$$

Consider now the sum of the flow on all edges that have an initial vertex in $\{1,\hdots,\head(y_j)\}$ and terminal vertex in $\{i_j+1,\hdots,n+1\}$. These edges are the edge labeled $x_{\head(y_j)}$ and all edges labeled $y_k$ for $k\in \inact(j^*)$. 
Since these edges form an edge cut so the sum of the flow on these edges is the sum of the net flow at vertices $1$ through $\head(y_j)$ so
$$\sum_{k\in \inact(j^*)}f_{y_k} \leq f_{x_{\head(j)}}+\sum_{k\in \inact(j^*)}f_{y_k}=\sum_{i=1}^{\head(y_j)}u_i<|\act(j)|.$$
Hence $\sum_{k\in \inact(j)}\skipstat(k)<|\act(j)|$ for all $j\in[d]$, showing that $\gamma$ is $G$-compatible.

Now we describe the inverse construction. Let $\gamma$ be a $G$-cyclic order and define the flow $(f_e)_{e\in E(G)}$ as follows. For every non-slack edge $y_j$, define $f_{y_j}=\skipstat(j)$.  For every slack edge $x_i$, define
\begin{equation}\label{equation:slack_flow_preserving}
    f_{x_i}=f_{x_{i-1}}+u_i+\sum_{k:\tail(y_k)=i}\skipstat(k) -\sum_{k:\head(y_k)=i}\skipstat(k)
\end{equation}
where $f_{x_0}=0$. (This notation is used to streamline the calculations; it is not part of the eventual flow.)  Note that Equation~\eqref{equation:slack_flow_preserving} is the necessary condition to have conservation of flow at the vertices $i=1,\dots, n$. We need to show that under this definition we have that $x_i \ge 0$ and that the conservation of flow is also happening at vertex $n+1$, that is, we need to show that
\begin{equation}\label{equation:conservation_at_last_vertex}
    f_{x_n}-\sum_{i=1}^nu_i+\sum_{k:\tail(y_k)=n+1}\skipstat(k)=0.
\end{equation}

We first apply the transformation obtained by adding Equation~\eqref{equation:slack_flow_preserving} for consecutive values $i'=1,\dots,i$ to get the following equivalent definition
\begin{equation}\label{equation:slack_flow_equivalent}
    f_{x_i}=\sum_{i'=1}^{i} u_{i'}-\sum_{k\in \inact(j^*)}\skipstat(k)
\end{equation}
where $j^*$ is the largest index of an edge leaving vertex $i$. Note that for $i=n$, Equation~\eqref{equation:slack_flow_equivalent} is equivalent to Equation~\eqref{equation:conservation_at_last_vertex}.
Also, we have that
\begin{align*}
     f_{x_i}&=\sum_{i'=1}^{i} u_i-\sum_{k\in \inact(j^*)}\skipstat(k)\\
     &=\sum_{i'=1}^{i}\big(\indeg(i)-1\big)-\sum_{k\in \inact(j^*)}\skipstat(k)\\
     &=|\act(j^*)|-1-\sum_{k\in \inact(j^*)}\skipstat(k)\\
     &\ge 0,
\end{align*}
where the last step follows by applying Equation~\eqref{eq:skip} for $j^*$.

Since in the two constructions above both directions are completely determined by the values of $y_j=\skipstat(j)$ for all $j=1,\dots,d$, they provide the desired bijection between the set 
of upper (lower) $G$-cyclic orders and the set of integer flows in $\F_{G}(0,u_2,\ldots,u_n, -\sum_{i=2}^n u_i)$. Therefore, we conclude that
\[\vol \F_{G} = K_G (0,u_2,\ldots,u_n, -\textstyle{\sum_{i=2}^n} u_i) = A^\uparrow_{G}={A}^{\downarrow}_G. \qedhere\]
\end{proof}

Since by Theorem \ref{thm:integral_equivalence} the polytope $\mathcal{B}_{M,\one}$ is integrally equivalent to the flow polytope $\F_{G}$, the enumeration of upper (or lower) $G$-cyclic orders gives the volume of polytopes of this type.

\begin{corollary} \label{corollary.volume_A}
For a column-convex matrix $M$ and its associated graph $G$,
$$\vol \mathcal{B}_{M,\one} = A^\uparrow_{G}=A^\downarrow_{G}.$$
\end{corollary}

We remark that although the sets $\mathcal{A}^\uparrow_G$ and $\mathcal{A}^{\downarrow}_G$ have the same cardinality, they are not the same set.
For the graph $G$ in Figure \ref{fig:G_and_gamma} there are $16$ upper $G$-cyclic orders and $16$ lower $G$-cyclic orders, and so $\vol \mathcal{F}_{G}=16$. One of the $(0,0,2,1,-3)$-integer flows on $G$ and its corresponding upper and lower $G$-cyclic orders are given in Figure~\ref{fig.flowcycbijection2}.

\section{Consecutive coordinate polytopes}
\label{sec:consec}

\subsection{Flow polytopes for consecutive coordinate polytopes}

Since the ordering of the defining inequalities of a polytope is irrelevant, we can expand the definition of $\mathcal{B}_{M,\bb}$ to apply to collections of subsets of $[d]$. For a collection $\C$ of subsets of $[d]=\{1,\ldots,d\}$ and a corresponding collection of integers $\bb=(b_I)_{I\in \C}\in \ZZ^{\C}$, we define the polytope
\begin{align*}
   \mathcal{B}_{\C,\bb}=\bigg\{(z_1,\ldots,z_d) \in \mathbb{R}_{\geq0}^d ~\bigg|~ 
\sum_{i\in I}z_{i} \leq b_{I}  \hbox{ for } I\in \C \bigg\}. 
\end{align*}
When $\bb=\one$, we denote $\mathcal{B}_{\C,\bb}$ simply by $\mathcal{B}_{\C}$.  We say that a collection $\C$ of subsets is \emph{non-redundant} if no two subsets $I$ and $I'$ of $\C$ satisfy $I\subseteq I'$. We say that a set consisting of the consecutive integers $\{i,i+1,\hdots,i'\}$ is an {\em interval} and denote it by $[i,i']$. 

\begin{remark}
After giving an ordering on a collection $\C$ of intervals, the polytope $\mathcal{B}_{\C,\bb}$ is equivalent to a polytope $\mathcal{B}_{M,\bb}$ where $M$ is a row-convex matrix. 
\end{remark}

We say that a matrix that is both row convex and column convex is {\em doubly convex}.

\begin{lemma}\label{lemma:row_convex_is_doubly_convex}For a positive integer $k$ and $\bb=(k,k,\dots,k)$, every polytope $\mathcal{B}_{M,\bb}$ associated to a row-convex matrix $M$ is also a polytope $\mathcal{B}_{M',\bb}$ associated to a doubly convex matrix $M'$.
\end{lemma}
\begin{proof}
If two intervals $I'$ and $I$ in $\C$ satisfy $I'\subseteq I$, we can remove $I'$ from $\C$ without impacting $\mathcal{B}_{\C,\bb}$ by Proposition~\ref{prop:redundantM}. Successively removing all such intervals yields a collection $\C'$ of non-redundant intervals $I=[i,i']$ that can be ordered lexicographically by the first entry. Keeping this order in the associated matrix gives a doubly convex matrix $M'$.
\end{proof}

For a graph $G$, two edges $(i,j)$ and $(k,l)$ that satisfy $i < k < l < j $ are said to be \emph{nested}. A graph without nested edges is said to be \emph{non-nested}.

\begin{lemma}\label{lemma:doubly_convex_gives_nonnesting}
The graph $G$ associated to a doubly convex matrix $M$ is non-nested.
\end{lemma}
\begin{proof}
The lexicographic order in the non-redundant rows of $M$ implies that the graph $G$ does not contain nested edges.
\end{proof}

\begin{figure}[t]
\begin{tikzpicture}[scale=1.0]
\begin{scope}[scale=0.7]
\node at (0,0){
$M=\begin{bmatrix}
1&1&0&0&0&0&0\\
0&1&1&1&0&0&0\\
0&0&1&1&1&1&0\\
0&0&0&0&1&1&1
\end{bmatrix}$
};
\end{scope}
\begin{scope}[xshift=120, yshift=-15, scale=1.2]
\node at (0.2,0.48){$G\;\!=$};
    \vertex[fill=orange, minimum size=4pt, label=below:{\tiny\textcolor{orange}{$1$}}](v1) at (1,0) {};
    \vertex[fill=orange, minimum size=4pt, label=below:{\tiny\textcolor{orange}{$2$}}](v2) at (2,0) {};
    \vertex[fill=orange, minimum size=4pt, label=below:{\tiny\textcolor{orange}{$3$}}](v3) at (3,0) {};
    \vertex[fill=orange, minimum size=4pt, label=below:{\tiny\textcolor{orange}{$4$}}](v4) at (4,0) {};
    \vertex[fill=orange, minimum size=4pt, label=below:{\tiny\textcolor{orange}{$5$}}](v5) at (5,0) {};

    \draw[thick, -stealth] (v1) to (v2);
    \draw[thick, -stealth] (v2) to (v3);
    \draw[thick, -stealth] (v3) to (v4);
    \draw[thick, -stealth] (v4) to (v5);
    \draw[thick, -stealth] (v1) to [out=60,in=120] (v2);
    \draw[thick, -stealth] (v1) to [out=60,in=120] (v3);
    \draw[thick, -stealth] (v2) to [out=60,in=120] (v4);
	\draw[-stealth, thick] (v2) .. controls (2.25,1.25) and (3.75,1.25) .. (v4);
    \draw[thick, -stealth] (v3) to [out=60,in=120] (v5);
	\draw[-stealth, thick] (v3) .. controls (3.25,1.25) and (4.75,1.25) .. (v5);
    \draw[thick, -stealth] (v4) to [out=60,in=120] (v5);
\end{scope}
\end{tikzpicture}
    \caption{The matrix $M$ and the graph $G$ associated to the collection of intervals $\C=\{[1,2], [2,4], [3,6], [5,7]\}\subseteq[7]$.}
    \label{fig.BS}
\end{figure}
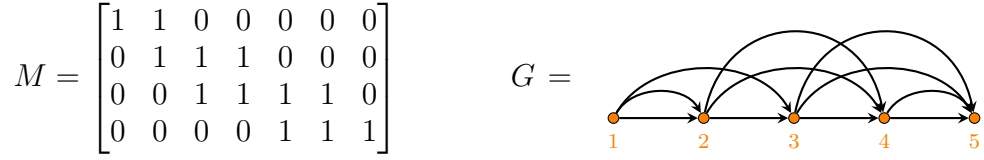

Figure~\ref{fig.BS} shows an example of a doubly convex matrix and its corresponding non-nested graph.

\smallskip
When $\C$ is a collection of intervals, $\mathcal{B}_{\C}$ is a \emph{consecutive coordinate polytope}, which is the main object of study in \cite{AyyerJosuatvergesRamassamy2018}. 

\begin{proposition}\label{proposition:intervals_are_row_convex}
Every consecutive coordinate polytope $\mathcal{B}_{\C}$ is integrally equivalent to a flow polytope $\F_G$ for a non-nested graph $G$.  
\end{proposition}

\begin{proof}
 By Lemma~\ref{lemma:row_convex_is_doubly_convex}, $\mathcal{B}_{\C}$ is integrally equivalent to $\mathcal{B}_{M,\one}$ for a doubly convex matrix $M$. By Theorem \ref{thm:integral_equivalence}, $\mathcal{B}_{M,\one}$ is integrally equivalent to $\mathcal{F}_G$ for its associated graph $G$, which is non-nested by Lemma~\ref{lemma:doubly_convex_gives_nonnesting}.
\end{proof}

Next, we give a converse of this result.

\begin{proposition} \label{proposition.non nested are consecutive}
Every flow polytope $\mathcal{F}_G$ of a non-nested graph $G$ is integrally equivalent to a consecutive coordinate polytope $\mathcal{B}_{\C}$.
\end{proposition}

\begin{proof}
Recall that in the integral equivalence of Theorem \ref{thm:integral_equivalence}, there is a defining equation of $\mathcal{B}_{\C}$ of the form 
\begin{equation}\label{eq:converse}
\sum_{i \in I}y_i=1
\end{equation}
for every slack edge $(r,r+1)$ in $G$ where $I$ is the set of indices of the edges that start at or before $r$ and end at or after $r+1$.

Suppose for the sake of contradiction that Equation~\eqref{eq:converse} satisfies $j<k<l$ with $j,l \in I$ and $k \not \in I$.
This would imply that $\head(y_k)\le r<\tail(y_j)$. Since $j<k$, the canonical order on the edges would imply that $\tail(y_j)<\head(y_k)$ so $y_k$ would be nested in $y_j$, a contradiction.
\end{proof}

\subsection{\texorpdfstring{$G$}{}-cyclic orders of non-nested graphs}

We show that when $G$ is a non-nested graph, the upper and lower $G$-cyclic orders coincide and are the same as total cyclic extensions of partial cyclic orders. We first show that all elements that are active at $j$ are smaller than all elements that are inactive at $j$ for all $j\in[d]$.

\begin{proposition}
\label{prop.active_before_inactive} Let $G$ be a non-nested graph with $d$ non-slack edges. Then for all $j\in[d]$ and for all $a\in \act(j)$ and $b \in \inact(j)$, we have that $a<b$. As a consequence, the set $$\{\max \act(j)\}\cup \inact(j)$$
is the interval $[\max \act(j),j]$ for all $j\in[d]$.
\end{proposition}

\begin{proof}
Suppose that for some $j$ it happens that $a\in \act(j)$ and $b \in \inact(j)$ with $b<a$. By the canonical labeling of the non-slack edges of $G$, $b<a$ means that $\tail(y_b)\le \tail(y_a)$. Since  $a\in \act(j)$ and $b \in \inact(j)$ we must have that $\head(y_a)\le j < \head(y_b)$. The canonical labeling of the non-slack edges then implies that $\tail(y_b)< \tail(y_a)$, which means that $y_b$ is nested inside $y_a$. This contradicts that $G$ is non-nested.
\end{proof}

\begin{proposition}
\label{prop.ofmaintheorem}
For a non-nested graph $G$, the sets $\mathcal{A}^{\uparrow}_G$ and  $\mathcal{A}^{\downarrow}_G$ are equal.
\end{proposition}

\begin{proof}
Let $G$ be a non-nested graph with an integer flow $(f_e)_{e\in E(G)}$. Consider the algorithm in the proof of Theorem~\ref{thm.volume} that constructs the upper and lover $G$-cyclic orders $\gamma^\uparrow$ and $\gamma_\downarrow$ from the flow. In Step~(2) of the algorithm, $j$ is inserted after a prescribed number of elements active at $j$. Since $G$ is non-nested, Proposition~\ref{prop.active_before_inactive} implies that all inactive elements come after all active elements, so the insertion of $j$ occurs in the same place in both the upper and lower $G$-cyclic order. Since the two bijections give the same $G$-cyclic orders, $\mathcal{A}^{\uparrow}_G$ and  $\mathcal{A}^{\downarrow}_G$ are equal.
\end{proof}

Because of this, when $G$ is a non-nested graph, we will denote the set $\mathcal{A}^{\uparrow}_G=\mathcal{A}^{\downarrow}_G$ by $\mathcal{A}_G$.

\begin{proposition}\label{prop:Gcompatibility_is_extension}
Let $G$ be a non-nested graph with $d$ non-slack edges. A total cyclic order $\gamma$ is $G$-compatible if and only if it is a total cyclic extension of the partial cyclic order whose chains are $(\max \act(j),\max \act(j)+1,\hdots,j-1,j)$ for $j\in[d]$.
\end{proposition}
\begin{proof}

First suppose $\gamma$ is $G$-compatible. We will use induction to show that for every $j \in [d]$, either $\max \act(j)=j-1$ or $(\max \act(j),\max \act(j)+1,\dots,j-1,j)$ is a chain in $\gamma$.

The base case when $j=1$ is true because $\max \act(1)=0=1-1$. The inductive hypothesis for $j-1$ along with Proposition~\ref{prop.active_before_inactive} implies either:

\begin{enumerate}
    \item $\max \act(j-1)=j-2$ and $\inact(j-1)=\{j-1\}$, or
    \item $(\max \act(j-1),\max \act(j-1)+1,\dots,j-2,j-1)$ is a chain in $\gamma$ and \newline $\inact(j-1)=[\max \act(j-1)+1,j-1]$.
\end{enumerate}

In both cases, the inequality in Equation \eqref{eq:skip} can then be rewritten as the equality
\begin{align}\label{eq.skipequality}
   u_{j-1}+\sum_{k\in \inact(j-1)} \skipstat(k)=|\act(j-1)|-1,  
\end{align}
where
$$u_{j-1}=\big\vert\big\{a\in \act(j-1)\mid (j-1, a, \max\act(j-1))\in \gamma\big\}\big\rvert.$$

For the inductive step, suppose $\max \act(j)\neq j-1$. Let $I=\inact(j-1)\cap \act(j)$; that is, $I$ is the set of elements inactive at $j-1$ that are active at $j$. Therefore $\act(j)=\act(j-1)\cup I$ and the number $N$ of $a\in\act(j)$ such that $(j-1,a,\max\act(j))\in\gamma$ is given by
$$N=u_{j-1}+\sum_{k\in I}\skipstat(k)+|I|.$$
Following Equation \eqref{eq.skipequality} we have
\begin{align*}
    N&= |\act(j-1)|-1 - \sum_{k\in \inact(j-1)\setminus I}\skipstat(k) +|I|\\
    &=|\act(j)|-1+\sum_{k\in \inact(j)\setminus \{j\}}\skipstat(k)\\
    &\ge \skipstat(j).
\end{align*}

If  $(\max\act(j),j,j-1)$ were a relation in $\gamma$ we would have that $\skipstat(j)\ge N+1$ which is a contradiction. Hence  $(\max\act(j),j-1,j)$ is a relation and by transitivity $(\max \act(j),\max \act(j)+1,\hdots,j-1,j)$ is a chain in $\gamma$. Therefore if $\gamma$ is $G$-compatible, it must be an extension of the partial cyclic order defined by the set of chains $(\max \act(j),$ $\max \act(j)+1,\hdots,j-1,j)$ for all $j\in[d]$.

\medskip
For the converse, suppose $\gamma$ is a
a total cyclic extension of the partial cyclic order whose chains are $(\max \act(j),\max \act(j)+1,\hdots,j-1,j)$ for $j\in[d]$. Again by induction on $j$ and a similar argument using Proposition \ref{prop.active_before_inactive} we see that Equation \eqref{eq:skip} is satisfied for every $j\in[d]$ and hence $\gamma$ is $G$-compatible. 
\end{proof}

\begin{example}
In Figure~\ref{fig.BS}, we can read from $G$ the value of $\max\act(j)$ for every $j$ by determining the largest edge index that lands on or before edge $y_j$ starts. This gives the following data:
$$\begin{array}{c|ccccccc}
j & 1 & 2 & 3 & 4 & 5 & 6 & 7\\ \hline
\max\act(j) & 0 & 0 & 1 & 1 & 2 & 2 & 4 \\
\end{array}
$$
The chains specified by Proposition~\ref{prop:Gcompatibility_is_extension} are $(0,1)$, $(0,1,2)$, $(1,2,3)$, $(1,2,3,4)$, $(2,3,4,5)$, $(2,3,4,5,6)$, and $(4,5,6,7)$. Therefore the $G$-cyclic orders are the ones that are total cyclic extensions of the partial cyclic orders determined by the chains $(0,1,2)$, $(1,2,3,4)$, $(2,3,4,5,6)$, and $(4,5,6,7)$.
\end{example}

Recall that if $\H$ is a set of chains then $\mathcal{A}_{\H}$ is the set of total cyclic extensions of the partial cyclic order whose chains are $\H$. Interpreting Proposition~\ref{prop:Gcompatibility_is_extension} in terms of collections of intervals we have the following.

\begin{proposition}
\label{prop:AJRresult}
Let $\C$ be a collection of intervals and let $G$ be its associated graph. Let $\H$ be the set of chains  $\big\{(i-1,i,\hdots,i') \mid [i,i']\in \C\big\}$. Then $\mathcal{A}_G=\mathcal{A}_{\H}$.
\end{proposition}
\begin{proof}
Order $\C$ in lexicographic order. If $[i,i']$ is the $r$-th interval, then row $r$ of the associated matrix $M$ has non-zero entries in columns $i$ through $i'$. This implies that in the associated graph $G$, edge $y_{i-1}$ is the last non-slack edge that terminates at or before vertex $r$ and the set of edges $y_i$ through $y_{i'}$ all terminate after vertex $r$. Applying Proposition \ref{prop:Gcompatibility_is_extension}, we notice that every vertex of $G$ contributes a unique maximal chain of the form $(i-1,i,\hdots,i')$ that a $G$-cyclic order must contain. 
\end{proof}

As a consequence of Proposition~\ref{prop:AJRresult}, we recover the following result of Ayyer, Josuat-Verg\`es, and Ramassamy.

\begin{corollary}[{\cite[Theorem 2.3]{AyyerJosuatvergesRamassamy2018}}]
Let $\C$ be a collection of intervals  and let $\H$ be the set of chains $\big\{(i-1,i,\hdots,i') \mid [i,i']\in \C\big\}$. Then $\vol \mathcal{B}_{\C}=A_{\H}$.
\end{corollary}

\medskip
\subsection{Counting non-redundant collections of intervals and graphs}~ \label{sec.counting collections}

The well-known sequence of Catalan numbers $\Cat(n):=\frac{1}{n+1}\binom{2n}{n}$ makes an appearance when enumerating non-redundant collections of intervals in $[d]$. 

\begin{proposition} \label{prop:number nonnesting}
There are $\Cat(d)$ non-redundant collections of intervals covering $[d]$. 
\end{proposition}

\begin{proof}
Such collections $\C$ are in bijection with antichains in the root poset of type $A_{d-1}$, $\C\mapsto A$ where $A=\{\ee_{i}-\ee_{i'} \,|\, [i,i']\in \C\ \text{ for } i\neq i'\}$. These antichains are counted by the Catalan number $\Cat(d)$, as shown, for example, in \cite[Corollary 1.4]{Ath} and \cite[Remark 2]{Rei97}.
\end{proof}

For a graph $G$ on vertex set $[n+1]$ with $d$ non-slack edges ordered canonically, define the \emph{edge cut set} $\{I_1,\hdots,I_n\}$ of subsets of $[d]$ where $I_r$ is the set of indices of the non-slack edges that start at or before vertex $r$ and end at or after $r+1$. We say that $G$ is \emph{non-redundant} if the edge cut set is non-redundant.

\begin{proposition} \label{prop: bijection nr collections nn graphs} 
There is a  bijection between non-redundant collections of $n$ intervals in $[d]$ and non-redundant non-nested graphs on $[n+1]$ with $d$ non-slack edges.
\end{proposition}

\begin{proof}
By the argument of Lemma~\ref{lemma:row_convex_is_doubly_convex}, to every non-redundant collection $S$ of intervals in $[d]$ we can bijectively associate a unique doubly convex matrix $M$ with $d$ columns whose rows are non-redundant and ordered in the canonical order. Using Definition \ref{def:associated} we can injectively associate to such matrix $M$ a graph $G$, that is non-nested by Lemma \ref{lemma:doubly_convex_gives_nonnesting}, with $d$ non-slack edges. 
Under this map the non-redundancy of $S$ translates into the non-redundancy of $G$. Following the arguments of the proof of Proposition \ref{proposition.non nested are consecutive} we see that from a non-redundant non-nested graph $G$ we can recover the non-redundant collection of intervals. The resultant bijection $S\mapsto G$ is such that the graph $G$ is the unique spinal graph with vertices $[|S|+1]$ and non-slack edges appears in all the intervals from the $i_j$-th interval to $i_j'$-th interval, considering intervals in their canonical order.
\end{proof}

\begin{corollary}
There are $\Cat(d)$ graphs with $d$ non-slack edges that are both non-nested and non-redundant.
\end{corollary}

\begin{proof}
The result follows from Proposition~\ref{prop:number nonnesting} and Proposition~\ref{prop: bijection nr collections nn graphs}.
\end{proof}

\begin{example}
For $d=3$ the $\Cat(3)=5$ non-redundant collections of intervals $\C$ are 
\begin{equation*}
\phantom{.}\hspace{.08in}
\{ \{1\}, \{2\}, \{3\} \} \hspace{.47in} 
\{[1,2], \{3\}\}   \hspace{.42in} 
\{ [1,3]\}  \hspace{.42in}  
\{\{1\},[2,3]\}  \hspace{.39in} 
\{[1,2], [2,3]\}
\end{equation*}
which correspond to these antichains in the root poset of type $A_2$:
\begin{equation*}
\phantom{.}\hspace{.75in}
\varnothing  \hspace{1.03in}
\{e_1-e_2\}  \hspace{.4in}
\{e_1-e_3\}  \hspace{.40in}
\{e_2-e_3\}  \hspace{.35in}
\{e_1-e_2, e_2-e_3\}
\end{equation*}
and to the following non-nested non-redundant graphs:
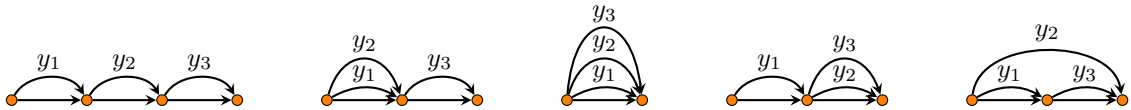
\begin{center}
\begin{tikzpicture}[scale=1.0]
    \vertex[fill=orange, minimum size=4pt](v1) at (0,0) {};
    \vertex[fill=orange, minimum size=4pt](v2) at (1,0) {};
	\vertex[fill=orange, minimum size=4pt](v3) at (2,0) {};
	\vertex[fill=orange, minimum size=4pt](v4) at (3,0) {};
	\draw[thick, -stealth](v1)--(v2);
	\draw[thick, -stealth](v2)--(v3);
	\draw[thick, -stealth](v3)--(v4);
	\draw[-stealth, thick] (v1) to [out=60,in=120] (v2);
	\draw[-stealth, thick] (v2) to [out=60,in=120] (v3);
	\draw[-stealth, thick] (v3) to [out=60,in=120] (v4);
	\node at (0.5,.5){\footnotesize{$y_1$}};
	\node at (1.5,.5){\footnotesize{$y_2$}};
	\node at (2.5,.5){\footnotesize{$y_3$}};
\end{tikzpicture}
\hspace{.3in}
\begin{tikzpicture}
    \vertex[fill=orange, minimum size=4pt](v1) at (0,0) {};
    \vertex[fill=orange, minimum size=4pt](v2) at (1,0) {};
	\vertex[fill=orange, minimum size=4pt](v3) at (2,0) {};
	\draw[thick, -stealth](v1)--(v2);
	\draw[thick, -stealth](v2)--(v3);
	\draw[-stealth, thick] (v1) to [out=30,in=150] (v2);
	\draw[-stealth, thick] (v1) .. controls (0.2, .7) and (0.75, .7) .. (v2);
	\draw[-stealth, thick] (v2) to [out=60,in=120] (v3); 
	\node at (0.5,.35){\footnotesize{$y_1$}};
	\node at (0.5,.75){\footnotesize{$y_2$}};
	\node at (1.5,.5){\footnotesize{$y_3$}};
\end{tikzpicture}
\hspace{.3in}
\begin{tikzpicture}
    \vertex[fill=orange, minimum size=4pt](v1) at (0,0) {};
    \vertex[fill=orange, minimum size=4pt](v2) at (1,0) {};
	\draw[thick, -stealth](v1)--(v2);
    \draw[-stealth, thick] (v1) to [out=30,in=150] (v2);
	\draw[-stealth, thick] (v1) .. controls (0.2, .7) and (0.75, .7) .. (v2);
	\draw[-stealth, thick] (v1) .. controls (0.1, 1.25) and (0.85, 1.25) .. (v2);
    \node at (0.5,.35){\footnotesize{$y_1$}};
	\node at (0.5,.75){\footnotesize{$y_2$}};
	\node at (0.5,1.15){\footnotesize{$y_3$}};
\end{tikzpicture}
\hspace{.3in}
\begin{tikzpicture}
    \vertex[fill=orange, minimum size=4pt](v1) at (0,0) {};
    \vertex[fill=orange, minimum size=4pt](v2) at (1,0) {};
	\vertex[fill=orange, minimum size=4pt](v3) at (2,0) {};
	\draw[thick, -stealth](v1)--(v2);
	\draw[thick, -stealth](v2)--(v3);
	\draw[-stealth, thick] (v1) to [out=60,in=120] (v2);
	\draw[-stealth, thick] (v2) .. controls (1.2, .7) and (1.75, .7) .. (v3);
    \draw[-stealth, thick] (v2) to [out=30,in=150] (v3); 
	\node at (0.5,.5){\footnotesize{$y_1$}};
	\node at (1.5,.35){\footnotesize{$y_2$}};
	\node at (1.5,.75){\footnotesize{$y_3$}};
\end{tikzpicture}
\hspace{.3in}
\begin{tikzpicture}
    \vertex[fill=orange, minimum size=4pt](v1) at (0,0) {};
    \vertex[fill=orange, minimum size=4pt](v2) at (1,0) {};
	\vertex[fill=orange, minimum size=4pt](v3) at (2,0) {};
	\draw[thick, -stealth](v1)--(v2);
	\draw[thick, -stealth](v2)--(v3);
	\draw[-stealth, thick] (v1) to [out=30,in=150] (v2);
	\draw[-stealth, thick] (v2) to [out=30,in=150] (v3);
	\draw[-stealth, thick] (v1) .. controls (0.1, 0.85) and (1.85, 0.85) .. (v3);
	\node at (0.5,.35){\footnotesize{$y_1$}};
	\node at (1.5,.35){\footnotesize{$y_3$}};
	\node at (1,0.9){\footnotesize{$y_2$}};
\end{tikzpicture}
\end{center}
\end{example}

This leads to the following open question.

\begin{question} \label{q.inherent catalan structure}
Since the family of non-redundant collections of intervals has an inherent structure (for instance, given by ordering the antichains by inclusion), does any property of the flow polytopes of the graphs corresponding to these collections behave well under this structure?
\end{question}


\section{Flow polytopes on distance graphs}\label{sec:distance}

We now focus on the study of a family of flow polytopes related to the combinatorics of Euler numbers and some of their generalizations.

\begin{figure}[t!]
\begin{tikzpicture}[scale=0.7]
\begin{scope}[xshift=0, yshift=0]
\node at (2.5,1.5){$G(1,4)$};
	\vertex[fill=orange, minimum size=4pt, label=below:{\tiny\textcolor{orange}{$1$}}](v1) at (1,0) {};
	\vertex[fill=orange, minimum size=4pt, label=below:{\tiny\textcolor{orange}{$2$}}](v2) at (2,0) {};
	\vertex[fill=orange, minimum size=4pt, label=below:{\tiny\textcolor{orange}{$3$}}](v3) at (3,0) {};
	\vertex[fill=orange, minimum size=4pt, label=below:{\tiny\textcolor{orange}{$4$}}](v4) at (4,0) {};
	\draw[-stealth, thick] (v1)--(v2);
	\draw[-stealth, thick] (v2)--(v3);
	\draw[-stealth, thick] (v3)--(v4);
	\draw[-stealth, thick] (v1) to [out=60,in=120] (v2);
	\draw[-stealth, thick] (v2) to [out=60,in=120] (v3);
	\draw[-stealth, thick] (v3) to [out=60,in=120] (v4);
\end{scope}
\begin{scope}[xshift=125, yshift=0]
\node at (3,1.5){$G(2,5)$};
		\vertex[fill=orange, minimum size=4pt, label=below:{\tiny\textcolor{orange}{$1$}}](v1) at (1,0) {};
	\vertex[fill=orange, minimum size=4pt, label=below:{\tiny\textcolor{orange}{$2$}}](v2) at (2,0) {};
	\vertex[fill=orange, minimum size=4pt, label=below:{\tiny\textcolor{orange}{$3$}}](v3) at (3,0) {};
	\vertex[fill=orange, minimum size=4pt, label=below:{\tiny\textcolor{orange}{$4$}}](v4) at (4,0) {};
	\vertex[fill=orange, minimum size=4pt, label=below:{\tiny\textcolor{orange}{$5$}}](v5) at (5,0) {};

	\draw[-stealth, thick] (v1)--(v2);
	\draw[-stealth, thick] (v2)--(v3);
	\draw[-stealth, thick] (v3)--(v4);
	\draw[-stealth, thick] (v4)--(v5);
	\draw[-stealth, thick] (v1) to [out=60,in=120] (v3);
	\draw[-stealth, thick] (v2) to [out=-60,in=-120] (v4);
	\draw[-stealth, thick] (v3) to [out=60,in=120] (v5);
\end{scope}
\begin{scope}[xshift=280, yshift=0]
\node at (5.5,1.5){$G(3,10)$};
	\vertex[fill=orange, minimum size=4pt, label=below:{\tiny\textcolor{orange}{$1$}}](v1) at (1,0) {};
	\vertex[fill=orange, minimum size=4pt, label=below:{\tiny\textcolor{orange}{$2$}}](v2) at (2,0) {};
	\vertex[fill=orange, minimum size=4pt, label=below:{\tiny\textcolor{orange}{$3$}}](v3) at (3,0) {};
	\vertex[fill=orange, minimum size=4pt, label=below:{\tiny\textcolor{orange}{$4$}}](v4) at (4,0) {};
	\vertex[fill=orange, minimum size=4pt, label=below:{\tiny\textcolor{orange}{$5$}}](v5) at (5,0) {};		
	\vertex[fill=orange, minimum size=4pt, label=below:{\tiny\textcolor{orange}{$6$}}](v6) at (6,0) {};
	\vertex[fill=orange, minimum size=4pt, label=below:{\tiny\textcolor{orange}{$7$}}](v7) at (7,0) {};
	\vertex[fill=orange, minimum size=4pt, label=below:{\tiny\textcolor{orange}{$8$}}](v8) at (8,0) {};
	\vertex[fill=orange, minimum size=4pt, label=below:{\tiny\textcolor{orange}{$9$}}](v9) at (9,0) {};
	\vertex[fill=orange, minimum size=4pt, label=below:{\tiny\textcolor{orange}{$10$}}](v10) at (10,0) {};

	\draw[-stealth, thick] (v1)--(v2);
	\draw[-stealth, thick] (v2)--(v3);
	\draw[-stealth, thick] (v3)--(v4);
	\draw[-stealth, thick] (v4)--(v5);
	\draw[-stealth, thick] (v5)--(v6);
	\draw[-stealth, thick] (v6)--(v7);
	\draw[-stealth, thick] (v7)--(v8);
	\draw[-stealth, thick] (v8)--(v9);
	\draw[-stealth, thick] (v9)--(v10);
\draw[-stealth, thick] (v1) to [out=60,in=120] (v4);
\draw[-stealth, thick] (v2) to [out=60,in=120] (v5);
\draw[-stealth, thick] (v3) to [out=60,in=120] (v6);
\draw[-stealth, thick] (v4) to [out=60,in=120] (v7);
\draw[-stealth, thick] (v5) to [out=60,in=120] (v8);
\draw[-stealth, thick] (v6) to [out=60,in=120] (v9);
\draw[-stealth, thick] (v7) to [out=60,in=120] (v10);
\end{scope}
\end{tikzpicture}
    \caption{Examples of distance graphs $G(k,d+k)$.}
    \label{fig:zigzagfamily}
\end{figure}
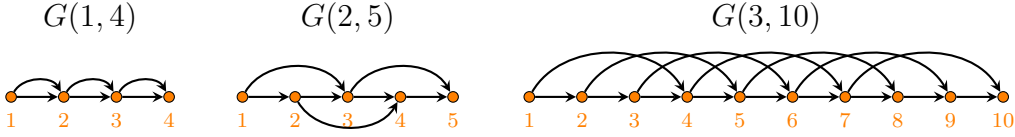

For positive integers $k$ and $d$, define the {\em distance graph} $G(k,d+k)$ to
be the graph with vertex set $[d+k]$ and $2d+k-1$ edges 
\[
\left\{(1,2),(2,3),\ldots,(d+k-1,d+k)\} \cup \{(1,k+1),(2,k+2),\ldots,(d,d+k)\right\}.
\]
See Figure~\ref{fig:zigzagfamily}. 
The graphs of the form $G(2,d+2)$ are called {\em zigzag graphs} because their flow polytopes with unit flow are integrally equivalent to the order and chain polytopes on the zigzag poset. 

\subsection{Vertices of flow polytopes of distance graphs}
We use the characterization of the vertices of flow polytopes in Proposition~\ref{prop: char verties F_G}  to enumerate the vertices of $\mathcal{F}_{G(k,d+k)}$. 
For the special cases when $k$ equals $1$ or $2$, the number of vertices of $\mathcal{F}_{G(1,d+1)}$ is $2^{d}$ and the number of vertices of $\mathcal{F}_{G(2,d+2)}$ is the Fibonacci number $F_{d+2}$. These formulas can be seen as part of the following general result. 

\begin{proposition} \label{prop: vertices G(k,d+k)}
For positive integers $d$ and $k$, the number $v_{k,d}$ of vertices of $\mathcal{F}_{G(k,d+k)}$ satisfies the recurrence $v_{k,d}=v_{k,d-1} + v_{k,d-k}$ for $d>  k$ with initial values $v_{k,d}=d+1$ for $d=1,\ldots,k$ and has generating function
\[
\sum_{d\geq 1} v_{k,d} x^d = \frac{x(2+x+x^2+\cdots + x^{k-1})}{1-x-x^k}.
\]
\end{proposition}

\begin{proof}
By Proposition~\ref{prop: char verties F_G}, the vertices of the polytope $\mathcal{F}_{G(k,d+k)}$ correspond to paths from vertex $1$ to vertex $d+k$ in $G(k,d+k)$. 
For $d=1,\ldots,k$ a path from $1$ to $d+k$ of $G(k,d+k)$ is either the path $1\to 2 \to \cdots \to d+k$ or it has exactly one non-slack edge out of the $d$ such edges. Thus $v_{k,d}=d+1$. For $d> k$, a path from vertex $1$ to $d+k$ is either of the form $P' \cup (d+k-1,d+k)$ where $P'$ is a path from $1$ to $d+k-1$ or $P'' \cup (d,d+k)$  where $P''$ is a path from $1$ to $d$, giving the desired recurrence for the number $v_{k,d}$. The generating series follows readily from the linear recurrence and the initial conditions.
\end{proof}

\subsection{Volumes of flow polytopes of distance graphs and \texorpdfstring{$k$}{}-Euler numbers}

It was shown in~\cite[Proposition 3.5]{BGHHKMY} that the volume of $\F_{G(2,d+2)}$ is the $d$-th Euler number $E_{d}$ 

Recall that $\C(k,d)$ denotes the set $\{[1,k], [2,k+1],\ldots, [d-k+1,d] \}$ of intervals of length $k$. Let $\C'$ be the collection $\C(k,d)$ together with the redundant intervals $[1,i]$ and $[d-k+1-i,d]$ for $i=1,\dots, k-1$. After ordering the intervals in $\C'$ lexicographically, $G(k,d+k)$ is the graph associated to $\C'$ by Definition \ref{def:associated}. By an application of Theorem \ref{thm:integral_equivalence}, the polytopes $\mathcal{B}_{\C'}$, $\mathcal{B}_{\C(k,d)}$, and $\mathcal{F}_{G(k,d+k)}$ are all integrally equivalent.
 
We use the simplified notation $\mathcal{A}_{k,d}$ for the set $\mathcal{A}_{\C(k,d)}$ of extensions of the partial cyclic order determined by $\C(k,d)$ and $A_{k,d}$ for its cardinality $A_{\C(k,d)}$. The next proposition follows from Proposition~\ref{prop:AJRresult}.

\begin{proposition}\label{prop:defAkd}
For any integers $d$ and $k \geq 1$ then
$$\vol \mathcal{F}_{G(k,d+k)} = \vol \mathcal{B}_{\C(k,d)}= A_{k,d}.$$
\end{proposition}

As a generalization of the case $k=2$, we see that the numbers $A_{2,d}$ can be viewed as analogue of the {\em Euler numbers}~\cite[A000111]{oeis}. We then call $A_{k,d}$ the \emph{$k$-Euler numbers}.\footnote{These numbers are different from other generalizations of Euler numbers such as \cite[Exercise 4.3.6]{MW} or \cite[A131454]{oeis}.} See Table~\ref{table:keuler}.

\begin{table}[t]
$$\begin{array}{crrrrrrrrrrrrrrr}
\hline
k \backslash d& 1&  2& 3& 4& 5& 6& 7& 8&9&10 & \\ \hline
1& 
    1& 2& 6& 24& 120& 720& 5040& 40320& 362880& 3628800 & \href{https://oeis.org/A000142}{\textup{A000142}}\\
2&	 
    1& 1& 2& 5& 16& 61& 272& 1385& 7936& 50521& \href{https://oeis.org/A000111}{\textup{A000111}}\\
3& 
    1& 1& 1& 2& 5& 14& 47& 182& 786& 3774& \href{https://oeis.org/A096402}{\textup{A096402}}\\
4&  
    1& 1& 1& 1& 2& 5& 14& 42& 146&574&\\ \hline
\end{array}
$$
\caption{Initial terms of the $k$-Euler numbers $A_{k,d}$.}
\label{table:keuler}
\end{table}

The $k$-Euler number $A_{k,d}$ can be interpreted as a number of integer flows on the distance graph $G(k,d)$.

\begin{corollary} \label{cor. k-Euler is kpf}
Let $k$ and $d$ be positive integers with $d>k$. 
We have that
\[
A_{k,d} = K_{G(k,d)}(1^{d-1},-d+1).
\]
\end{corollary}

\begin{proof}
By Proposition~\ref{prop:defAkd} and Theorem~\ref{thm.lidskii} we have that 
\[
A_{k,d} = K_{G(k,d+k)}(0^k,1^{d-1},-d+1).
\]
The net flow on the first $k$ vertices of $G(k,d+k)$ is zero for the integer flows counted on the right hand side. Thus the support of such integer flows is on the subgraph of $G(k,d+k)$ of the last $d$ vertices, which is isomorphic to $G(k,d)$, with netflow $(1^{d-1},-d+1)$. Conversely, every integer flow of $G(k,d)$ with netflow $(1^{d-1},-d+1)$ can be extended to an integer flow of $G(k,d+k)$ with netflow $(0^k,1^{d-1},-d+1)$. Thus we have the identity
\[
K_{G(k,d+k)}(0^k,1^{d-1},-d+1) = K_{G(k,d)}(1^{d-1},-d+1),
\]
which gives the desired result.
\end{proof}


\subsection{\texorpdfstring{$k$}{}-Entringer numbers}
\label{subsec.Entringer}

By partitioning the set of integer flows on $G(k,d)$ with netflow $\vv=(1^{d-1},-d+1)$ based on the flow on each of the edges in the edge cut separating the first $d+1$ vertices from the last $k-1$ vertices, we provide a new combinatorial interpretation of a refinement of the $k$-Euler numbers first found by Ayyer, Josuat-Verg\`es, and Ramassamy  \cite{AyyerJosuatvergesRamassamy2018}. 

\begin{figure}[t]
\begin{center}
\begin{tikzpicture}[scale=1.0]
\begin{scope}[xshift=0, yshift=0, scale=0.7]
    \vertex[fill=orange, minimum size=4pt](v1) at (0,0) {};
    \vertex[fill=orange, minimum size=4pt](v2) at (1,0) {};
	\vertex[fill=orange, minimum size=4pt](v3) at (2,0) {};
	\vertex[fill=orange, minimum size=4pt](v4) at (3,0) {};
	\vertex[fill=orange, minimum size=4pt](v5) at (4,0) {};
	\vertex[fill=orange, minimum size=4pt](v6) at (5,0) {};
	\vertex[fill=orange, minimum size=4pt](v7) at (6,0) {};
	\vertex[fill=orange, minimum size=4pt](v8) at (7,0) {};
	\vertex[fill=orange, minimum size=4pt](v9) at (8,0) {};
	\vertex[fill=orange, minimum size=4pt](v10) at (9,0){};
	\vertex[fill=orange, minimum size=4pt](v11) at (10,0){};
	\vertex[fill=orange, minimum size=4pt](v12) at (11,0){};
	\draw[thick, -stealth](v1)--(v2);
	\draw[thick, -stealth](v2)--(v3);
	\draw[thick, -stealth](v3)--(v4);
	\draw[thick, -stealth](v4)--(v5);
	\draw[thick, -stealth](v5)--(v6);
	\draw[thick, -stealth](v6)--(v7);
	\draw[thick, -stealth](v7)--(v8);
	\draw[thick, -stealth](v8)--(v9);
	\draw[thick, -stealth](v9)--(v10);
	\draw[thick, -stealth](v10)--(v11);
	\draw[thick, -stealth](v11)--(v12);
    \draw[thick, -stealth] (v1) to [out=60,in=120] (v5);
    \draw[thick, -stealth] (v2) to [out=60,in=120] (v6);
    \draw[thick, -stealth] (v3) to [out=60,in=120] (v7);
    \draw[thick, -stealth] (v4) to [out=60,in=120] (v8);
	\draw[thick, -stealth] (v5) to [out=60,in=120] (v9);
	\draw[thick, -stealth] (v6) to [out=60,in=120] (v10);
	\draw[thick, -stealth] (v7) to [out=60,in=120] (v11);
	\draw[thick, -stealth] (v8) to [out=60,in=120] (v12);
	\node at (7,1.33){\textcolor{blue}{\footnotesize$s_4$}};
	\node at (8,1.33){\textcolor{blue}{\footnotesize$s_3$}};
	\node at (9,1.33){\textcolor{blue}{\footnotesize$s_2$}};
	\node at (8.3,0.33){\textcolor{blue}{\footnotesize$s_1$}};
    \node at (0,-0.5){\textcolor{red}{\footnotesize$1$}};
	\node at (1,-0.5){\textcolor{red}{\footnotesize$1$}};
	\node at (2,-0.5){\textcolor{red}{\footnotesize$1$}};
	\node at (3,-0.5){\textcolor{red}{\footnotesize$1$}};
	\node at (4,-0.5){\textcolor{red}{\footnotesize$1$}};
	\node at (5,-0.5){\textcolor{red}{\footnotesize$1$}};
	\node at (6,-0.5){\textcolor{red}{\footnotesize$1$}};
	\node at (7,-0.5){\textcolor{red}{\footnotesize$1$}};
	\node at (8,-0.5){\textcolor{red}{\footnotesize$1$}};
	\node at (9,-0.5){\textcolor{red}{\footnotesize$1$}};
	\node at (10,-0.5){\textcolor{red}{\footnotesize$1$}};
	\node at (11,-0.5){\textcolor{red}{\footnotesize$-11$}};
\end{scope}

\begin{scope}[xshift=250, yshift=0, scale=0.7]
	\vertex[fill=orange, minimum size=4pt](v1) at (1,0) {};
	\vertex[fill=orange, minimum size=4pt](v2) at (2,0) {};
	\vertex[fill=orange, minimum size=4pt](v3) at (3,0) {};
	\vertex[fill=orange, minimum size=4pt](v4) at (4,0) {};
	\vertex[fill=orange, minimum size=4pt](v5) at (5,0) {};
	\vertex[fill=orange, minimum size=4pt](v6) at (6,0) {};
	\vertex[fill=orange, minimum size=4pt](v7) at (7,0) {};
	\vertex[fill=orange, minimum size=4pt](v8) at (8,0) {};
	\vertex[fill=orange, minimum size=4pt](v9) at (9,0) {};
	\draw[thick, -stealth](v1)--(v2);
	\draw[thick, -stealth](v2)--(v3);
	\draw[thick, -stealth](v3)--(v4);
	\draw[thick, -stealth](v4)--(v5);
	\draw[thick, -stealth](v5)--(v6);
	\draw[thick, -stealth](v6)--(v7);
	\draw[thick, -stealth](v7)--(v8);
	\draw[thick, -stealth](v8)--(v9);
	\draw[thick, -stealth] (v1) to [out=60,in=120] (v5);
	\draw[thick, -stealth] (v2) to [out=60,in=120] (v6);
	\draw[thick, -stealth] (v3) to [out=60,in=120] (v7);
	\draw[thick, -stealth] (v4) to [out=60,in=120] (v8);
	\draw[thick, -stealth] (v5) to [out=60,in=120] (v9);
	\node at (1,-0.5){\textcolor{red}{\footnotesize$1$}};
    \node at (2,-0.5){\textcolor{red}{\footnotesize$1$}};
	\node at (3,-0.5){\textcolor{red}{\footnotesize$1$}};
	\node at (4,-0.5){\textcolor{red}{\footnotesize$1$}};
	\node at (5,-0.5){\textcolor{red}{\footnotesize$1$}};
	\node[rotate=-20] at (6.2,-0.5){\textcolor{red}{\tiny$1\!-\!s_4$}};
	\node[rotate=-20] at (7.2,-0.5){\textcolor{red}{\tiny$1\!-\!s_3$}};
	\node[rotate=-20] at (8.2,-0.5){\textcolor{red}{\tiny$1\!-\!s_2$}};
	\node[rotate=-20] at (9.2,-0.5){\textcolor{red}{\tiny$1\!-\!s_1$}};
\end{scope}
\end{tikzpicture}
\end{center}
\caption{The two integer flow interpretations of $E_{(s_1,s_2,s_3,s_4)}$ in the graphs $G(4,12)$ and $G(4,9)$ where $(s_1,s_2,s_3,s_4)\in T^4_{9}$, as described in Proposition~\ref{prop: other interpretation E_s}.
}
\label{fig.flow_interpretations E_s}
\end{figure}
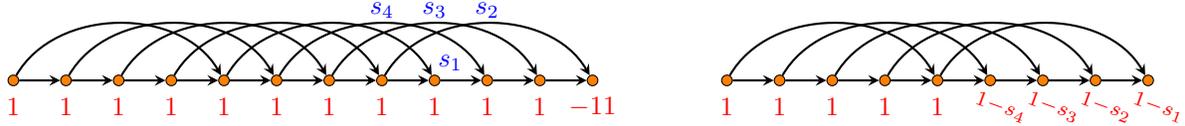

The refinement is indexed by the vectors
$$T_N^k =  \{(s_1,\ldots,s_k) \in \ZZ_{\geq0}^{k}\mid s_1+\cdots+s_k = N \},$$
where $N=d-k+1 \geq 0$. Given $\ss=(s_1,\ldots,s_k)\in T_N^k$, define $\E_{\ss}$ to be the set of integral $\vv$-flows on $G(k,d)$ whose flow on the slack edge $(d-k+1,d-k+2)$ is $s_1$ and whose flow on the last $k-1$ non-slack edges from right to left are $s_2,\ldots, s_{k}$. This is illustrated on the left side of Figure~\ref{fig.flow_interpretations E_s}. 
Define the {\em $k$-Entringer number} indexed by $\ss\in T_N^k$ to be the number $E_{\ss}=|\E_{\ss}|$ of such integer flows.  For a fixed $k$ and $N$, the numbers $E_{\ss}$ for $\ss\in T_N^k$ can be arranged into a $(k-1)$-dimensional array in the shape of a simplex, as in Figure~\ref{fig:boustrophedon_recursion}.

When $k=2$ the number $E_{(s_1, N-s_1)}$ coincides with the (classical) Entringer number $E_{N,s_1}$. When $k=1$, there is no refinement and $A_{1,d}=E_{(d)}=d!$. This agrees with the observation that the consecutive coordinate polytope $\mathcal{B}_{C(1,d)}$ is the $d$-hypercube.

The sets $\{\mathcal{E}_\ss \mid \ss\in T_N^k \}$ partition the set of integral $\vv$-flows on the distance graph $G(k,d)$, so it follows that the $k$-Euler number $A_{k,d}$ is the sum of $k$-Entringer numbers:
\begin{equation}\label{eqn.genEnt}
A_{k,d}=\sum_{ \ss\in T_N^k} E_{\ss}.
\end{equation}

\begin{remark}
In the language of~\cite{AyyerJosuatvergesRamassamy2018}, the $k$-Entringer number $a_{(i_1,\ldots,i_k)}$ is defined as the number of total cyclic orders in $\mathcal{A}_{k,d}$ such that there are $i_j-1$ numbers between $d-k+j$ and $d-k+1+j$ in the total cyclic order, for $j=1,\ldots, k-1$.
Via the bijection in Theorem~\ref{thm.volume}, this is precisely the number of integer flows on $G(k,d)$ with net flow $\vv=(1^{d-1}, -d+1)$  whose flow on the last $k-1$ non-slack edges from left to right are $i_1-1, i_2-1,\ldots, i_{k-1}-1$.
In other words,
$$E_{(s_1,s_2,\ldots,s_k)} = a_{(s_k+1,\ldots, s_2+1,s_1+1)},$$
from which we see that Equation~\eqref{eqn.genEnt} is equivalent to Equation (7.3) in~\cite{AyyerJosuatvergesRamassamy2018}.
\end{remark}

Next, we show that the $k$-Entringer numbers can also be viewed as a value of a Kostant partition function. 
This viewpoint becomes useful when we show Section~\ref{sec:logconcavity} that the numbers $E_\ss$ are log-concave along root directions.

\begin{proposition} \label{prop: other interpretation E_s}
Let $k\geq 2$ and $N=d-k+1\geq 0$. For $\ss\in T_N^k$,
\[
E_{\ss} = K_{G(k,N)}(1^{d-2k+1},1-s_k,\ldots,1-s_2,1-s_1).
\]
\end{proposition}

\begin{proof}
The restriction of the graph $G(k,d)$ to its first $d-k+1$ vertices gives a bijection between the set of integer flows in $\mathcal{E}_\ss$ on $G(k,d)$ with net flow vector $(1^{d-1},-d+1)$ and the set $\widetilde{\mathcal{E}}_\ss$ of integer flows on $G(k,d-k+1)$ with net flow vector $(1^{d-2k+1},1-s_{k},\ldots,1-s_2,1-s_1)$ since the flow on every edge of $G(k,d)$ not in $G(k,d-k+1)$ is fixed by the choice of $\ss$. (See Figure~\ref{fig.flow_interpretations E_s}.)
\end{proof}

We remark that when $s_1=0$, then $E_\ss=0$ because of the bijection in Proposition~\ref{prop: other interpretation E_s} and the fact that the net flow into vertex $d-k+1$ in the original graph $G(k,d)$ was $1$.

\subsection{The boustrophedon recursion for \texorpdfstring{$k$}{}-Entringer numbers}

By further exploiting the recursive nature of distance graphs, we next show that the $k$-Entringer numbers can be computed recursively on $k$ levels.
In other words, each $k$-Entringer number $E_\ss$ where $\sum_{i=1}^k s_i = N$ can be expressed as a partial sum of $k$-Entringer numbers indexed by entries in each of the simplices $T_N^k, T_{N-1}^k,\ldots, T_{N-k+1}^k$. 

\begin{theorem}[The $k$-boustrophedon recursion]
\label{thm:kBoustrophedon}
Let $k\geq2$, $N=d-k+1\geq1$, and $(s_1,\ldots,s_k)\in T_N^k$. For $j=1,\ldots, k-1$, we have
$$E_{(s_1,\ldots,s_k)} 
	= \sum_{\mathbf{u} } E_{(u_{j+1}+s_{j+1}, s_{j+2} ,\ldots, s_k, u_1,u_2,\ldots, u_{j})}, $$
a sum over weak compositions $\mathbf{u}=(u_1,\dots,u_{j+1})\vDash N-j -\sum_{i=j+1}^{k}s_i$ that satisfies the inequalities
$u_1+\cdots+u_h \leq s_1+\cdots+s_h-h$ for $h=1,\ldots, j$.
When $N=0$, $E_{(0,\ldots,0)}=1$.
\end{theorem}

\begin{figure}[t!]
\centering
\raisebox{5pt}{
\begin{tikzpicture}
\begin{scope}[scale=0.7]
    \vertex[fill=orange, minimum size=4pt](v0) at (0,0) {};
	\vertex[fill=orange, minimum size=4pt](v1) at (1,0) {};
	\vertex[fill=orange, minimum size=4pt](v2) at (2,0) {};
	\vertex[fill=orange, minimum size=4pt](v3) at (3,0) {};
	\vertex[fill=orange, minimum size=4pt](v4) at (4,0) {};
	\vertex[fill=orange, minimum size=4pt](v5) at (5,0) {};
	\vertex[fill=orange, minimum size=4pt](v6) at (6,0) {};
	\vertex[fill=orange, minimum size=4pt](v7) at (7,0) {};
	\vertex[fill=orange, minimum size=4pt](v8) at (8,0) {};
	\vertex[fill=orange, minimum size=4pt](v9) at (9,0) {};
	\draw[thick, -stealth](v0)--(v1);
	\draw[thick, -stealth](v1)--(v2);
	\draw[thick, -stealth](v2)--(v3);
	\draw[thick, -stealth](v3)--(v4);
	\draw[thick, -stealth](v4)--(v5);
	\draw[thick, -stealth](v5)--(v6);
	\draw[thick, -stealth](v6)--(v7);
	\draw[thick, -stealth](v7)--(v8);
	\draw[thick, -stealth](v8)--(v9);
    \draw[thick, -stealth] (v0) to [out=60,in=120] (v4);
    \draw[thick, -stealth] (v1) to [out=60,in=120] (v5);
    \draw[thick, -stealth] (v2) to [out=60,in=120] (v6);
	\draw[thick, -stealth] (v3) to [out=60,in=120] (v7);
	\draw[thick, -stealth] (v4) to [out=60,in=120] (v8);
	\draw[thick, -stealth] (v5) to [out=60,in=120] (v9);
    \node at (7.5,1.4){\textcolor{blue}{\footnotesize$u_1$}};
	\node at (6.5,1.4){\textcolor{blue}{\footnotesize$u_2$}};	
	\node at (7.35,0.3){\textcolor{blue}{\footnotesize$u_3$}};
    \node at  (0,-0.5){\textcolor{red}{\footnotesize$1$}};
	\node at  (1,-0.5){\textcolor{red}{\footnotesize$1$}};
	\node at (2,-0.5){\textcolor{red}{\footnotesize$1$}};
	\node at  (3,-0.5){\textcolor{red}{\footnotesize$1$}};
	\node at (4,-0.5){\textcolor{red}{\footnotesize$1$}};
	\node at (5,-0.5){\textcolor{red}{\footnotesize$1$}};
	\node[rotate=-20] at (6.2,-0.5){\textcolor{red}{\tiny$1\!-\!s_4$}};
	\node[rotate=-20] at (7.2,-0.5){\textcolor{red}{\tiny$1\!-\!s_3$}};
	\node[rotate=-20] at (8.2,-0.5){\textcolor{red}{\tiny$1\!-\!s_2$}};
	\node[rotate=-20] at (9.2,-0.5){\textcolor{red}{\tiny$1\!-\!s_1$}};
\end{scope}
\end{tikzpicture}}
\qquad 
\begin{tikzpicture}
\begin{scope}[scale=0.7]
    \vertex[fill=orange, minimum size=4pt](v0) at (0,0) {};
	\vertex[fill=orange, minimum size=4pt](v1) at (1,0) {};
	\vertex[fill=orange, minimum size=4pt](v2) at (2,0) {};
	\vertex[fill=orange, minimum size=4pt](v3) at (3,0) {};
	\vertex[fill=orange, minimum size=4pt](v4) at (4,0) {};
	\vertex[fill=orange, minimum size=4pt](v5) at (5,0) {};
	\vertex[fill=orange, minimum size=4pt](v6) at (6,0) {};
	\vertex[fill=orange, minimum size=4pt](v7) at (7,0) {};
	\draw[thick, -stealth](v0)--(v1);
	\draw[thick, -stealth](v1)--(v2);
	\draw[thick, -stealth](v2)--(v3);
	\draw[thick, -stealth](v3)--(v4);
	\draw[thick, -stealth](v4)--(v5);
	\draw[thick, -stealth](v5)--(v6);
	\draw[thick, -stealth](v6)--(v7);
    \draw[thick, -stealth] (v0) to [out=60,in=120] (v4);
    \draw[thick, -stealth] (v1) to [out=60,in=120] (v5);
    \draw[thick, -stealth] (v2) to [out=60,in=120] (v6);
	\draw[thick, -stealth] (v3) to [out=60,in=120] (v7);
    \node at (0,-0.5){\textcolor{red}{\footnotesize$1$}};
    \node at (1,-0.5){\textcolor{red}{\footnotesize$1$}};
	\node at (2,-0.5){\textcolor{red}{\footnotesize$1$}};
	\node at (3,-0.5){\textcolor{red}{\footnotesize$1$}};
	\node[rotate=-20] at (4.3,-0.5){\textcolor{red}{\tiny$1\!-\!u_2$}};
	\node[rotate=-20] at (5.3,-0.5){\textcolor{red}{\tiny$1\!-\!u_1$}};
	\node[rotate=-20] at (6.3,-0.5){\textcolor{red}{\tiny$1\!-\!s_4$}};
	\node[rotate=-20] at (7.7,-0.6){\textcolor{red}{\tiny$1\!-\!s_3\!-\!u_3\!$}};
\end{scope}
\end{tikzpicture}
    \caption{An example of the two interpretations of $E_{(s_1,s_2,s_3,s_4)}(u_1,u_2)$ in the graphs $G(4,10)$ and $G(4,8)$ where $ (s_1,s_2,s_3,s_4)\in T^4_{10}$, as described in Theorem~\ref{thm:kBoustrophedon}.}
    \label{fig:proof_kBoustrophedon}
\end{figure}
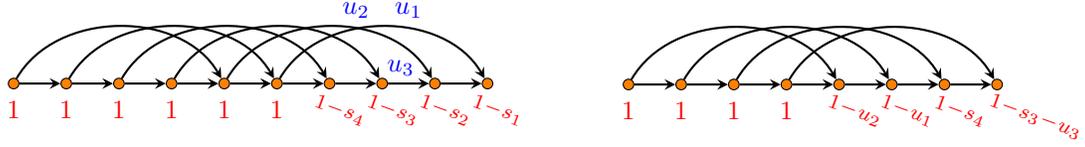

\begin{proof}
Interpret $E_\ss$ as the number of integer flows on $G(k,N)$ with net flow \[\vv=(1^{d-2k+1}, 1-s_k,\ldots, 1-s_1).\]
In the case $N=0$ so that $d=k-1$, then \[E_{(0,\ldots, 0)}= K_{G(k,0)}(1,\ldots, 1,-k+2)=1\] because $G(k,0)$ is the empty graph and there is only one integer flow.

Now let $N>0$ so that $d\geq k$, and fix $j\in[k-1]$.
Let us further refine $E_\ss$ according to the flows on the $j+1$ edges in the vertical edge cut between the vertices $N-j$ and $N-j+1$ in $G(k,N)$. Let $\widetilde{\mathcal{E}}_\ss$  be defined as in the bijection in Proposition~\ref{prop: other interpretation E_s}. Let $\widetilde{\mathcal{E}}_\ss(u_1,\ldots, u_{j+1})$ denote the subset of integer flows in $\widetilde{\mathcal{E}}_\ss$  whose flows on the last $j$ non-slack edges of $G(k,N)$ are $u_1,\ldots, u_j$ from right to left, and is $u_{j+1}$ on the slack edge between vertices $N-j$ and $N-j+1$.  (See Figure~\ref{fig:proof_kBoustrophedon} for an illustration in the case $k=4$, $d=5$, $j=2$.)
By examining the vertical edge cut between vertices $N-j$ and $N-j+1$, we see that the nonnegative integers $u_1,\ldots, u_j$ satisfy the inequalities
\[u_1+\cdots+u_h \leq s_1+\cdots+s_h -h\] for $h=1,\ldots, j$.
In addition, 
$$u_1+\cdots+u_{j+1}
= s_1+\cdots+s_j-j
= N-j - (s_{j+1} + \cdots+s_k).
$$
Let $E_\ss(u_1,\ldots, u_{j+1})$ be the cardinality of $\widetilde{\mathcal{E}}_\ss(u_1,\ldots, u_{j+1})$. 
The sets $\widetilde{\mathcal{E}}_\ss(u_1,\ldots, u_{j+1})$ partition $\widetilde{\mathcal{E}}_\ss$, so
$$E_\ss = \sum_{\mathbf{u}} E_\ss(u_1,\ldots, u_{j+1}). $$

As in the proof of Proposition~\ref{prop: other interpretation E_s}, the restriction of the graph $G(k, N)$ to the graph $G(k, N-j)$ gives a bijection between integer flows in $\widetilde{\mathcal{E}}_\ss(u_1,\ldots, u_{j+1})$ on $G(k, N)$ and integer flows on $G(k, N-j)$ with net flow 
$$(1^{d-2k+1-j}, 1-u_j, \ldots, 1-u_1, 1-s_k,\ldots, 1-s_{j+2}, 1-s_{j+1}-u_{j+1})$$
since the flow on each edge of $G(k,N)$ not in $G(k,N-j)$ is fixed by the choice of $(u_1,\ldots, u_{j+1})$. 
Thus we have shown that
$$E_\ss(u_1,\ldots, u_{j+1}) = K_{G(k, N-j)} (1^{d-2k+1-j}, 1-u_j, \ldots, 1-u_1, 1-s_k,\ldots, 1-s_{j+2}, 1-s_{j+1}-u_{j+1}),$$
which by Proposition~\ref{prop: other interpretation E_s} is counted by $E_{(u_{j+1}+s_{j+1}, u_{j+2}, \ldots, s_k, u_1, u_2,\ldots, u_j)}$, as desired.
\end{proof}

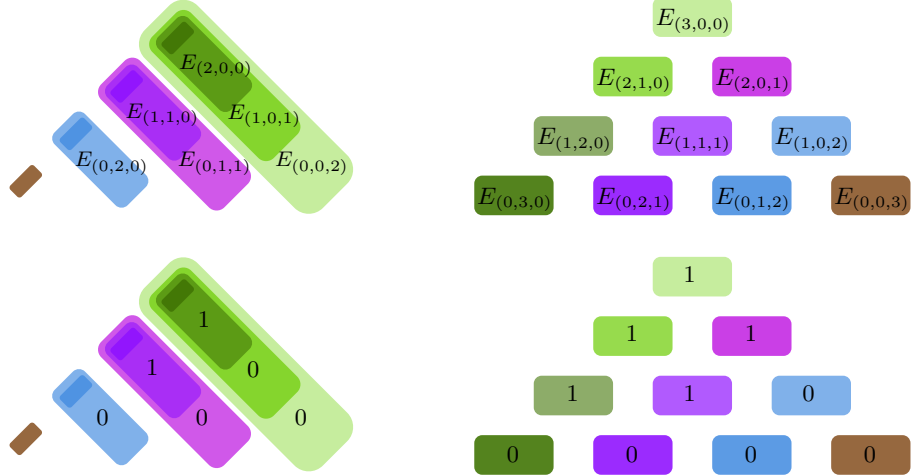
\begin{figure}
    \centering
\tikzset{every picture/.style={line width=0.75pt}} 
\begin{tikzpicture}[x=0.75pt,y=0.75pt,yscale=-1,xscale=1]

\draw  [draw opacity=0][fill={rgb, 255:red, 184; green, 233; blue, 134 }  ,fill opacity=0.8 ] (176.08,34.22) .. controls (180.05,30.25) and (186.5,30.25) .. (190.47,34.22) -- (265.93,109.68) .. controls (269.9,113.65) and (269.9,120.1) .. (265.93,124.07) -- (253.7,136.3) .. controls (249.73,140.27) and (243.28,140.27) .. (239.31,136.3) -- (163.85,60.84) .. controls (159.88,56.87) and (159.88,50.42) .. (163.85,46.45) -- cycle ;
\draw  [draw opacity=0][fill={rgb, 255:red, 126; green, 211; blue, 33 }  ,fill opacity=0.9 ] (179.42,38) .. controls (181.59,35.83) and (185.11,35.83) .. (187.28,38) -- (240.71,91.43) .. controls (242.88,93.6) and (242.88,97.12) .. (240.71,99.29) -- (228.93,111.07) .. controls (226.76,113.24) and (223.24,113.24) .. (221.07,111.07) -- (167.64,57.64) .. controls (165.47,55.47) and (165.47,51.95) .. (167.64,49.78) -- cycle ;
\draw  [draw opacity=0][fill={rgb, 255:red, 74; green, 144; blue, 226 }  ,fill opacity=0.9 ] (450,124) .. controls (450,121.79) and (451.79,120) .. (454,120) -- (486,120) .. controls (488.21,120) and (490,121.79) .. (490,124) -- (490,136) .. controls (490,138.21) and (488.21,140) .. (486,140) -- (454,140) .. controls (451.79,140) and (450,138.21) .. (450,136) -- cycle ;
\draw  [draw opacity=0][fill={rgb, 255:red, 126; green, 211; blue, 33 }  ,fill opacity=0.8 ] (390,64) .. controls (390,61.79) and (391.79,60) .. (394,60) -- (426,60) .. controls (428.21,60) and (430,61.79) .. (430,64) -- (430,76) .. controls (430,78.21) and (428.21,80) .. (426,80) -- (394,80) .. controls (391.79,80) and (390,78.21) .. (390,76) -- cycle ;
\draw  [draw opacity=0][fill={rgb, 255:red, 184; green, 233; blue, 134 }  ,fill opacity=0.8 ] (420,34) .. controls (420,31.79) and (421.79,30) .. (424,30) -- (456,30) .. controls (458.21,30) and (460,31.79) .. (460,34) -- (460,46) .. controls (460,48.21) and (458.21,50) .. (456,50) -- (424,50) .. controls (421.79,50) and (420,48.21) .. (420,46) -- cycle ;
\draw  [draw opacity=0][fill={rgb, 255:red, 65; green, 117; blue, 5 }  ,fill opacity=0.6 ] (360,94) .. controls (360,91.79) and (361.79,90) .. (364,90) -- (396,90) .. controls (398.21,90) and (400,91.79) .. (400,94) -- (400,106) .. controls (400,108.21) and (398.21,110) .. (396,110) -- (364,110) .. controls (361.79,110) and (360,108.21) .. (360,106) -- cycle ;
\draw  [draw opacity=0][fill={rgb, 255:red, 189; green, 16; blue, 224 }  ,fill opacity=0.7 ] (153.64,62.07) .. controls (155.81,59.91) and (159.33,59.91) .. (161.5,62.07) -- (215.78,116.36) .. controls (217.95,118.53) and (217.95,122.05) .. (215.78,124.22) -- (204,136) .. controls (201.83,138.17) and (198.31,138.17) .. (196.14,136) -- (141.86,81.71) .. controls (139.69,79.54) and (139.69,76.03) .. (141.86,73.86) -- cycle ;
\draw  [draw opacity=0][fill={rgb, 255:red, 74; green, 144; blue, 226 }  ,fill opacity=0.7 ] (127.57,88.93) .. controls (129.27,87.23) and (132.01,87.23) .. (133.7,88.93) -- (164.38,119.61) .. controls (166.08,121.3) and (166.08,124.04) .. (164.38,125.74) -- (155.19,134.94) .. controls (153.49,136.63) and (150.75,136.63) .. (149.05,134.94) -- (118.37,104.26) .. controls (116.68,102.56) and (116.68,99.82) .. (118.37,98.12) -- cycle ;
\draw  [draw opacity=0][fill={rgb, 255:red, 65; green, 117; blue, 5 }  ,fill opacity=0.9 ] (330,124) .. controls (330,121.79) and (331.79,120) .. (334,120) -- (366,120) .. controls (368.21,120) and (370,121.79) .. (370,124) -- (370,136) .. controls (370,138.21) and (368.21,140) .. (366,140) -- (334,140) .. controls (331.79,140) and (330,138.21) .. (330,136) -- cycle ;
\draw  [draw opacity=0][fill={rgb, 255:red, 144; green, 19; blue, 254 }  ,fill opacity=0.9 ] (390,124) .. controls (390,121.79) and (391.79,120) .. (394,120) -- (426,120) .. controls (428.21,120) and (430,121.79) .. (430,124) -- (430,136) .. controls (430,138.21) and (428.21,140) .. (426,140) -- (394,140) .. controls (391.79,140) and (390,138.21) .. (390,136) -- cycle ;
\draw  [draw opacity=0][fill={rgb, 255:red, 189; green, 16; blue, 224 }  ,fill opacity=0.8 ] (450,64) .. controls (450,61.79) and (451.79,60) .. (454,60) -- (486,60) .. controls (488.21,60) and (490,61.79) .. (490,64) -- (490,76) .. controls (490,78.21) and (488.21,80) .. (486,80) -- (454,80) .. controls (451.79,80) and (450,78.21) .. (450,76) -- cycle ;
\draw  [draw opacity=0][fill={rgb, 255:red, 144; green, 19; blue, 254 }  ,fill opacity=0.7 ] (420,94) .. controls (420,91.79) and (421.79,90) .. (424,90) -- (456,90) .. controls (458.21,90) and (460,91.79) .. (460,94) -- (460,106) .. controls (460,108.21) and (458.21,110) .. (456,110) -- (424,110) .. controls (421.79,110) and (420,108.21) .. (420,106) -- cycle ;
\draw  [draw opacity=0][fill={rgb, 255:red, 139; green, 87; blue, 42 }  ,fill opacity=0.9 ] (510,124) .. controls (510,121.79) and (511.79,120) .. (514,120) -- (546,120) .. controls (548.21,120) and (550,121.79) .. (550,124) -- (550,136) .. controls (550,138.21) and (548.21,140) .. (546,140) -- (514,140) .. controls (511.79,140) and (510,138.21) .. (510,136) -- cycle ;
\draw  [draw opacity=0][fill={rgb, 255:red, 74; green, 144; blue, 226 }  ,fill opacity=0.7 ] (480,94) .. controls (480,91.79) and (481.79,90) .. (484,90) -- (516,90) .. controls (518.21,90) and (520,91.79) .. (520,94) -- (520,106) .. controls (520,108.21) and (518.21,110) .. (516,110) -- (484,110) .. controls (481.79,110) and (480,108.21) .. (480,106) -- cycle ;
\draw  [draw opacity=0][fill={rgb, 255:red, 74; green, 144; blue, 226 }  ,fill opacity=0.9 ] (450,255) .. controls (450,252.79) and (451.79,251) .. (454,251) -- (486,251) .. controls (488.21,251) and (490,252.79) .. (490,255) -- (490,267) .. controls (490,269.21) and (488.21,271) .. (486,271) -- (454,271) .. controls (451.79,271) and (450,269.21) .. (450,267) -- cycle ;
\draw  [draw opacity=0][fill={rgb, 255:red, 126; green, 211; blue, 33 }  ,fill opacity=0.8 ] (390,195) .. controls (390,192.79) and (391.79,191) .. (394,191) -- (426,191) .. controls (428.21,191) and (430,192.79) .. (430,195) -- (430,207) .. controls (430,209.21) and (428.21,211) .. (426,211) -- (394,211) .. controls (391.79,211) and (390,209.21) .. (390,207) -- cycle ;
\draw  [draw opacity=0][fill={rgb, 255:red, 184; green, 233; blue, 134 }  ,fill opacity=0.8 ] (420,165) .. controls (420,162.79) and (421.79,161) .. (424,161) -- (456,161) .. controls (458.21,161) and (460,162.79) .. (460,165) -- (460,177) .. controls (460,179.21) and (458.21,181) .. (456,181) -- (424,181) .. controls (421.79,181) and (420,179.21) .. (420,177) -- cycle ;
\draw  [draw opacity=0][fill={rgb, 255:red, 65; green, 117; blue, 5 }  ,fill opacity=0.6 ] (360,225) .. controls (360,222.79) and (361.79,221) .. (364,221) -- (396,221) .. controls (398.21,221) and (400,222.79) .. (400,225) -- (400,237) .. controls (400,239.21) and (398.21,241) .. (396,241) -- (364,241) .. controls (361.79,241) and (360,239.21) .. (360,237) -- cycle ;
\draw  [draw opacity=0][fill={rgb, 255:red, 65; green, 117; blue, 5 }  ,fill opacity=0.9 ] (330,255) .. controls (330,252.79) and (331.79,251) .. (334,251) -- (366,251) .. controls (368.21,251) and (370,252.79) .. (370,255) -- (370,267) .. controls (370,269.21) and (368.21,271) .. (366,271) -- (334,271) .. controls (331.79,271) and (330,269.21) .. (330,267) -- cycle ;
\draw  [draw opacity=0][fill={rgb, 255:red, 144; green, 19; blue, 254 }  ,fill opacity=0.9 ] (390,255) .. controls (390,252.79) and (391.79,251) .. (394,251) -- (426,251) .. controls (428.21,251) and (430,252.79) .. (430,255) -- (430,267) .. controls (430,269.21) and (428.21,271) .. (426,271) -- (394,271) .. controls (391.79,271) and (390,269.21) .. (390,267) -- cycle ;
\draw  [draw opacity=0][fill={rgb, 255:red, 189; green, 16; blue, 224 }  ,fill opacity=0.8 ] (450,195) .. controls (450,192.79) and (451.79,191) .. (454,191) -- (486,191) .. controls (488.21,191) and (490,192.79) .. (490,195) -- (490,207) .. controls (490,209.21) and (488.21,211) .. (486,211) -- (454,211) .. controls (451.79,211) and (450,209.21) .. (450,207) -- cycle ;
\draw  [draw opacity=0][fill={rgb, 255:red, 144; green, 19; blue, 254 }  ,fill opacity=0.7 ] (420,225) .. controls (420,222.79) and (421.79,221) .. (424,221) -- (456,221) .. controls (458.21,221) and (460,222.79) .. (460,225) -- (460,237) .. controls (460,239.21) and (458.21,241) .. (456,241) -- (424,241) .. controls (421.79,241) and (420,239.21) .. (420,237) -- cycle ;
\draw  [draw opacity=0][fill={rgb, 255:red, 139; green, 87; blue, 42 }  ,fill opacity=0.9 ] (510,255) .. controls (510,252.79) and (511.79,251) .. (514,251) -- (546,251) .. controls (548.21,251) and (550,252.79) .. (550,255) -- (550,267) .. controls (550,269.21) and (548.21,271) .. (546,271) -- (514,271) .. controls (511.79,271) and (510,269.21) .. (510,267) -- cycle ;
\draw  [draw opacity=0][fill={rgb, 255:red, 74; green, 144; blue, 226 }  ,fill opacity=0.7 ] (480,225) .. controls (480,222.79) and (481.79,221) .. (484,221) -- (516,221) .. controls (518.21,221) and (520,222.79) .. (520,225) -- (520,237) .. controls (520,239.21) and (518.21,241) .. (516,241) -- (484,241) .. controls (481.79,241) and (480,239.21) .. (480,237) -- cycle ;
\draw  [draw opacity=0][fill={rgb, 255:red, 139; green, 87; blue, 42 }  ,fill opacity=0.9 ] (105.87,115.13) .. controls (106.49,114.51) and (107.51,114.51) .. (108.13,115.13) -- (111.54,118.54) .. controls (112.17,119.17) and (112.17,120.18) .. (111.54,120.81) -- (101.76,130.59) .. controls (101.13,131.22) and (100.11,131.22) .. (99.49,130.59) -- (96.08,127.19) .. controls (95.46,126.56) and (95.46,125.54) .. (96.08,124.92) -- cycle ;
\draw  [draw opacity=0][fill={rgb, 255:red, 144; green, 19; blue, 254 }  ,fill opacity=0.6 ] (154.19,65.19) .. controls (155.88,63.5) and (158.63,63.5) .. (160.32,65.19) -- (191,95.87) .. controls (192.69,97.56) and (192.69,100.31) .. (191,102) -- (181.8,111.2) .. controls (180.11,112.89) and (177.36,112.89) .. (175.67,111.2) -- (144.99,80.52) .. controls (143.3,78.83) and (143.3,76.08) .. (144.99,74.39) -- cycle ;
\draw  [draw opacity=0][fill={rgb, 255:red, 65; green, 117; blue, 5 }  ,fill opacity=0.6 ] (179.9,40.65) .. controls (181.59,38.96) and (184.33,38.96) .. (186.03,40.65) -- (216.71,71.33) .. controls (218.4,73.03) and (218.4,75.77) .. (216.71,77.46) -- (207.51,86.66) .. controls (205.82,88.36) and (203.07,88.36) .. (201.38,86.66) -- (170.7,55.98) .. controls (169,54.29) and (169,51.54) .. (170.7,49.85) -- cycle ;
\draw  [draw opacity=0][fill={rgb, 255:red, 65; green, 117; blue, 5 }  ,fill opacity=0.9 ] (182.49,42.82) .. controls (183.12,42.19) and (184.14,42.19) .. (184.76,42.82) -- (188.17,46.22) .. controls (188.8,46.85) and (188.8,47.87) .. (188.17,48.49) -- (178.39,58.28) .. controls (177.76,58.9) and (176.74,58.9) .. (176.12,58.28) -- (172.71,54.87) .. controls (172.08,54.25) and (172.08,53.23) .. (172.71,52.6) -- cycle ;
\draw  [draw opacity=0][fill={rgb, 255:red, 144; green, 19; blue, 254 }  ,fill opacity=0.9 ] (156.75,67.27) .. controls (157.38,66.64) and (158.4,66.64) .. (159.02,67.27) -- (162.43,70.68) .. controls (163.05,71.3) and (163.05,72.32) .. (162.43,72.94) -- (152.64,82.73) .. controls (152.02,83.35) and (151,83.35) .. (150.38,82.73) -- (146.97,79.32) .. controls (146.34,78.7) and (146.34,77.68) .. (146.97,77.05) -- cycle ;
\draw  [draw opacity=0][fill={rgb, 255:red, 74; green, 144; blue, 226 }  ,fill opacity=0.9 ] (131.16,91.2) .. controls (131.78,90.58) and (132.8,90.58) .. (133.43,91.2) -- (136.83,94.61) .. controls (137.46,95.23) and (137.46,96.25) .. (136.83,96.88) -- (127.05,106.66) .. controls (126.42,107.29) and (125.41,107.29) .. (124.78,106.66) -- (121.37,103.26) .. controls (120.75,102.63) and (120.75,101.61) .. (121.37,100.99) -- cycle ;
\draw  [draw opacity=0][fill={rgb, 255:red, 184; green, 233; blue, 134 }  ,fill opacity=0.8 ] (176.08,164.22) .. controls (180.05,160.25) and (186.5,160.25) .. (190.47,164.22) -- (265.93,239.68) .. controls (269.9,243.65) and (269.9,250.1) .. (265.93,254.07) -- (253.7,266.3) .. controls (249.73,270.27) and (243.28,270.27) .. (239.31,266.3) -- (163.85,190.84) .. controls (159.88,186.87) and (159.88,180.42) .. (163.85,176.45) -- cycle ;
\draw  [draw opacity=0][fill={rgb, 255:red, 126; green, 211; blue, 33 }  ,fill opacity=0.9 ] (179.42,168) .. controls (181.59,165.83) and (185.11,165.83) .. (187.28,168) -- (240.71,221.43) .. controls (242.88,223.6) and (242.88,227.12) .. (240.71,229.29) -- (228.93,241.07) .. controls (226.76,243.24) and (223.24,243.24) .. (221.07,241.07) -- (167.64,187.64) .. controls (165.47,185.47) and (165.47,181.95) .. (167.64,179.78) -- cycle ;
\draw  [draw opacity=0][fill={rgb, 255:red, 189; green, 16; blue, 224 }  ,fill opacity=0.7 ] (153.64,192.07) .. controls (155.81,189.91) and (159.33,189.91) .. (161.5,192.07) -- (215.78,246.36) .. controls (217.95,248.53) and (217.95,252.05) .. (215.78,254.22) -- (204,266) .. controls (201.83,268.17) and (198.31,268.17) .. (196.14,266) -- (141.86,211.71) .. controls (139.69,209.54) and (139.69,206.03) .. (141.86,203.86) -- cycle ;
\draw  [draw opacity=0][fill={rgb, 255:red, 74; green, 144; blue, 226 }  ,fill opacity=0.7 ] (127.57,218.93) .. controls (129.27,217.23) and (132.01,217.23) .. (133.7,218.93) -- (164.38,249.61) .. controls (166.08,251.3) and (166.08,254.04) .. (164.38,255.74) -- (155.19,264.94) .. controls (153.49,266.63) and (150.75,266.63) .. (149.05,264.94) -- (118.37,234.26) .. controls (116.68,232.56) and (116.68,229.82) .. (118.37,228.12) -- cycle ;
\draw  [draw opacity=0][fill={rgb, 255:red, 139; green, 87; blue, 42 }  ,fill opacity=0.9 ] (105.87,245.13) .. controls (106.49,244.51) and (107.51,244.51) .. (108.13,245.13) -- (111.54,248.54) .. controls (112.17,249.17) and (112.17,250.18) .. (111.54,250.81) -- (101.76,260.59) .. controls (101.13,261.22) and (100.11,261.22) .. (99.49,260.59) -- (96.08,257.19) .. controls (95.46,256.56) and (95.46,255.54) .. (96.08,254.92) -- cycle ;
\draw  [draw opacity=0][fill={rgb, 255:red, 144; green, 19; blue, 254 }  ,fill opacity=0.6 ] (154.19,195.19) .. controls (155.88,193.5) and (158.63,193.5) .. (160.32,195.19) -- (191,225.87) .. controls (192.69,227.56) and (192.69,230.31) .. (191,232) -- (181.8,241.2) .. controls (180.11,242.89) and (177.36,242.89) .. (175.67,241.2) -- (144.99,210.52) .. controls (143.3,208.83) and (143.3,206.08) .. (144.99,204.39) -- cycle ;
\draw  [draw opacity=0][fill={rgb, 255:red, 65; green, 117; blue, 5 }  ,fill opacity=0.6 ] (179.9,170.65) .. controls (181.59,168.96) and (184.33,168.96) .. (186.03,170.65) -- (216.71,201.33) .. controls (218.4,203.03) and (218.4,205.77) .. (216.71,207.46) -- (207.51,216.66) .. controls (205.82,218.36) and (203.07,218.36) .. (201.38,216.66) -- (170.7,185.98) .. controls (169,184.29) and (169,181.54) .. (170.7,179.85) -- cycle ;
\draw  [draw opacity=0][fill={rgb, 255:red, 65; green, 117; blue, 5 }  ,fill opacity=0.9 ] (182.49,172.82) .. controls (183.12,172.19) and (184.14,172.19) .. (184.76,172.82) -- (188.17,176.22) .. controls (188.8,176.85) and (188.8,177.87) .. (188.17,178.49) -- (178.39,188.28) .. controls (177.76,188.9) and (176.74,188.9) .. (176.12,188.28) -- (172.71,184.87) .. controls (172.08,184.25) and (172.08,183.23) .. (172.71,182.6) -- cycle ;
\draw  [draw opacity=0][fill={rgb, 255:red, 144; green, 19; blue, 254 }  ,fill opacity=0.9 ] (156.75,197.27) .. controls (157.38,196.64) and (158.4,196.64) .. (159.02,197.27) -- (162.43,200.68) .. controls (163.05,201.3) and (163.05,202.32) .. (162.43,202.94) -- (152.64,212.73) .. controls (152.02,213.35) and (151,213.35) .. (150.38,212.73) -- (146.97,209.32) .. controls (146.34,208.7) and (146.34,207.68) .. (146.97,207.05) -- cycle ;
\draw  [draw opacity=0][fill={rgb, 255:red, 74; green, 144; blue, 226 }  ,fill opacity=0.9 ] (131.16,221.2) .. controls (131.78,220.58) and (132.8,220.58) .. (133.43,221.2) -- (136.83,224.61) .. controls (137.46,225.23) and (137.46,226.25) .. (136.83,226.88) -- (127.05,236.66) .. controls (126.42,237.29) and (125.41,237.29) .. (124.78,236.66) -- (121.37,233.26) .. controls (120.75,232.63) and (120.75,231.61) .. (121.37,230.99) -- cycle ;

\draw (390.03,124.4) node [anchor=north west][inner sep=0.75pt]  [font=\scriptsize]  {$E_{( 0,2,1)}$};
\draw (330.01,124.4) node [anchor=north west][inner sep=0.75pt]  [font=\scriptsize]  {$E_{( 0,3,0)}$};
\draw (450,124.4) node [anchor=north west][inner sep=0.75pt]  [font=\scriptsize]  {$E_{( 0,1,2)}$};
\draw (510,124.4) node [anchor=north west][inner sep=0.75pt]  [font=\scriptsize]  {$E_{( 0,0,3)}$};
\draw (360,94.4) node [anchor=north west][inner sep=0.75pt]  [font=\scriptsize]  {$E_{( 1,2,0)}$};
\draw (420,94.4) node [anchor=north west][inner sep=0.75pt]  [font=\scriptsize]  {$E_{( 1,1,1)}$};
\draw (480,94.4) node [anchor=north west][inner sep=0.75pt]  [font=\scriptsize]  {$E_{( 1,0,2)}$};
\draw (390,64.4) node [anchor=north west][inner sep=0.75pt]  [font=\scriptsize]  {$E_{( 2,1,0)}$};
\draw (450,64.4) node [anchor=north west][inner sep=0.75pt]  [font=\scriptsize]  {$E_{( 2,0,1)}$};
\draw (420,34.4) node [anchor=north west][inner sep=0.75pt]  [font=\scriptsize]  {$E_{( 3,0,0)}$};
\draw (127,106.4) node [anchor=north west][inner sep=0.75pt]  [font=\tiny]  {$E_{( 0,2,0)}$};
\draw (179,106.4) node [anchor=north west][inner sep=0.75pt]  [font=\tiny]  {$E_{( 0,1,1)}$};
\draw (228,106.4) node [anchor=north west][inner sep=0.75pt]  [font=\tiny]  {$E_{( 0,0,2)}$};
\draw (152.64,82.13) node [anchor=north west][inner sep=0.75pt]  [font=\tiny]  {$E_{( 1,1,0)}$};
\draw (204,82.4) node [anchor=north west][inner sep=0.75pt]  [font=\tiny]  {$E_{( 1,0,1)}$};
\draw (179.64,57.18) node [anchor=north west][inner sep=0.75pt]  [font=\tiny]  {$E_{( 2,0,0)}$};
\draw (525,255.4) node [anchor=north west][inner sep=0.75pt]  [font=\scriptsize]  {$0$};
\draw (375,224.4) node [anchor=north west][inner sep=0.75pt]  [font=\scriptsize]  {$1$};
\draw (434,164.4) node [anchor=north west][inner sep=0.75pt]  [font=\scriptsize]  {$1$};
\draw (405,194.4) node [anchor=north west][inner sep=0.75pt]  [font=\scriptsize]  {$1$};
\draw (465,194.4) node [anchor=north west][inner sep=0.75pt]  [font=\scriptsize]  {$1$};
\draw (434,224.4) node [anchor=north west][inner sep=0.75pt]  [font=\scriptsize]  {$1$};
\draw (345,255.4) node [anchor=north west][inner sep=0.75pt]  [font=\scriptsize]  {$0$};
\draw (405,255.4) node [anchor=north west][inner sep=0.75pt]  [font=\scriptsize]  {$0$};
\draw (494,224.4) node [anchor=north west][inner sep=0.75pt]  [font=\scriptsize]  {$0$};
\draw (465,255.4) node [anchor=north west][inner sep=0.75pt]  [font=\scriptsize]  {$0$};
\draw (189,187.4) node [anchor=north west][inner sep=0.75pt]  [font=\scriptsize]  {$1$};
\draw (163,211.4) node [anchor=north west][inner sep=0.75pt]  [font=\scriptsize]  {$1$};
\draw (238,236.4) node [anchor=north west][inner sep=0.75pt]  [font=\scriptsize]  {$0$};
\draw (188,236.4) node [anchor=north west][inner sep=0.75pt]  [font=\scriptsize]  {$0$};
\draw (138,236.4) node [anchor=north west][inner sep=0.75pt]  [font=\scriptsize]  {$0$};
\draw (214,212.4) node [anchor=north west][inner sep=0.75pt]  [font=\scriptsize]  {$0$};
\end{tikzpicture}
    \caption{Boustrophedon recursion for $k$-Entringer numbers.}
    \label{fig:boustrophedon_recursion}
\end{figure}

The following corollary is the special case of the $k$-boustrophedon recursion when $j=1$. It recovers the boustrophedon recurrence of Ayyer, Josuat-Verg\`es and Ramassamy~\cite[Theorem 7.4]{AyyerJosuatvergesRamassamy2018}. See Figure~\ref{fig:boustrophedon_recursion} for an illustration.
Alternatively, Theorem~\ref{thm:kBoustrophedon} can be derived by applying Corollary~\ref{thm:boustrophedon_recurrence} $j$ times.
\begin{corollary}[{\cite[Theorem 7.4]{AyyerJosuatvergesRamassamy2018}}]
\label{thm:boustrophedon_recurrence}
Let $k\geq2$, $N=d-k+1\geq1$, and $(s_1,\ldots,s_k) \in T^k_N$. Then  
\[
E_{(s_1,\ldots,s_k)} = 
\begin{cases}
1 & \hbox{if } (s_1,\ldots,s_k)=(0,\ldots,0),\\
	\sum_{t=0}^{s_1-1}  E_{(s_2+t,s_3,\ldots,s_k,s_1-t-1)}, &\hbox{if } s_1>0,\\
	0, & \hbox{otherwise}.
\end{cases}
\]
\end{corollary}

\begin{corollary}
When $\ss=(d,0,\ldots,0)$, the $k$-Entringer number $E_\ss$ is also a $k$-Euler number:
$$E_{(d,0,\ldots,0)} = \sum_{\ss\in T_N^k} E_\ss = A_{k,d},$$
where $N=d-k+1$.
Thus the $k$-Entringer number $E_{(d,0,\ldots,0)}$ is also the volume of the flow polytope $\F_{G(k,d+k)}$.
\end{corollary}
\begin{proof}
Applying Theorem~\ref{thm:kBoustrophedon} with $s=(d, 0,\ldots,0)$ and $j=k-1$ then 
$$E_{(d,0,\ldots,0)} = \sum_{\mathbf{u}} E_{(u_k+0, u_1,\ldots, u_{k-1})} $$
is a sum over all compositions $\mathbf{u}=(u_k,u_1,\ldots, u_{k-1})$ such that $u_1+\cdots+ u_k = d-k+1 = N$.
By Equation~\eqref{eqn.genEnt} we have that $A_{k,d} = \sum_{\ss \in T_N^k}E_\ss $.
Putting these equations together gives the desired result.
\end{proof}

\begin{example}\label{eg.simplices}
Let $k=3$ and $N=d-k+1$. We present the $3$-Entringer numbers of $\ss\in T_N^3$ for $N=3,4,5$:

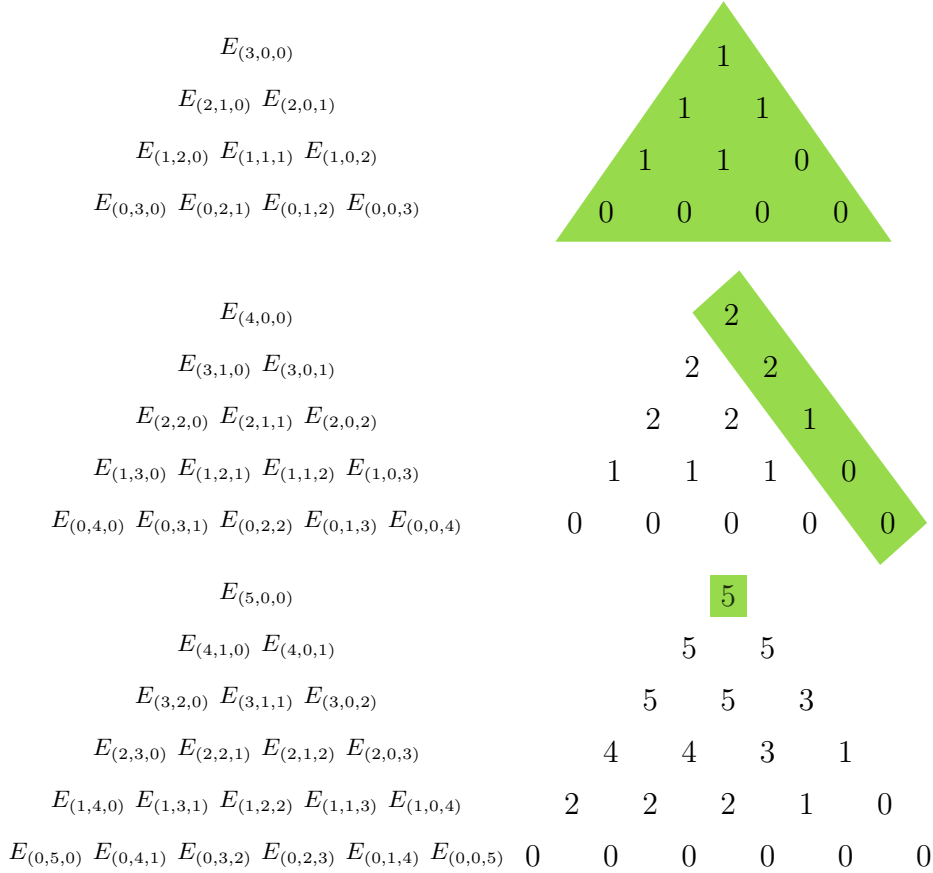
\begin{center}
    \begin{tikzpicture}[scale=0.4]
\begin{scope}[xshift=60, xscale=1.4]	
\node at (0,0){\scriptsize$E_{(0,3,0)}$};
\node at (2,0){\scriptsize$E_{(0,2,1)}$};
\node at (4,0){\scriptsize$E_{(0,1,2)}$};
\node at (6,0){\scriptsize$E_{(0,0,3)}$};
\node at (1,1.73){\scriptsize$E_{(1,2,0)}$};
\node at (3,1.73){\scriptsize$E_{(1,1,1)}$};
\node at (5,1.73){\scriptsize$E_{(1,0,2)}$};
\node at (2,3.46){\scriptsize$E_{(2,1,0)}$};
\node at (4,3.46){\scriptsize$E_{(2,0,1)}$};
\node at (3,5.19){\scriptsize$E_{(3,0,0)}$};
\end{scope}

\begin{scope}[xshift=380, xscale=1.3, shift={(3.5,-0.19)}]	
\draw[draw opacity=0, fill={rgb, 255:red, 126; green, 211; blue, 33 }  ,fill opacity=0.8 ]  (-1.3,-1)--(7.3,-1)--(3,7)--cycle;
\node at (0,0) {$0$};
\node at (2,0) {$0$};
\node at (4,0) {$0$};
\node at (6,0) {$0$};
\node at (1,1.73) {$1$};
\node at (3,1.73) {$1$};
\node at (5,1.73) {$0$};
\node at (2,3.46) {$1$};
\node at (4,3.46) {$1$};
\node at (3,5.19) {$1$};
\end{scope}

\begin{scope}[xshift=20,yshift=-300, xscale=1.4]	
\node at (0,0){\scriptsize$E_{(0,4,0)}$};
\node at (2,0){\scriptsize$E_{(0,3,1)}$};
\node at (4,0){\scriptsize$E_{(0,2,2)}$};
\node at (6,0){\scriptsize$E_{(0,1,3)}$};
\node at (8,0){\scriptsize$E_{(0,0,4)}$};
\node at (1,1.73){\scriptsize$E_{(1,3,0)}$};
\node at (3,1.73){\scriptsize$E_{(1,2,1)}$};
\node at (5,1.73){\scriptsize$E_{(1,1,2)}$};
\node at (7,1.73){\scriptsize$E_{(1,0,3)}$};
\node at (2,3.46){\scriptsize$E_{(2,2,0)}$};
\node at (4,3.46){\scriptsize$E_{(2,1,1)}$};
\node at (6,3.46){\scriptsize$E_{(2,0,2)}$};
\node at (3,5.19){\scriptsize$E_{(3,1,0)}$};
\node at (5,5.19){\scriptsize$E_{(3,0,1)}$};
\node at (4,6.92){\scriptsize$E_{(4,0,0)}$};
\end{scope}

\begin{scope}[xshift=480,yshift=-300, xscale=1.3]	

\draw[draw opacity=0, fill={rgb, 255:red, 126; green, 211; blue, 33 }  ,fill opacity=0.8 ]  (3,7)--(7.8,-1.4)--(9,0)--(4.2,8.4)--cycle;
\node at (0,0){$0$};
\node at (2,0){$0$};
\node at (4,0){$0$};
\node at (6,0){$0$};
\node at (8,0){$0$};
\node at (1,1.73){$1$};
\node at (3,1.73){$1$};
\node at (5,1.73){$1$};
\node at (7,1.73){$0$};
\node at (2,3.46){$2$};
\node at (4,3.46){$2$};
\node at (6,3.46){$1$};
\node at (3,5.19){$2$};
\node at (5,5.19){$2$};
\node at (4,6.92){$2$};
\end{scope}

\begin{scope}[xshift=0,yshift=-600, xscale=1.4, shift={(-0.5,-0.5)}]	
\node at (0,0) {\scriptsize$E_{(0,5,0)}$};
\node at (2,0) {\scriptsize$E_{(0,4,1)}$};
\node at (4,0) {\scriptsize$E_{(0,3,2)}$};
\node at (6,0) {\scriptsize$E_{(0,2,3)}$};
\node at (8,0) {\scriptsize$E_{(0,1,4)}$};
\node at (10,0) {\scriptsize$E_{(0,0,5)}$};
\node at (1,1.73) {\scriptsize$E_{(1,4,0)}$};
\node at (3,1.73) {\scriptsize$E_{(1,3,1)}$};
\node at (5,1.73) {\scriptsize$E_{(1,2,2)}$};
\node at (7,1.73) {\scriptsize$E_{(1,1,3)}$};
\node at (9,1.73) {\scriptsize$E_{(1,0,4)}$};
\node at (2,3.46) {\scriptsize$E_{(2,3,0)}$};
\node at (4,3.46) {\scriptsize$E_{(2,2,1)}$};
\node at (6,3.46) {\scriptsize$E_{(2,1,2)}$};
\node at (8,3.46) {\scriptsize$E_{(2,0,3)}$};
\node at (3,5.19) {\scriptsize$E_{(3,2,0)}$};
\node at (5,5.19) {\scriptsize$E_{(3,1,1)}$};
\node at (7,5.19) {\scriptsize$E_{(3,0,2)}$};
\node at (4,6.92) {\scriptsize$E_{(4,1,0)}$};
\node at (6,6.92) {\scriptsize$E_{(4,0,1)}$};
\node at (5,8.65) {\scriptsize$E_{(5,0,0)}$};
\end{scope}

\begin{scope}[xshift=440,yshift=-615, xscale=1.3]	
\node at (0,0) {$0$};
\node at (2,0) {$0$};
\node at (4,0) {$0$};
\node at (6,0) {$0$};
\node at (8,0) {$0$};
\node at (10,0) {$0$};
\node at (1,1.73) {$2$};
\node at (3,1.73) {$2$};
\node at (5,1.73) {$2$};
\node at (7,1.73) {$1$};
\node at (9,1.73) {$0$};
\node at (2,3.46) {$4$};
\node at (4,3.46) {$4$};
\node at (6,3.46) {$3$};
\node at (8,3.46) {$1$};
\node at (3,5.19) {$5$};
\node at (5,5.19) {$5$};
\node at (7,5.19) {$3$};
\node at (4,6.92) {$5$};
\node at (6,6.92) {$5$};
\node[fill={rgb, 255:red, 126; green, 211; blue, 33 },fill opacity=0.8] at (5,8.65) {$5$};
\end{scope}
\end{tikzpicture}
\end{center}

By Corollary~\ref{cor. k-Euler is kpf}, for $\ss=(s_1,s_2,s_3)\in T_N^3$, the $3$-Entringer numbers $E_\ss$ enumerate integral $(1^{d-1},-d+1)$-flows on $G(3,d)$ whose flows on the last two non-slack edges $(d-1,d+2)$ and $(d,d+3)$ are $s_3$ and $s_2$, and whose flow on the slack edge $(d-2,d-1)$ is $s_1$.

By Theorem~\ref{thm:kBoustrophedon}, we can express any Entringer number indexed by $\ss\in T_5^3$ as partial sums of $3$-Entringer numbers indexed by entries $T_4^3$ and by $T_3^3$. For example, 
$$5=E_{(5,0,0)} 
    = \sum_{\mathbf{s}\in T_4^3: s_2=0} E_{\mathbf{s}}
    = \sum_{\mathbf{t}\in T_3^3} E_{\mathbf{t}}.
$$
We see that $E_{(5,0,0)}=5$ is simultaneously
\begin{enumerate}
    \item[(i)] the $3$-Entringer number indexed by the top entry of $T_5^3$,
    \item[(ii)] the sum of $3$-Entringer numbers indexed by entries along the right edge of $T_4^3$,
    \item[(iii)] and the sum of $3$-Entringer numbers indexed by all the entries of $T_3^3$.
\end{enumerate}
This verifies the result of Equation~\eqref{eqn.genEnt} which states that $k$-Euler numbers are refined by $k$-Entringer numbers 
$$\vol\F_{G(3,8)} = A_{3,5} = \sum_{\ss\in T_3^3} E_\ss = 5.$$
\end{example}


\subsection{Log-concavity of the \texorpdfstring{$k$}{}-Entringer numbers}
\label{sec:logconcavity}

We can use the machinery of flow polytopes to study log-concavity properties of the $k$-Entringer numbers. In the case of $k=2$, Benedetti et al.~{\cite[Corollary 7.6]{BGHHKMY}} already proved  that the sequence $E_{(0,N)},\ldots, E_{(N,N)}$ of Entringer numbers is log-concave.

The following log-concavity result follows from the general Lidskii volume formula, Theorem~\ref{thm.genlidskii}. In this formula, the numbers $K_G(\ss-\tt)$ are \emph{mixed volumes} (see~\cite[Section 3.4]{BV} for example). The following result is then a consequence of the  Aleksandrov--Fenchel inequalities~\cite{Alexandrov, Fe2, Fe1}. Alternatively, the result also follows from work of Huh et al.~\cite[Proposition 11]{HMMStD} on \emph{Lorentzian polynomials}. See Section~\ref{subsec.Lorentzian}.

\begin{lemma}[{Huh et al.~\cite[Proposition 11]{HMMStD}}] \label{log-concavity-kpf}
The numbers $K_G(\ss-\tt)$ appearing in Equation~\eqref{eq: lidskii outdegree} are log-concave along root directions.
That is, 
$$K_G(\ss-\tt)^2 \geq K_G(\ss-\tt - \ee_i + \ee_{j})\cdot K_G(\ss-\tt + \ee_i-\ee_{j})$$
for each $1\leq i<j \leq n$. 
\end{lemma}

\logconcavitythm

\begin{figure}
    \centering
\begin{tikzpicture}
\begin{scope}[xshift=0, yshift=0, scale=0.6]
	\vertex[fill=orange, minimum size=4pt](v1) at (1,0) {};
	\vertex[fill=orange, minimum size=4pt](v2) at (2,0) {};
	\vertex[fill=orange, minimum size=4pt](v3) at (3,0) {};
	\vertex[fill=orange, minimum size=4pt](v4) at (4,0) {};
	\vertex[fill=orange, minimum size=4pt](v5) at (5,0) {};
	\vertex[fill=orange, minimum size=4pt](v6) at (6,0) {};
	\vertex[fill=orange, minimum size=4pt](v7) at (7,0) {};
	\vertex[fill=orange, minimum size=4pt](v8) at (8,0) {};
	\vertex[fill=orange, minimum size=4pt](v9) at (9,0) {};
	\vertex[fill=orange, minimum size=4pt](v10) at (10,0) {};
	\vertex[fill=orange, minimum size=4pt](v11) at (11,0) {};
	\vertex[fill=orange, minimum size=4pt](v12) at (12,0) {};
	\vertex[fill=orange, minimum size=4pt](v13) at (13,0) {};
	\draw[thick, -stealth](v1)--(v2);
	\draw[thick, -stealth](v2)--(v3);
	\draw[thick, -stealth](v3)--(v4);
	\draw[thick, -stealth](v4)--(v5);
	\draw[thick, -stealth](v5)--(v6);
	\draw[thick, -stealth](v6)--(v7);
	\draw[thick, -stealth](v7)--(v8);
	\draw[thick, -stealth](v8)--(v9);
	\draw[thick, -stealth](v9)--(v10);
	\draw[thick, -stealth](v10)--(v11);
	\draw[thick, -stealth](v11)--(v12);
	\draw[thick, -stealth](v12)--(v13);
	\draw[thick, -stealth] (v1) to [out=60,in=120] (v5);
	\draw[thick, -stealth] (v2) to [out=60,in=120] (v6);
	\draw[thick, -stealth] (v3) to [out=60,in=120] (v7);
	\draw[thick, -stealth] (v4) to [out=60,in=120] (v8);
	\draw[thick, -stealth] (v5) to [out=60,in=120] (v9);
	\draw[thick, -stealth] (v6) to [out=60,in=120] (v10);
	\draw[thick, -stealth] (v7) to [out=60,in=120] (v11);
	\draw[thick, -stealth] (v8) to [out=60,in=120] (v12);
	\draw[thick, -stealth] (v9) to [out=60,in=120] (v13);
	\node at (1,-0.5){\textcolor{red}{\footnotesize$1$}};
	\node at (2,-0.5){\textcolor{red}{\footnotesize$1$}};
	\node at (3,-0.5){\textcolor{red}{\footnotesize$1$}};
	\node at (4,-0.5){\textcolor{red}{\footnotesize$1$}};
	\node at (5,-0.5){\textcolor{red}{\footnotesize$0$}};
	\node[rotate=0] at (6,-0.5){\textcolor{red}{\footnotesize$0$}};
	\node[rotate=0] at (7,-0.5){\textcolor{red}{\footnotesize$0$}};
	\node[rotate=0] at (8,-0.5){\textcolor{red}{\footnotesize$0$}};
	\node[rotate=0] at (9,-0.5){\textcolor{red}{\footnotesize$0$}};
	\node[rotate=0] at (10,-0.5){\textcolor{red}{\footnotesize$0$}};
	\node[rotate=0] at (11,-0.5){\textcolor{red}{\footnotesize$0$}};
	\node[rotate=0] at (12,-0.5){\textcolor{red}{\footnotesize$0$}};
	\node[rotate=0] at (13,-0.5){\textcolor{red}{\footnotesize$-4$}};
\end{scope}
\begin{scope}[xshift=250, yshift=0, scale=0.6]
	\vertex[fill=orange, minimum size=4pt](v1) at (1,0) {};
	\vertex[fill=orange, minimum size=4pt](v2) at (2,0) {};
	\vertex[fill=orange, minimum size=4pt](v3) at (3,0) {};
	\vertex[fill=orange, minimum size=4pt](v4) at (4,0) {};
	\vertex[fill=orange, minimum size=4pt](v5) at (5,0) {};
	\vertex[fill=orange, minimum size=4pt](v6) at (6,0) {};
	\vertex[fill=orange, minimum size=4pt](v7) at (7,0) {};
	\vertex[fill=orange, minimum size=4pt](v8) at (8,0) {};
	\vertex[fill=orange, minimum size=4pt](v9) at (9,0) {};
	\vertex[fill=orange, minimum size=4pt](v10) at (10,0) {};
	\vertex[fill=orange, minimum size=4pt](v11) at (11,0) {};
	\vertex[fill=orange, minimum size=4pt](v12) at (12,0) {};
	\vertex[fill=orange, minimum size=4pt](v13) at (13,0) {};	
	\draw[thick, -stealth](v1)--(v2);
	\draw[thick, -stealth](v2)--(v3);
	\draw[thick, -stealth](v3)--(v4);
	\draw[thick, -stealth](v4)--(v5);
	\draw[thick, -stealth](v5)--(v6);
	\draw[thick, -stealth](v6)--(v7);
	\draw[thick, -stealth](v7)--(v8);
	\draw[thick, -stealth](v8)--(v9);
	\draw[thick, -stealth](v9)--(v10);
	\draw[thick, -stealth](v10)--(v11);
	\draw[thick, -stealth](v11)--(v12);
	\draw[thick, -stealth](v12)--(v13);	
	\draw[thick, -stealth] (v1) to [out=60,in=120] (v5);
	\draw[thick, -stealth] (v2) to [out=60,in=120] (v6);
	\draw[thick, -stealth] (v3) to [out=60,in=120] (v7);
	\draw[thick, -stealth] (v4) to [out=60,in=120] (v8);
	\draw[thick, -stealth] (v5) to [out=60,in=120] (v9);
	\draw[thick, -stealth] (v6) to [out=60,in=120] (v10);
	\draw[thick, -stealth] (v7) to [out=60,in=120] (v11);
	\draw[thick, -stealth] (v8) to [out=60,in=120] (v12);
	\draw[thick, -stealth] (v9) to [out=60,in=120] (v13);
    \node[rotate=20] at (0.6,-0.5){\textcolor{red}{\tiny$s_1\!-\!1$}};
	\node[rotate=20] at (1.6,-0.5){\textcolor{red}{\tiny$s_2\!-\!1$}};
	\node[rotate=20] at (2.6,-0.5){\textcolor{red}{\tiny$s_3\!-\!1$}};
	\node[rotate=20] at (3.6,-0.5){\textcolor{red}{\tiny$s_4\!-\!1$}};
	\node at (5,-0.5){\textcolor{red}{\footnotesize$1$}};
	\node at (6,-0.5){\textcolor{red}{\footnotesize$1$}};
	\node at (7,-0.5){\textcolor{red}{\footnotesize$1$}};
	\node at (8,-0.5){\textcolor{red}{\footnotesize$1$}};
	\node at (9,-0.5){\textcolor{red}{\footnotesize$1$}};
	\node at (10,-0.5){\textcolor{red}{\footnotesize$0$}};
	\node at (11,-0.5){\textcolor{red}{\footnotesize$0$}};
	\node at (12,-0.5){\textcolor{red}{\footnotesize$0$}};
	\node at (13,-0.5){\textcolor{red}{\footnotesize$0$}};	
\end{scope}
\begin{scope}[xshift=250, yshift=-50, scale=0.6]
	\vertex[fill=orange, minimum size=4pt](v1) at (1,0) {};
	\vertex[fill=orange, minimum size=4pt](v2) at (2,0) {};
	\vertex[fill=orange, minimum size=4pt](v3) at (3,0) {};
	\vertex[fill=orange, minimum size=4pt](v4) at (4,0) {};
	\vertex[fill=orange, minimum size=4pt](v5) at (5,0) {};
	\vertex[fill=orange, minimum size=4pt](v6) at (6,0) {};
	\vertex[fill=orange, minimum size=4pt](v7) at (7,0) {};
	\vertex[fill=orange, minimum size=4pt](v8) at (8,0) {};
	\vertex[fill=orange, minimum size=4pt](v9) at (9,0) {};
	\draw[thick, -stealth](v1)--(v2);
	\draw[thick, -stealth](v2)--(v3);
	\draw[thick, -stealth](v3)--(v4);
	\draw[thick, -stealth](v4)--(v5);
	\draw[thick, -stealth](v5)--(v6);
	\draw[thick, -stealth](v6)--(v7);
	\draw[thick, -stealth](v7)--(v8);
	\draw[thick, -stealth](v8)--(v9);
	\draw[thick, -stealth] (v1) to [out=60,in=120] (v5);
	\draw[thick, -stealth] (v2) to [out=60,in=120] (v6);
	\draw[thick, -stealth] (v3) to [out=60,in=120] (v7);
	\draw[thick, -stealth] (v4) to [out=60,in=120] (v8);
	\draw[thick, -stealth] (v5) to [out=60,in=120] (v9);
    \node[rotate=20] at (0.6,-0.5){\textcolor{red}{\tiny$s_1\!-\!1$}};
	\node[rotate=20] at (1.6,-0.5){\textcolor{red}{\tiny$s_2\!-\!1$}};
	\node[rotate=20] at (2.6,-0.5){\textcolor{red}{\tiny$s_3\!-\!1$}};
	\node[rotate=20] at (3.6,-0.5){\textcolor{red}{\tiny$s_4\!-\!1$}};
	\node at (5,-0.5){\textcolor{red}{\footnotesize$1$}};
	\node at (6,-0.5){\textcolor{red}{\footnotesize$1$}};
	\node at (7,-0.5){\textcolor{red}{\footnotesize$1$}};
	\node at (8,-0.5){\textcolor{red}{\footnotesize$1$}};
	\node at (9,-0.5){\textcolor{red}{\footnotesize$1$}};
\end{scope}
\end{tikzpicture}
    \caption{Left:  The graph $G(4,13)$ with netflow $(1^4,0^9,-4)$ is an example of the type of flow polytope used in the proof of Theorem~\ref{thm.log-concavity}. Right: The  integer flows (mixed volumes) in the Lidskii formula for the volume of such polytope on the graph $G(4,13)$ and its subgraph $G(4,9)$ correspond to one of the interpretations of $E_{(s_1,s_2,s_3,s_4)}$ after reversing the graph and considering the subgraph of vertices with nonzero flow (compare with Figure~\ref{fig.flow_interpretations E_s}).}
    \label{fig:proof log concavity E_s}
\end{figure}
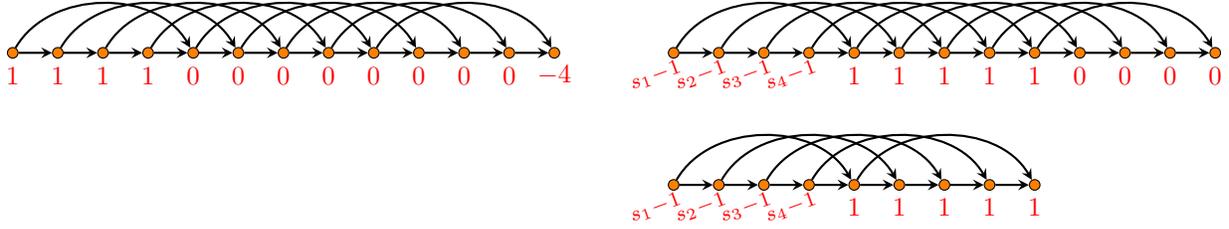

\begin{proof}
We show that $E_{\ss}$ equals $K_G(\ss-\tt)$ on the right side of Equation~\eqref{eq: lidskii outdegree} for a certain graph $G$ and net flow $\aa$.

The flow polytope of the graph $G(k,n-k+2)$ has dimension $d=n-2k+2$. 
Applying Lemma~\ref{thm.genlidskii} to $\F_{G(k,d+1)}(1^k,0^{d-k},-k)$, we have that $\tt=(1^{d-k+1}, 0^{k-1})$ and
\begin{equation}\label{eq.lidskii k-Springer}
\vol \F_{G(k,d+1)}(1^k,0^{d-k},-k) \,=\, \sum_{\ss} \binom{N}{\ss} \cdot 1 \cdot K_{G(k,d+1)}(\ss-\tt),
\end{equation}
where the zeros in the net flow $\aa = (1^k,0^{d-k},-k)$ restrict the sum to be over compositions $\ss=(s_1,\ldots,s_k)$ of $N$. 
Thus the Kostant partition functions appearing on the RHS are
\[
K_{G(k,d+1)}(\ss - \tt) = K_{G(k,d+1)}(s_1-1,\ldots,s_k-1,(-1)^{d-2k+1},0^k).
\]
 By the same argument as in the proof of Corollary~\ref{cor. k-Euler is kpf}, since the net flow on the last $k$ vertices of $G(k,d+1)$ is zero for the integer flows counted on the right hand side we have that 
\[
K_{G(k,d+1)}(s_1-1,\ldots,s_k-1,(-1)^{d-2k+1},0^k)
 = K_{G(k,d-k+1)}(s_1-1,\ldots,s_k-1,(-1)^{d-2k+1})
\]
A visualization of this operation is given in Figure~\ref{fig:proof log concavity E_s}. Since reversing the direction of the edges of the graph $G(k,d-k+1)$ yields a graph isomorphic to $G(k,d-k+1)$, then by reversing the flow, the Kostant partition remains unchanged \cite[Corollary 2.4]{MM19}. We conclude that 
\[
K_{G(k,d-k+1)}(s_1-1,\ldots,s_k-1,(-1)^{d-2k+1}) = K_{G(k,d-k+1)}(1^{d-2k+1},1-s_k,\ldots,1-s_1).
\]
By Proposition~\ref{prop: other interpretation E_s} the number of integer flows on the right hand side above are counted by the $k$-Entringer numbers. We have shown that 
\begin{equation} \label{E_s in Lidskii formula}
K_{G(k,d+1)}(\ss -\tt) = E_{\bf s};
\end{equation}
the desired log-concavity for $E_{\ss}$ follows from Lemma~\ref{log-concavity-kpf}.
\end{proof}

\begin{figure}
\begin{center}

\tikzset{every picture/.style={line width=0.75pt}} 

\begin{tikzpicture}[x=0.75pt,y=0.75pt,yscale=-1,xscale=1]

\draw  [color={rgb, 255:red, 155; green, 155; blue, 155 }  ,draw opacity=0.4 ] (197.49,84.61) -- (212.32,109.75) -- (182.65,109.75) -- cycle ;
\draw  [color={rgb, 255:red, 155; green, 155; blue, 155 }  ,draw opacity=0.4 ] (226.68,84.89) -- (241.51,110.03) -- (211.85,110.03) -- cycle ;
\draw  [color={rgb, 255:red, 155; green, 155; blue, 155 }  ,draw opacity=0.4 ] (256.35,84.61) -- (271.18,109.75) -- (241.51,109.75) -- cycle ;
\draw  [color={rgb, 255:red, 155; green, 155; blue, 155 }  ,draw opacity=0.4 ] (241.51,110.03) -- (256.35,135.17) -- (226.68,135.17) -- cycle ;
\draw  [color={rgb, 255:red, 155; green, 155; blue, 155 }  ,draw opacity=0.4 ] (271.18,110.03) -- (286.02,135.17) -- (256.35,135.17) -- cycle ;
\draw  [color={rgb, 255:red, 155; green, 155; blue, 155 }  ,draw opacity=0.4 ] (211.85,110.03) -- (226.68,135.17) -- (197.01,135.17) -- cycle ;
\draw  [color={rgb, 255:red, 155; green, 155; blue, 155 }  ,draw opacity=0.4 ] (182.18,110.03) -- (197.01,135.17) -- (167.34,135.17) -- cycle ;
\draw  [color={rgb, 255:red, 155; green, 155; blue, 155 }  ,draw opacity=0.4 ] (286.02,135.17) -- (300.85,160.32) -- (271.18,160.32) -- cycle ;
\draw  [color={rgb, 255:red, 155; green, 155; blue, 155 }  ,draw opacity=0.4 ] (256.35,135.17) -- (271.18,160.32) -- (241.51,160.32) -- cycle ;
\draw  [color={rgb, 255:red, 155; green, 155; blue, 155 }  ,draw opacity=0.4 ] (226.68,135.17) -- (241.51,160.32) -- (211.85,160.32) -- cycle ;
\draw  [color={rgb, 255:red, 155; green, 155; blue, 155 }  ,draw opacity=0.4 ] (197.01,135.17) -- (211.85,160.32) -- (182.18,160.32) -- cycle ;
\draw  [color={rgb, 255:red, 155; green, 155; blue, 155 }  ,draw opacity=0.4 ] (286.02,84.89) -- (300.85,110.03) -- (271.18,110.03) -- cycle ;
\draw  [color={rgb, 255:red, 155; green, 155; blue, 155 }  ,draw opacity=0.4 ] (211.85,160.32) -- (226.68,185.46) -- (197.01,185.46) -- cycle ;
\draw  [color={rgb, 255:red, 155; green, 155; blue, 155 }  ,draw opacity=0.4 ] (241.51,160.32) -- (256.35,185.46) -- (226.68,185.46) -- cycle ;
\draw  [color={rgb, 255:red, 155; green, 155; blue, 155 }  ,draw opacity=0.4 ] (271.18,160.32) -- (286.02,185.46) -- (256.35,185.46) -- cycle ;
\draw  [color={rgb, 255:red, 155; green, 155; blue, 155 }  ,draw opacity=0.4 ] (315.69,84.89) -- (330.52,110.03) -- (300.85,110.03) -- cycle ;
\draw  [color={rgb, 255:red, 155; green, 155; blue, 155 }  ,draw opacity=0.4 ] (300.37,110.31) -- (315.21,135.45) -- (285.54,135.45) -- cycle ;
\draw  [color={rgb, 255:red, 126; green, 211; blue, 33 }  ,draw opacity=0.4 ][line width=1.5]  (197.01,84.34) -- (315.69,84.89) -- (359.71,160.59) -- (315.21,236.57) -- (196.53,236.57) -- (152.03,160.04) -- cycle ;
\draw  [color={rgb, 255:red, 155; green, 155; blue, 155 }  ,draw opacity=0.4 ] (330.04,110.31) -- (344.88,135.45) -- (315.21,135.45) -- cycle ;
\draw  [color={rgb, 255:red, 155; green, 155; blue, 155 }  ,draw opacity=0.4 ] (315.69,135.17) -- (330.52,160.32) -- (300.85,160.32) -- cycle ;
\draw  [color={rgb, 255:red, 155; green, 155; blue, 155 }  ,draw opacity=0.4 ] (344.88,135.45) -- (359.71,160.59) -- (330.04,160.59) -- cycle ;
\draw  [color={rgb, 255:red, 155; green, 155; blue, 155 }  ,draw opacity=0.4 ] (300.37,161.14) -- (315.21,186.28) -- (285.54,186.28) -- cycle ;
\draw  [color={rgb, 255:red, 155; green, 155; blue, 155 }  ,draw opacity=0.4 ] (330.04,161.14) -- (344.88,186.28) -- (315.21,186.28) -- cycle ;
\draw  [color={rgb, 255:red, 155; green, 155; blue, 155 }  ,draw opacity=0.4 ] (315.21,186.28) -- (330.04,211.43) -- (300.37,211.43) -- cycle ;
\draw  [color={rgb, 255:red, 155; green, 155; blue, 155 }  ,draw opacity=0.4 ] (300.37,211.43) -- (315.21,236.57) -- (285.54,236.57) -- cycle ;
\draw  [color={rgb, 255:red, 155; green, 155; blue, 155 }  ,draw opacity=0.4 ] (286.02,185.46) -- (300.85,210.6) -- (271.18,210.6) -- cycle ;
\draw  [color={rgb, 255:red, 155; green, 155; blue, 155 }  ,draw opacity=0.4 ] (256.35,185.46) -- (271.18,210.6) -- (241.51,210.6) -- cycle ;
\draw  [color={rgb, 255:red, 155; green, 155; blue, 155 }  ,draw opacity=0.4 ] (226.68,185.46) -- (241.51,210.6) -- (211.85,210.6) -- cycle ;
\draw  [color={rgb, 255:red, 155; green, 155; blue, 155 }  ,draw opacity=0.4 ] (270.71,211.43) -- (285.54,236.57) -- (255.87,236.57) -- cycle ;
\draw  [color={rgb, 255:red, 155; green, 155; blue, 155 }  ,draw opacity=0.4 ] (241.04,211.43) -- (255.87,236.57) -- (226.2,236.57) -- cycle ;
\draw  [color={rgb, 255:red, 155; green, 155; blue, 155 }  ,draw opacity=0.4 ] (197.01,185.46) -- (211.85,210.6) -- (182.18,210.6) -- cycle ;
\draw  [color={rgb, 255:red, 155; green, 155; blue, 155 }  ,draw opacity=0.4 ] (211.37,211.43) -- (226.2,236.57) -- (196.53,236.57) -- cycle ;
\draw  [color={rgb, 255:red, 155; green, 155; blue, 155 }  ,draw opacity=0.4 ] (182.18,160.32) -- (197.01,185.46) -- (167.34,185.46) -- cycle ;
\draw  [color={rgb, 255:red, 155; green, 155; blue, 155 }  ,draw opacity=0.4 ] (166.87,134.9) -- (181.7,160.04) -- (152.03,160.04) -- cycle ;
\draw [color={rgb, 255:red, 144; green, 19; blue, 254 }  ,draw opacity=0.8 ][fill={rgb, 255:red, 245; green, 166; blue, 35 }  ,fill opacity=1 ][line width=1.5]    (226.2,236.57) -- (329,62.33) ;
\draw [shift={(330.52,59.74)}, rotate = 480.54] [color={rgb, 255:red, 144; green, 19; blue, 254 }  ,draw opacity=0.8 ][line width=1.5]    (14.21,-4.28) .. controls (9.04,-1.82) and (4.3,-0.39) .. (0,0) .. controls (4.3,0.39) and (9.04,1.82) .. (14.21,4.28)   ;
\draw [color={rgb, 255:red, 65; green, 117; blue, 5 }  ,draw opacity=0.8 ][fill={rgb, 255:red, 245; green, 166; blue, 35 }  ,fill opacity=1 ][line width=1.5]    (359.71,160.59) -- (125.36,160.05) ;
\draw [shift={(122.36,160.04)}, rotate = 360.13] [color={rgb, 255:red, 65; green, 117; blue, 5 }  ,draw opacity=0.8 ][line width=1.5]    (14.21,-4.28) .. controls (9.04,-1.82) and (4.3,-0.39) .. (0,0) .. controls (4.3,0.39) and (9.04,1.82) .. (14.21,4.28)   ;
\draw [color={rgb, 255:red, 74; green, 144; blue, 226 }  ,draw opacity=0.8 ][fill={rgb, 255:red, 245; green, 166; blue, 35 }  ,fill opacity=1 ][line width=1.5]    (315.21,236.57) -- (213.36,62.33) ;
\draw [shift={(211.85,59.74)}, rotate = 419.69] [color={rgb, 255:red, 74; green, 144; blue, 226 }  ,draw opacity=0.8 ][line width=1.5]    (14.21,-4.28) .. controls (9.04,-1.82) and (4.3,-0.39) .. (0,0) .. controls (4.3,0.39) and (9.04,1.82) .. (14.21,4.28)   ;
\draw [color={rgb, 255:red, 144; green, 19; blue, 254 }  ,draw opacity=0.8 ][fill={rgb, 255:red, 245; green, 166; blue, 35 }  ,fill opacity=1 ][line width=1.5]    (403,133) -- (415.86,110.73) ;
\draw [shift={(417.36,108.13)}, rotate = 480] [color={rgb, 255:red, 144; green, 19; blue, 254 }  ,draw opacity=0.8 ][line width=1.5]    (14.21,-4.28) .. controls (9.04,-1.82) and (4.3,-0.39) .. (0,0) .. controls (4.3,0.39) and (9.04,1.82) .. (14.21,4.28)   ;
\draw [color={rgb, 255:red, 74; green, 144; blue, 226 }  ,draw opacity=0.8 ][fill={rgb, 255:red, 245; green, 166; blue, 35 }  ,fill opacity=1 ][line width=1.5]    (417.83,173.14) -- (404.52,150.58) ;
\draw [shift={(403,148)}, rotate = 419.46000000000004] [color={rgb, 255:red, 74; green, 144; blue, 226 }  ,draw opacity=0.8 ][line width=1.5]    (14.21,-4.28) .. controls (9.04,-1.82) and (4.3,-0.39) .. (0,0) .. controls (4.3,0.39) and (9.04,1.82) .. (14.21,4.28)   ;
\draw [color={rgb, 255:red, 65; green, 117; blue, 5 }  ,draw opacity=0.8 ][fill={rgb, 255:red, 245; green, 166; blue, 35 }  ,fill opacity=1 ][line width=1.5]    (422,199) -- (395.33,199) ;
\draw [shift={(392.33,199)}, rotate = 360] [color={rgb, 255:red, 65; green, 117; blue, 5 }  ,draw opacity=0.8 ][line width=1.5]    (14.21,-4.28) .. controls (9.04,-1.82) and (4.3,-0.39) .. (0,0) .. controls (4.3,0.39) and (9.04,1.82) .. (14.21,4.28)   ;

\draw (179,204.4) node [anchor=north west][inner sep=0.75pt]  [font=\scriptsize]  {$0$};
\draw (208,204.4) node [anchor=north west][inner sep=0.75pt]  [font=\scriptsize]  {$0$};
\draw (237,204.4) node [anchor=north west][inner sep=0.75pt]  [font=\scriptsize]  {$0$};
\draw (267,204.4) node [anchor=north west][inner sep=0.75pt]  [font=\scriptsize]  {$0$};
\draw (297,204.4) node [anchor=north west][inner sep=0.75pt]  [font=\scriptsize]  {$0$};
\draw (326,204.4) node [anchor=north west][inner sep=0.75pt]  [font=\scriptsize]  {$0$};
\draw (311,180.4) node [anchor=north west][inner sep=0.75pt]  [font=\scriptsize]  {$0$};
\draw (194,179.4) node [anchor=north west][inner sep=0.75pt]  [font=\scriptsize]  {$2$};
\draw (222,179.4) node [anchor=north west][inner sep=0.75pt]  [font=\scriptsize]  {$2$};
\draw (252,179.4) node [anchor=north west][inner sep=0.75pt]  [font=\scriptsize]  {$2$};
\draw (282,179.4) node [anchor=north west][inner sep=0.75pt]  [font=\scriptsize]  {$1$};
\draw (208,154.4) node [anchor=north west][inner sep=0.75pt]  [font=\scriptsize]  {$4$};
\draw (223,129.4) node [anchor=north west][inner sep=0.75pt]  [font=\scriptsize]  {$5$};
\draw (238,154.4) node [anchor=north west][inner sep=0.75pt]  [font=\scriptsize]  {$4$};
\draw (267,155.4) node [anchor=north west][inner sep=0.75pt]  [font=\scriptsize]  {$3$};
\draw (296,155.4) node [anchor=north west][inner sep=0.75pt]  [font=\scriptsize]  {$1$};
\draw (252,130.4) node [anchor=north west][inner sep=0.75pt]  [font=\scriptsize]  {$5$};
\draw (281,130.4) node [anchor=north west][inner sep=0.75pt]  [font=\scriptsize]  {$3$};
\draw (238,104.4) node [anchor=north west][inner sep=0.75pt]  [font=\scriptsize]  {$5$};
\draw (267,104.4) node [anchor=north west][inner sep=0.75pt]  [font=\scriptsize]  {$5$};
\draw (252,79.4) node [anchor=north west][inner sep=0.75pt]  [font=\scriptsize]  {$5$};
\draw (333,43.4) node [anchor=north west][inner sep=0.75pt]  [font=\scriptsize,color={rgb, 255:red, 144; green, 19; blue, 254 }  ,opacity=1 ]  {$\ee_{1} -\ee_{2}$};
\draw (119,140.4) node [anchor=north west][inner sep=0.75pt]  [font=\scriptsize,color={rgb, 255:red, 65; green, 117; blue, 5 }  ,opacity=1 ]  {$\ee_{2} -\ee_{3}$};
\draw (213,43.4) node [anchor=north west][inner sep=0.75pt]  [font=\scriptsize,color={rgb, 255:red, 74; green, 144; blue, 226 }  ,opacity=1 ]  {$\ee_{1} -\ee_{3}$};
\draw (424,112.4) node [anchor=north west][inner sep=0.75pt]  [font=\footnotesize]  {$3^2=E_{( 2,1,2)}^{2} \geq E_{( 1,2,2)} E_{( 3,0,2)} =2\cdot 3$};
\draw (424,151.4) node [anchor=north west][inner sep=0.75pt]  [font=\footnotesize]  {$3^2=E_{( 2,1,2)}^{2} \geq E_{( 3,1,1)} E_{( 1,1,3)} =5\cdot 1$};
\draw (424,191.4) node [anchor=north west][inner sep=0.75pt]  [font=\footnotesize]  {$3^2=E_{( 2,1,2)}^{2} \geq E_{( 2,2,1)} E_{( 2,0,3)} =4\cdot 1$};

\end{tikzpicture}
\end{center}
\caption{Example of the log-concavity of the $k$-Entringer numbers along root directions.}
\label{fig: log concavity k-Entringer}
\end{figure}
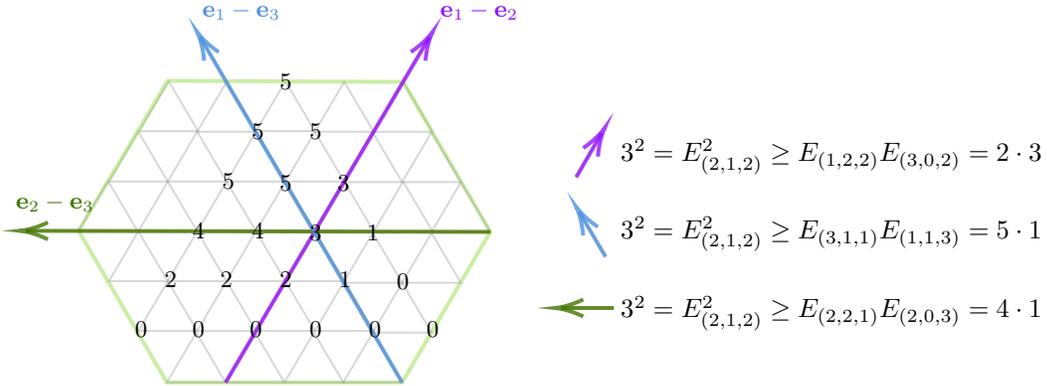

\begin{example}
Figure~\ref{fig: log concavity k-Entringer} highlights the log-concavity along root directions for the $3$-Entringer numbers that are indexed by compositions of $5$. For each root direction $\ee_1-\ee_2$, $\ee_1-\ee_3$, and $\ee_2-\ee_3$, the two $3$-Entringer numbers on either side of $E_{(2,1,2)}$ in that direction multiply to no more than $E_{(2,1,2)}^2$.
\end{example}

\subsection{\texorpdfstring{$k$}{}-Springer numbers} 
In the proof of Theorem~\ref{thm.log-concavity} we studied the volume of the flow polytope of $G(k,d+1)$ with netflow $(1^k,0^{d-k},-k)$. Moreover, in  \cite[Theorem 6.11]{BGHHKMY}, Benedetti et al.\ showed that in the case $k=2$ this volume  equals the  \emph{Springer number}  $S_{n-1}$ \cite[A001586]{oeis}. This motivates the definition of a generalization of Springer numbers that are a multinomial transform of the numbers $E_{\ss}$.

\begin{definition}[$k$-Springer numbers]
For $k \geq 1$ and $N=d-k+1\geq0$, let 
\[
S_{k,d} = \vol \mathcal{F}_{G(k,d+1    )}(1^k,0^{d-k},-k).
\]
\end{definition}

See Table~\ref{table:kspringer} for some values of this sequence.

\begin{proposition} \label{prop.equation kSpringer}
Let $k\geq 1$ and $N=d-k+1\geq0$. We have that
\[
S_{k,d} = \sum_{\ss \in T^k_N} \binom{d-k+1}{\ss} E_{\ss}.
\]
\end{proposition}

\begin{proof}
By Equation~\eqref{eq.lidskii k-Springer} and the definition of $S_{k,d}$, we have
\[
S_{k,d} = \sum_{\ss \in T^k_N} \binom{d-k+1}{\ss} \cdot K_{G(k,d+1)}(\ss -\tt).
\]
The result then follows by applying  Equation~\eqref{E_s in Lidskii formula} to the LHS above.
\end{proof}

\begin{table}[t!]
$$\begin{array}{crrrrrrrrrrrr}
\hline
k\backslash d& 1& 2& 3& 4& 5& 6& 7& 8&9&10 & \\ \hline
1&1&2&6&24&120&720&5040&40320&362880& 3628800 &\href{https://oeis.org/A000142}{\textup{A000142}}\\
2&	 1& 1& 3& 11& 57& 361& 2763 &24611& 250737&2873041 & \href{https://oeis.org/A001586}{\textup{A001586}}\\
3&1 &1 & 1& 3& 16& 88& 625& 5527& 55760&640540\\
4&1 &1 & 1& 1& 3& 16& 125& 927&8357&91735\\ \hline
\end{array}
$$
\caption{Table of the $k$-Springer numbers $S_{k,d}$.}
\label{table:kspringer}
\end{table}


\subsection{A Lorentzian polynomial related to \texorpdfstring{$k$}{}-Entringer numbers} \label{subsec.Lorentzian}

In Theorem~\ref{thm.log-concavity} we showed that the $k$-Entringer numbers satisfy certain log-concavity relations. This is part of a more general story: Br\"and\'en and Huh \cite{BH} defined {\em Lorentzian polynomials} as a generalization of volume polynomials in algebraic geometry and {\em stable polynomials} in optimization. The latter are  a generalization  of real-rooted polynomials in a multivarite setting introduced by Borcea and Br\"and\'en in  \cite{BB09}. 
The following result is implicit in \cite[Section 2]{MS} and it follows from the fact that volume polynomials of Minkowski sums of convex bodies are Lorentzian \cite[Theorem 9.1]{BH} and that the flow polytope $\mathcal{F}_G(\aa)$ is a Minkowski sum of flow poytopes \cite[Section 3.4]{BV}.

\begin{lemma} \label{lem.volume flow is Lorentzian}
Let $G$ be a directed graph on the vertex set $[n+1]$ with $m$ edges, such that the out-degree of each vertex in $[n]$ is at least one. The polynomial $\vol \mathcal{F}_G({\bf a})$ from Equation~\eqref{eq: lidskii outdegree} is Lorentzian in the variables $a_1,\ldots,a_n$.
\end{lemma}

For $k\geq 1$ and $N=d-k+1\geq 0$, let $A_{k,d}(x_1,\ldots,x_k)$ be the polynomial
\[
A_{k,d}(x_1,\ldots,x_k) := \sum_{\ss \in T^k_N}  E_{\ss} \cdot  x_1^{s_1}\cdots x_k^{s_k}. 
\]
The \emph{normalization operator} $\mathsf{N}$ on $\mathbb{R}[x_1,\ldots,x_k]$ acts by sending $x_1^{s_1}\cdots x_k^{s_k} \mapsto \frac{x_1^{s_1}}{s_1!}\cdots \frac{x_k^{s_k}}{s_k!}$. We collect some straightforward identities of the polynomials $A_{k,d}({\bf x})$.

\begin{proposition} \label{prop.easy.props.multivariable.Adk}
For $k \geq 1$ and $d-k+1\geq0$, we have that
\begin{itemize}
    \item[(a)] $A_{k,d}(1,\ldots,1) = A_{k,d}$,
    \item[(b)] $(d-k+1)!\cdot \mathsf{N}(A_{k,d}(1,\ldots,1)) = S_{k,d}$,
    \item[(c)] for nonnegative integers $x_1,\ldots,x_k$
    \[
    (d-k+1)!\cdot \mathsf{N}(A_{k,d}(x_1,\ldots,x_n)) = \vol \mathcal{F}_{G(k,d+1)}(x_1,\ldots,x_k,0^{d-k},-{\textstyle \sum_{i=1}^k x_i}).
    \]
\end{itemize}
\end{proposition}

\begin{proof}
Part (a) follows from the definition of $A_{k,d}({\bf x})$ and Equation~\eqref{eqn.genEnt}. Part (b) follows from  Proposition~\ref{prop.equation kSpringer}. Lastly, part (c) follows from the Lidskii formula Theorem~\ref{thm.genlidskii} and Equation~\eqref{E_s in Lidskii formula}.
\end{proof}

\begin{corollary}
Let $k\geq 1$ and $d-k+1\geq 0$. The quantity $\mathsf{N}(A_{k,d}({\bf x}))$ is Lorentzian.
\end{corollary}

\begin{proof}
This follows from Lemma~\ref{lem.volume flow is Lorentzian} and Proposition~\ref{prop.easy.props.multivariable.Adk} (c).
\end{proof}


\section{The \texorpdfstring{$h^*$}{}-polynomial of flow polytopes}\label{sec:h_star}

In this article we have studied the volume of flow polytopes $\mathcal{F}_G$. A popular refinement of the volume of integral polytope is known as the \emph{$h^*$-polynomial} of an integral polytope. We review these definitions and refer to \cite[Chapters 3 \& 10]{BR} for more background on this subject. 

\begin{definition} \label{def: ehrhart_poly_and_h}
For an integral polytope $\P \subset \mathbb{R}^n$ and a nonnegative integer $t\geq 1$ let $L_{\P}(t):=\# t\P \cap \mathbb{Z}^n$, where $t\P$ is the $t$-th dilation of $\P$.
It is known that for an integral $d$-polytope $\P$ the function $L_{\P}(t)$ is a polynomial of degree $d$ called the \emph{Ehrhart polynomial} of $\P$.
For an integral $d$-polytope $\P$, the \emph{$h^*$-polynomial} of $\P$ is the polynomial $h^*_{\P}(z)$  defined by
 \[
 1 + \sum_{t \geq 1} L_{\P}(t)z^t = \frac{h^*_{\P}(z)}{(1-z)^{d+1}}.
 \]
\end{definition}

Stanley \cite{St_hstar} showed that the coefficients of $h^*_{\P}(t)$ are nonnegative integers. Moreover, from the definition of $h^*_{\P}(z)$ and the fact that the leading term of  $L_{\P}(t)$ is $\vol(\P)/d!$ (see, for example, \cite[Lemma~3.19]{BR}) one can see that $h^*_{\P}(1)=\vol(\P)$. Thus the $h^*$-polynomial gives a refinement of the volume of $\P$. An interesting open problem is to find a combinatorial rule to compute the $h^*$-vector of flow polytopes $\mathcal{F}_G$ for any graph $G$. Special cases have been solved in \cite{Khstar, MMS, BGMY}. In \cite{AyyerJosuatvergesRamassamy2018} the authors found a rule to compute the $h^*$-vectors of the consecutive coordinates polytopes $\mathcal{B}_{\C}$. Since by Proposition~\ref{proposition:intervals_are_row_convex} these polytopes are flow polytopes  corresponding to certain non-nested graphs then we have a rule for the $h^*$-vectors of this Catalan family of flow polytopes.

Given a total cyclic order $\gamma$ of $\{0,1,\dots, n\}$, let $\pi(\gamma)$ be the permutation of $[n]$ obtained by reading the elements following zero in clockwise order. Given a word $\pi=\pi_1\pi_2\cdots \pi_n$, we say $\pi$ has a a {\em descent} at position $i$ if  $\pi_i>\pi_{i+1}$. The number of descents of $\pi$ is denoted by $\des(\pi)$.

\begin{example}
For the total cyclic order $\gamma$ depicted in Figure \ref{fig:G_and_gamma} we have that $\pi(\gamma)=125364$ with descents at the third and fifth positions so that $\des(\pi(\gamma))=2$.
\end{example}

\begin{theorem}[{Ayyer--Josuat-Verg\`es--Ramassamy  \cite{AyyerJosuatvergesRamassamy2018}}] \label{thm.hstar B_C}
For a collection $\C$ of intervals in $[d]$, the $h^*$-polynomial of the polytope $\mathcal{B}_{\C}$ is
\[
h^*_{\mathcal{B}_{\C}}(z) = \sum_{ \gamma \in \A_{\C}} z^{\des(\pi(\gamma))}.
\]
\end{theorem}

Recall that for a non-nested graph $G$, the set of upper and  lower $G$-cyclic orders are equal by Proposition~\ref{prop.ofmaintheorem}, and are denoted by  $\mathcal{A}_G$.
 
\begin{theorem}\label{theorem:hstar_non_nested}
For a non-nested graph $G$, the $h^*$-polynomial of the flow polytope $\mathcal{F}_{G}$ is
\[
h^*_{\mathcal{F}_{G}}(z) = \sum_{ \gamma \in \mathcal{A}_G} z^{\des(\pi(\gamma))}.
\]
\end{theorem}

\begin{proof}
Given a non-nested graph $G$ by Proposition~\ref{proposition.non nested are consecutive} there exists a collection $\C$ of intervals in $[d]$ such that $\mathcal{F}_G$ is integrally equivalent to $\mathcal{B}_{\C}$. Thus $h^*_{\mathcal{F}_G}(z)=h^*_{\mathcal{B}_{\C}}(z)$. The result then follows by Theorem~\ref{thm.hstar B_C} since in this case $\mathcal{A}_{\C}=\mathcal{A}_G$. 
\end{proof}

\begin{example}
For the collection $\C=\{[1,2],[2,4],[3,6],[5,7]\} \subset [7]$ and its associated graph in Figure~\ref{fig.BS}, the $h^*$-polynomial is 
\[
h^*_{\mathcal{F}_{G}}(z) = h^*_{\mathcal{B}_{\C}}(z)  = 1 + 12z + 25z^2 + 10z^3.
\]
\end{example}

Given a statistic $\stat:\C\rightarrow \RR$ on a set $\C$ we denote by $P_{\C,\stat}(z)=\sum_{x\in \C}z^{\stat(x)}$ the \emph{generating polynomial}  of $\stat$ on $\C$.
In general, when a graph $G$ contains nested pairs of edges,  $\mathcal{A}^{\uparrow}_G$ and $\mathcal{A}^{\downarrow}_G$ are not equal and the descent statistic gives different polynomials $P_{\mathcal{A}^{\uparrow}_G,\des}(z)\neq P_{\mathcal{A}^{\downarrow}_G,\des}(z)$. (See Example~\ref{ex: counterex_des_h*} below.) While Theorem~\ref{theorem:hstar_non_nested} cannot hold in this case, computations of the $h^*$-polynomial of flow polytopes using the Lidskii formula for lattice points (\cite[Theorem 38]{BV} and \cite[Equation (1.2)]{MM19}) suggest that this polynomial lies ``in between" $P_{\mathcal{A}^{\uparrow}_G,\des}(z)$ and $P_{\mathcal{A}^{\downarrow}_G,\des}(z)$. 

For two polynomials $a(z)=\sum_{i\ge0} a_i z^i$ and $b(z)=\sum_{i\ge0} b_i z^i$, we say that $a(z)$ \emph {is dominated by} $b(z)$ and write $a(z) \lhd b(z)$ if $\sum_{i=0}^k a_i \leq \sum_{i=0}^k b_i$ for every $k\ge 0$.

\hstarconjecture

This conjecture has been computationally verified for all simple spinal graphs $G$ with $\leq 7$ vertices.

\begin{example} \label{ex: counterex_des_h*}
Let $G$ be the graph from 
Figure~\ref{fig:G_and_gamma}. The polynomials
\[
P_{\mathcal{A}^\downarrow_G,\des}(z) = h_{\F_G}^*(z) = 1 + 7z + 7z^2 + z^3 \textup{ and } P_{\mathcal{A}^{\uparrow}_G,\des}(z) = 1 + 9z + 6z^2
\]
satisfy Conjecture \ref{conj: hstar}.
\end{example}

\begin{example}
Let $G$ be the complete graph $k_7$ on seven vertices. The polynomials
\begin{align*}
    P_{\mathcal{A}^\downarrow_{k_7},\des}(z) &= 1 + 15z + 55z^2 + 59z^3 + 10z^4,\\
    h^*_{\mathcal{F}_{k_7}}(z) &= 1 + 16z + 58z^2 + 56z^3 + 9z^4, \textup{ and}\\
    P_{\mathcal{A}^\uparrow_{k_7},\des}(z) &= 1 + 18z + 64z^2 + 51z^3 + 6z^4
\end{align*}
satisfy Conjecture \ref{conj: hstar}.

\end{example}


\section{Open questions and further work} \label{sec:further_work}

\subsection{Graphs with equivalent flow polytopes}

In Section~\ref{sec:integralequivalence} we showed several operations on graphs that yield integrally equivalent flow polytopes $\mathcal{F}_G$. We would like a full characterization of graphs with integrally equivalent flow polytopes. See Question~\ref{question:equivalent}.

\subsection{Structure and enumeration of non-redundant column-convex \texorpdfstring{$\{0,1\}$}{} matrices}

In Section~\ref{sec.counting collections} we showed that there are Catalan many non-redundant collections of intervals of $[d]$ by putting them in bijection with antichains in the type $A$ root poset. The latter have an inherent poset structure (ordered by inclusion) and we would like to know whether the flow polytopes corresponding to each antichain are compatible with this poset. See Question~\ref{q.inherent catalan structure}.

\subsection{Generating functions for \texorpdfstring{$k$}{}-Euler and \texorpdfstring{$k$}{}-Entringer numbers}

The Euler numbers $E_n$, Springer numbers $S_n$ and Entringer numbers $E_{n,r}$ have beautiful generating functions. The following are from \cite[Proposition 1.6.1]{Stanley2012}, \cite[\href{https://oeis.org/A001586}{\textup{A001586}}]{oeis}, and \cite[Exercise 141]{Stanley2012}.
\[
\sum_{n\geq 0} E_n \frac{x^n}{n!} = \sec x+\tan x \qquad \sum_{n\geq 0} S_n \frac{x^n}{n!} = \frac{1}{\cos x - \sin x}
\]
\[
\sum_{m,n \geq 0} E_{m+n,[m,n]} \frac{x^m}{m!} \frac{y^n}{n!} = \frac{\cos x + \sin x}{\cos(x+y)}, \quad\text{where}  \quad [m,n] = \begin{cases} m &\text{ if } m+n \text{ is odd},\\ n &\text{ if } m+n \text{ is even}\end{cases}
\]
We propose generalizing these results for the $k$-Euler numbers $A_{k,d}$, the $k$-Springer numbers $S_{k,d}$, and the $k$-Entringer numbers $E_{\ss}$.

\begin{question}
Is there a closed form for the generating function for  $k$-Euler, $k$-Springer, and $k$-Entringer numbers?
\end{question}

\subsection{A general rule to compute \texorpdfstring{$h^*$}{}-polynomial of flow polytopes}

Theorem~\ref{theorem:hstar_non_nested} gave a rule to compute the $h^*$-polynomial of flow polytopes $\mathcal{F}_G$ of non-nested graphs $G$ using the rule to compute the $h^*$-polynomials of consecutive coordinate polytopes in \cite{AyyerJosuatvergesRamassamy2018}. We would like to find a formula to compute the $h^*$-polynomial for all flow polytopes $\mathcal{F}_G$ and $\mathcal{F}_G({\bf a})$.  Here is a list of special cases that are known.

\begin{itemize}
    \item[(i)] For planar graphs $G$, it is shown in \cite{MMS} that the polytope $\mathcal{F}_G$ is integrally equivalent to an {\em order polytope} of a poset \cite{Stanley86}. Order polytopes of posets have a rule for their $h^*$-polynomial by descents of linear extensions due to Stanley \cite[Theorem 3.15.8]{Stanley2012}.
    \item[(ii)] A special case of the previous case is the order polytope of the zigzag poset, which is integrally equivalent to the flow polytope of the planar graph $G(2,d)$. Coons and Sullivant~\cite[Theorem 1.9]{CS} give a new combinatorial interpretation for the coefficients of the $h^*$-polynomial via shellings of the canonical triangulation of the order polytope.
    \item[(iii)] For non-nested graphs $G$, which includes the distance graphs $G(k,d+k)$, the $h^*$-polynomial of $\mathcal{F}_G$ is given by descents on $G$-cyclic orders (Theorem~\ref{theorem:hstar_non_nested}).
    \item[(iv)] For {\em $\nu$-caracol graphs}, which are certain graphs indexed by lattice paths, the authors in \cite{BGMY}  give a combinatorial formula involving $\nu$-Narayana numbers for computing the $h^*$-polynomial in two ways: via shellings of a generalization of the {\em Tamari lattice}, and of principal order ideals in {\em Young's lattice}.
\end{itemize}

Still, these cases do not cover some important flow polytopes like the {\em Chan--Robbins--Yuen polytope}  $\mathcal{F}_{k_{n+1}}$ \cite{CRY}, the {\em Tesler polytope} $\mathcal{F}_{k_{n+1}}({\bf 1})$ \cite{MMR}, the flow polytope of the {\em caracol graph}  with netflow ${\bf 1}$ \cite{BGHHKMY} or the {\em Pitman-Stanley polytope} \cite{PS}. Some positivity  results on the $h^*$-polynomial and Ehrhart series of flow polytopes first appeared in \cite{Khstar}. One approach to compute $h^*$-polynomials is by studying triangulations of the polytopes.  Danilov--Karzanov--Koshevoy in \cite{DKK} studied  regular unimodular trinagulatons of flow polytopes $\mathcal{F}_G$ and M\'eszaros--Morales and Kapoor--M\'esz\'aros--Setiabrata \cite{MM19,KMS19} studied subdivisions of flow polytopes $\mathcal{F}_G({\bf b})$ into products of simplices (which can each be further triangulated) related to the Lidskii formulas (see Theorem~\ref{thm.genlidskii}).

\section*{Acknowledgements}
This work was supported by the
American Institute of Mathematics through their SQuaRE program. We are very appreciative
of their support and funding which made this research collaboration possible. We thank Carolina Benedetti and Pamela E. Harris for fruitful discussions throughout the process. We also thank Arvind Ayyer, Matthieu Josuat-Verg\`es, and Sanjay Ramassamy for telling us about \cite{AyyerJosuatvergesRamassamy2018}. C.\ R.\ H.\ Hanusa was partially supported by a PSC-CUNY Award, jointly funded by The Professional Staff Congress and The City University of New York. A.\ H.\ Morales was partially supported by NSF Grant DMS-1855536, M.\ Yip was partially supported by Simons Collaboration Grant 429920.

\bibliographystyle{plain}
\bibliography{h_star}

\end{document}